\documentclass[10pt,reqno,twoside]{amsart}

\usepackage{amsfonts}
\usepackage{paralist}
  \usepackage{graphics}
  \usepackage{epsfig}
 \usepackage[colorlinks=true]{hyperref}
\usepackage{hyperref}

\numberwithin{equation}{section}
\usepackage{latexsym}
\usepackage{amsmath}
\usepackage{amssymb}
\usepackage{mathrsfs}
\usepackage{graphicx}
\usepackage{ifthen}
\usepackage[T1]{fontenc}
\usepackage[latin1]{inputenc}
\usepackage{color}
\newtheorem{Satz}{Theorem}[section]

\newtheorem{Def}[Satz]{Definition}

\newcommand{\oo}{\overline \Omega}

\newcommand{\po}{\partial\Omega}

\newcommand{\dist}{\text {dist}}

\newcommand{\diag}{\text {diag}}

\newcommand{\diam}{\text {diam}}

\newcommand{\supp}{\text {supp}}

\begin{document}

\title[Concentrating   solutions]
{Concentrating solutions for an anisotropic planar elliptic  Neumann problem with Hardy-H\'{e}non weight and large exponent}

\author[Yibin Zhang]{Yibin Zhang}
\address{College of Sciences, Nanjing Agricultural University, Nanjing 210095, China}
\email{yibin10201029@njau.edu.cn}



\begin{abstract}
Let $\Omega$ be a bounded  domain in $\mathbb{R}^2$ with smooth boundary,
we study
the following anisotropic elliptic  Neumann problem
with Hardy-H\'{e}non weight
$$
\begin{cases}
-\nabla(a(x)\nabla u)+a(x)u=a(x)|x-q|^{2\alpha}u^p,\,\,\,\,
u>0\,\,\,\,\,
\textrm{in}\,\,\,\,\,
\Omega,\\[2mm]
\frac{\partial u}{\partial\nu}=0\,\,
\qquad\quad\qquad\qquad\qquad
\qquad\qquad\qquad\qquad
\,\ \ \,\,\,\,
\textrm{on}\,\,\,
\partial\Omega,
\end{cases}
$$
where $\nu$
denotes the outer unit normal vector to $\partial\Omega$, $q\in\overline{\Omega}$,
$\alpha\in(-1,+\infty)\setminus\mathbb{N}$, $p>1$ is a large exponent
and $a(x)$ is a  positive smooth function.
We investigate the effect of the interaction between anisotropic coefficient
$a(x)$ and   singular source $q$ on the existence of concentrating solutions.
We show that if
$q\in\Omega$ is a strict local maximum point of $a(x)$,
there exists a family of  positive solutions with arbitrarily many
interior spikes accumulating to
$q$; while if $q\in\partial\Omega$ is a strict local maximum point of $a(x)$
and satisfies $\langle\nabla a(q),\,\nu(q)\rangle=0$,
such a problem has a family of  positive solutions with arbitrarily many mixed
interior and boundary spikes accumulating to $q$. In particular, we find
that concentration at singular source $q$ is always possible whether $q\in\overline{\Omega}$
is an isolated local maximum point of $a(x)$ or not.
\end{abstract}

\maketitle

\section{Introduction}
This paper deals with the existence and profile of solutions
for the following anisotropic elliptic Neumann problem
\begin{equation}\label{1.1}
\left\{\aligned
&-\nabla(a(x)\nabla u)+a(x)u=a(x)|x-q|^{2\alpha}u^p,\,\,\,\,
u>0\,\,\,\,\,
\textrm{in}\,\,\,\,\,
\Omega,\\[1mm]
&\frac{\partial u}{\partial\nu}=0\,\,
\qquad\quad\qquad\qquad\qquad
\qquad\qquad\qquad\qquad
\,\ \ \,\,\,\,
\textrm{on}\,\,\,
\partial\Omega,
\endaligned\right.
\end{equation}
where $\Omega$ is a bounded domain in $\mathbb{R}^2$
with smooth boundary,
$\nu$
denotes the outer unit normal vector to $\po$,
$q\in\oo$,
$\alpha\in(-1,+\infty)\setminus\mathbb{N}$, $p>1$ is a large
 exponent
and $a(x)$ is a  positive smooth function over $\oo$.
The term $|\cdot|^{2\alpha}$ in equation (\ref{1.1}) is called
the Hardy weight if $-1<\alpha<0$, whereas the H\'{e}non weight  if $\alpha>0$
(see \cite{H1,H}).
We are interested in  solutions of problem (\ref{1.1})
which exhibit the {\it concentration phenomenon} as
the exponent $p$ tends to  infinity.

This work is strongly  motivated by
some extensive research involving   the case $\alpha=0$ in
equation (\ref{1.1}):
\begin{equation}\label{1.4}
\left\{\aligned
&-\nabla(a(x)\nabla u)+a(x)u=a(x)u^p,\,\,\,\,
u>0\,\,\,\,\,
\textrm{in}\,\,\,\,\,
\Omega,\\[1mm]
&\frac{\partial u}{\partial\nu}=0\,\,
\qquad\quad\qquad\qquad\qquad
\qquad\qquad
\,\ \ \,\,\,\,\,
\textrm{on}\,\,\,
\partial\Omega,
\endaligned\right.
\end{equation}
where $\Omega$ is  a smooth bounded domain in $\mathbb{R}^n$ with $n\geq2$ and $p>1$.
Equation (\ref{1.4}) has a strong biological meaning because it
arises from the study of
steady states for the logarithmic Keller-Segel system in
chemotaxis (see \cite{KS}):
\begin{equation}\label{1.5}
\aligned
\left\{\aligned
&C_1\Delta \psi-\chi\nabla\cdot(\psi\nabla \log \omega)=0\,\,\ \,\textrm{in}\quad\,\,\,\,\mathcal{D},\\[1mm]
&C_2\Delta \omega-a\omega+b\psi=0\,\,\,\,
\,\quad\,\quad\quad\,\textrm{in}\quad\,\,\,\,\mathcal{D},\\
&\frac{\partial \omega}{\partial \nu}=\frac{\partial \psi}{\partial \nu}=0\,\,\
\ \ \,\,\,\,
\,\quad\,\quad\,\quad\,\quad\,\,
\textrm{on}\quad\,\,\partial\mathcal{D},\\
&\frac1{|\mathcal{D}|}\int_{\mathcal{D}}\psi(x)dx=\bar{\psi}>0\,\,\,\ \,\,\quad\,\,\,\,\,\textrm{(prescribed),}
\endaligned\right.
\endaligned
\end{equation}
where $\mathcal{D}$ is a smooth  bounded domain in $\mathbb{R}^N(N\geq2)$
and the constants $C_1$, $C_2$, $a$, $b$ and $\chi$ are positive.
Testing the first equation in (\ref{1.5}) against
$\nabla(\log \psi-p\log\omega)$ with $p=\chi/C_1$, we find
$$
\aligned
\int_{\mathcal{D}}\psi|\nabla(\log \psi-p\log\omega)|^2=0,
\endaligned
$$
i.e. $\psi=\lambda\omega^p$ for some  constant $\lambda>0$. Furthermore, setting
$\varepsilon^2=C_2/a$, $\gamma=(b\lambda/a)^{\frac1{p-1}}$ and $\upsilon=\gamma\omega$,
we have that $\upsilon$ satisfies the Lin-Ni-Takagi problem in \cite{LNT,NT1,NT2}, namely
the  singularly perturbed elliptic  Neumann equation
\begin{equation}\label{1.6}
\aligned
\left\{\aligned
&-\varepsilon^2\Delta
\upsilon+\upsilon=\upsilon^p,\,\,\,\,\,\,\,
\upsilon>0\,\,\,\,\,
\textrm{in}\,\,\,\,\,\,
\mathcal{D},\\[1mm]
&\frac{\partial \upsilon}{\partial\nu}=0\,\ \
\,\,\,\qquad\qquad
\ \,\ \ \ \ \ \
\ \ \,\,\,\,
\textrm{on}\,\,\,\,
\partial\mathcal{D}.
\endaligned\right.
\endaligned
\end{equation}
This  equation has attracted considerable attention in the past three years because
its solutions  exhibit a  variety of concentration phenomena not only at one or more
points but also on higher dimensional subsets of $\overline{\mathcal{D}}$
as either
$\varepsilon$ tends to zero or  $p$ approaches the $(h+1)$-th  critical exponent $p_{h+1}^*$,
where  $p_{N-1}^*=+\infty$ and $p_{h+1}^*=(N-h+2)/(N-h-2)$  for any $0\leq h\leq N-3$.
The reader can refer to \cite{CW,GW1,GWW} for the subcritical case $p<p_{1}^*$,
to \cite{AM,R,W,WY} for the critical case $p=p_{1}^*$,
to \cite{CN,DMP,MW,RW1,W3} for the almost first critical case $p=p_1^*\pm d$
with $0<d\rightarrow0$,
to \cite{MM1,M1,MM2,MM3} for the
$(h+1)$-th subcritical case $p<p_{h+1}^*$ with $1\leq h\leq N-2$,
to \cite{DMM} for the
$(h+1)$-th critical case $p=p_{h+1}^*$ with $1\leq h\leq N-7$,
and to \cite{DMM1} for the almost
$(h+1)$-th critical case $p=p_{h+1}^*\pm d$
with  $1\leq h\leq N-7$ and $0<d\rightarrow0$.
In particular,  problem  (\ref{1.6}) admits a solution with arbitrarily many
mixed interior and boundary spikes, which has been shown in
\cite{GW1} for  $p<p_{1}^*$ fixed   but  $\varepsilon>0$ small enough,
and in \cite{MW} for $\varepsilon=1$  but
 $p\rightarrow p_{1}^*=+\infty$ with $N=2$.

It is very interesting  to point out that as a slight
but natural generalization of  (\ref{1.6})$|_{\varepsilon=1}$,
equation (\ref{1.4})  can  also be  viewed as a special case
of (\ref{1.6})$|_{\varepsilon=1}$
when  the domain $\mathcal{D}$ has some rotational symmetries.
Indeed,
take $n\geq2$ and $n\geq m\geq1$ as fixed integers.
Let $\Omega$  be  a smooth bounded domain in $\mathbb{R}^n$ such that
\begin{equation*}\label{}
\aligned
\overline{\Omega}\subset
\{
(x_1,\ldots,x_m,x')\in\mathbb{R}^m\times\mathbb{R}^{n-m}|\,
\,\,x_i>0,\,\,\,i=1,\ldots,m
\}.
\endaligned
\end{equation*}
Fix $k_1,\ldots,k_m\in\mathbb{N}\setminus\{0\}$ with $h:=k_1+\dots+k_m$
and set
\begin{equation*}\label{}
\aligned
\mathcal{D}:=
\big\{
(y_1,\ldots,y_m,x')\in\mathbb{R}^{k_1+1}\times\cdots\times\mathbb{R}^{k_m+1}\times\mathbb{R}^{n-m}|\,
\,(|y_1|,\ldots,|y_m|,x')\in\Omega
\big\}.
\endaligned
\end{equation*}
Then $\mathcal{D}$ is a smooth bounded domain in $\mathbb{R}^N$ with $N:=h+n$,
which is  invariant under the action of the group $\Upsilon:=\mathcal{O}(k_1+1)\times\cdots\times \mathcal{O}(k_m+1)$
on $R^N$ given by
$$
\aligned
(g_1,\ldots,g_m)(y_1,\ldots,y_m,x'):=
(g_1y_1,\ldots,g_my_m,x'),
\endaligned
$$
for each $g_i\in\mathcal{O}(k_i+1)$, $y_i\in\mathbb{R}^{k_i+1}$
and $x'\in\mathbb{R}^{n-m}$.
Note that $\mathcal{O}(k_i+1)$ is the group of linear isometries of $\mathbb{R}^{k_i+1}$
and $\mathbb{S}^{k_i}$ is the unit sphere in $\mathbb{R}^{k_i+1}$.
For any $p>1$   we shall seek $\Upsilon$-invariant solutions of problem (\ref{1.6})$|_{\varepsilon=1}$,
i.e. solutions $\upsilon$ of the form
\begin{equation}\label{1.7}
\aligned
\upsilon(y_1,\ldots,y_m,x')=u(|y_1|,\ldots,|y_m|,x').
\endaligned
\end{equation}
A simple calculation implies that $\upsilon$ solves problem (\ref{1.6})$|_{\varepsilon=1}$
 if and only if $u$ satisfies
\begin{equation}\label{1.8}
\left\{\aligned
&-\Delta u-\sum_{i=1}^m\frac{k_i}{\,x_i\,}\frac{\partial u}{\partial x_i}+u=u^p,
\,\,\,\,\,
u>0\,\,\,\,\,\,
\textrm{in}\,\,\,\,\,
\Omega,\\
&\frac{\partial u}{\partial\nu}=0\,\,\,\,\,
\ \,\ \,\qquad\qquad\qquad\qquad\,\,\,\,
\,\ \,\,\quad\quad\,\,\
\textrm{on}\,\,\,\,
\partial\Omega.
\endaligned\right.
\end{equation}
Thus  if we take anisotropic coefficient
\begin{equation}\label{1.9}
\aligned
a(x)=a(x_1,\ldots,x_m,x'):=x_1^{k_1}\dots x_m^{k_m},
\endaligned
\end{equation}
then problem (\ref{1.8}) can be rewritten as equation (\ref{1.4}).
Hence by considering rotational symmetry of $\mathcal{D}$,
 a fruitful approach for seeking solutions
of problem (\ref{1.6})$|_{\varepsilon=1}$ with concentration
along some $h$-dimensional minimal submanifolds of $\overline{\mathcal{D}}$
diffeomorphic to
$\mathbb{S}^{k_1}\times\ldots\times\mathbb{S}^{k_m}$
is to reduce it to produce  point-wise spiky or bubbling solutions of the
anisotropic equation (\ref{1.4}) in the domain
$\Omega$ of lower dimension. This approach, together
with some finite dimensional reduction arguments, has recently been taken
to construct solutions of (\ref{1.6})$|_{\varepsilon=1}$ concentrating
along an $h$-dimensional minimal submanifold of $\partial\mathcal{D}$,
which can be found in \cite{MP} for the almost
$(h+1)$-th critical case $p=p_{h+1}^*\pm d$
with $1\leq h\leq N-3$ and $0<d\rightarrow0$,
and in \cite{Z} for the slightly $(N-1)$-th subcritical case
$p\rightarrow p_{N-1}^*=+\infty$.

This paper is devoted to studying the existence and  concentration
behavior of  spiky solutions to the
anisotropic planar equation (\ref{1.1}) with large exponent
$p$ and
Hardy-H\'{e}non weight $|\cdot-q|^{2\alpha}$
 involving  $\alpha\in(-1,+\infty)\setminus\mathbb{N}$ only.
We try to use a finite dimensional reduction to investigate the effect of the interaction between anisotropic coefficient
$a(x)$ and   singular source $q$ on the existence
of concentrating solutions to problem (\ref{1.1}) when $p$ goes to $+\infty$.
As a result, we prove that if
$q\in\Omega$ is a strict local maximum point of $a(x)$,
there exists a family of  positive solutions with arbitrarily many
interior spikes accumulating to
$q$; while if $q\in\partial\Omega$ is a strict local maximum point of $a(x)$
and satisfies $\partial_{\nu}a(q):=\langle\nabla a(q),\,\nu(q)\rangle=0$,
such a problem has a family of  positive solutions with arbitrarily many mixed
interior and boundary spikes accumulating to $q$.
To state our results more precisely, we  first introduce some notations.
Let
$$
\aligned
\Delta_au=\frac1{a(x)}\nabla(a(x)\nabla u)=\Delta u+\nabla\log a(x)\nabla u
\endaligned
$$
and $G(x,y)$ be the   Green's function satisfying
\begin{equation}\label{1.2}
\left\{\aligned
&-\Delta_aG(x,y)+G(x,y)=\delta_y(x),\,\,\,\,\,\,\,
x\in\Omega,\\
&\frac{\partial G}{\partial\nu_x}(x,y)=0,
\qquad\,\,
\qquad\qquad\qquad\,\,
x\in\partial\Omega,
\endaligned\right.
\end{equation}
for each $y\in\oo$, then its
regular part  is defined depending on whether $y$ lies in
the domain or on its boundary as
\begin{equation}\label{1.3}
\aligned
H(x,y)=\left\{\aligned
&G(x,y)+\frac{1}{2\pi}\log|x-y|,\,\quad\,y\in\Omega,\\[1mm]
&G(x,y)+\frac1{\pi}\log|x-y|,\,\ \quad\,y\in\po.
\endaligned\right.
\endaligned
\end{equation}

\vspace{1mm}
\vspace{1mm}
\vspace{1mm}

\noindent{\bf Theorem 1.1.} {\it
Let $\alpha\in(-1,+\infty)\setminus\mathbb{N}$ and assume that $q\in\Omega$ is a strict local maximum
point of $a(x)$.
Then for any integer $m\geq1$, there exists $p_m>0$
such that for any $p>p_m$,  problem {\upshape (\ref{1.1})} has a
family of  positive solutions $u_p$
with  $m+1$ different interior spikes
which accumulate to $q$ as $p\rightarrow+\infty$. More precisely,
$$
\aligned
u_p(x)=&
\sum\limits_{i=1}^{m}\frac{1}{\gamma\mu_i^{2/(p-1)}\big|\xi^p_i-q\big|^{2\alpha/(p-1)}}\left[\,\log
\frac1{(\varepsilon^2\mu_i^2+|x-\xi^p_i|^2)^2}
+8\pi H(x,\xi^p_i)\right]\\
&+\frac{1}{\gamma\mu_0^{2/(p-1)}}\left[\,\log
\frac1{(\varepsilon^2\mu_0^2+|x-q|^{2(1+\alpha)})^{2}}
+8\pi(1+\alpha) H(x,q)\right]
+o\left(1\right),
\endaligned
$$
where $o(1)\rightarrow0$, as $p\rightarrow+\infty$, on each compact subset of
$\oo\setminus\{q,\xi^p_1,\ldots,\xi^p_m\}$,
the parameters $\gamma$, $\varepsilon$, $\mu_0$ and $\mu_i$,
$i=1,\ldots,m$  satisfy
$$
\aligned
\gamma=p^{p/(p-1)}\varepsilon^{2/(p-1)},\,\,
\quad\qquad\quad\,\,
\varepsilon=e^{-p/4},
\quad\qquad\quad\,\,
\frac1{C}<\mu_0<Cp^{C},
\quad\qquad\quad\,\,
\frac1{C}<\mu_i<Cp^{C},
\endaligned
$$
for some $C>0$,
and $(\xi^p_1,\ldots,\xi^p_m)\in (\Omega\setminus\{q\} )^m$ satisfies
$$
\aligned
\xi_{i}^p\rightarrow q\,\ \ \,\,\textrm{for all}\,\,\,i,\,\,
\ \ \ \,
\textrm{and}\,\,\ \ \,\,\,|\xi_{i}^p-\xi_{k}^p|>1/p^{2(m+1+\alpha)^2}\,\,\ \,\,\forall
\,\,\,
i\neq k.
\endaligned
$$
In particular, for any  $d>0$,
as $p\rightarrow+\infty$,
$$
\aligned
pa(x)|x-q|^{2\alpha}u_{p}^{p+1}\rightharpoonup8\pi e(m+1+\alpha)a(q)\delta_{q}
\,\quad\textrm{weakly in the sense of measure in}\,\,\,\,
\overline{\Omega},
\endaligned
$$
$$
\aligned
u_{p}\rightarrow0
\,\quad\textrm{uniformly in}\,\,\,\,
\overline{\Omega}\setminus B_{d}(q),
\endaligned
$$
but
$$
\aligned
\sup\limits_{x\in B_{1/p^{4(m+1+\alpha)^2}}(q)}u_p(x)\rightarrow\sqrt{e}
\qquad\quad
\textrm{and}
\qquad\quad
\sup\limits_{x\in B_{1/p^{4(m+1+\alpha)^2}}(\xi_i^p)}u_p(x)\rightarrow\sqrt{e},\,
\,\,\ \,i=1,\ldots,m.
\endaligned
$$
}

\vspace{1mm}
\vspace{1mm}
\vspace{1mm}

%
%
%

\noindent{\bf Theorem 1.2.}\,\,\,{\it
Let $\alpha\in(-1,+\infty)\setminus\mathbb{N}$ and  assume that $q\in\po$ is a strict local maximum point of $a(x)$
and satisfies $\partial_{\nu}a(q):=\langle\nabla a(q),\,\nu(q)\rangle=0$.
Then for any integers $m\geq1$ and $0\leq l\leq m$,
there exists $p_m>0$ such that for any $p>p_m$,  problem {\upshape (\ref{1.1})} has a
family of  positive solutions $u_p$
with $m-l+1$ different boundary spikes and $l$ different interior spikes
which accumulate to $q$ as $p\rightarrow+\infty$. More precisely,
$$
\aligned
u_p(x)=&
\sum\limits_{i=1}^{m}\frac{1}{\gamma\mu_i^{2/(p-1)}\big|\xi^p_i-q\big|^{2\alpha/(p-1)}}\left[\,\log
\frac1{(\varepsilon^2\mu_i^2+|x-\xi^p_i|^2)^2}
+c_i H(x,\xi^p_i)\right]\\
&+\frac{1}{\gamma\mu_0^{2/(p-1)}}\left[\,\log
\frac1{(\varepsilon^2\mu_0^2+|x-q|^{2(1+\alpha)})^{2}}
+4\pi(1+\alpha) H(x,q)\right]
+o\left(1\right),
\endaligned
$$
where $o(1)\rightarrow0$, as $p\rightarrow+\infty$, on each compact subset of
$\oo\setminus\{q,\xi^p_1,\ldots,\xi^p_m\}$,
the parameters $\gamma$, $\varepsilon$, $\mu_0$ and $\mu_i$,
$i=1,\ldots,m$  satisfy
$$
\aligned
\gamma=p^{p/(p-1)}\varepsilon^{2/(p-1)},\,\,
\quad\qquad\quad\,\,
\varepsilon=e^{-p/4},
\quad\qquad\quad\,\,
\frac1{C}<\mu_0<Cp^{C},
\quad\qquad\quad\,\,
\frac1{C}<\mu_i<Cp^{C},
\endaligned
$$
for some $C>0$,
$(\xi^p_1,\ldots,\xi^p_m)\in\Omega^l\times(\po)^{m-l}$  satisfies
$$
\aligned
\xi^p_i\rightarrow q
\quad\forall\,\,i,
\qquad
|\xi^p_i-\xi^p_k|>1/p^{2(m+1+\alpha)^2}
\quad\forall\,\,i\neq k,
\qquad
\textrm{and}
\qquad
\dist(\xi^p_i,\po)>1/p^{2(m+1+\alpha)^2}
\quad\forall\,\,i=1,\ldots,l,
\endaligned
$$
and $c_i=8\pi$ for $i=1,\ldots,l$, while  $c_i=4\pi$
for $i=l+1,\ldots,m$.
In particular, for any  $d>0$,
as $p\rightarrow+\infty$,
$$
\aligned
pa(x)|x-q|^{2\alpha}u_{p}^{p+1}\rightharpoonup4\pi e(m+l+1+\alpha)a(q)\delta_{q}
\,\quad\textrm{weakly in the sense of measure in}\,\,\,\,
\overline{\Omega},
\endaligned
$$
$$
\aligned
u_{p}\rightarrow0
\,\quad\textrm{uniformly in}\,\,\,\,
\overline{\Omega}\setminus B_{d}(q),
\endaligned
$$
but
$$
\aligned
\sup\limits_{x\in\oo\cap B_{1/p^{4(m+1+\alpha)^2}}(q)}u_p(x)\rightarrow\sqrt{e}
\qquad\quad
\textrm{and}
\qquad\quad
\sup\limits_{x\in\oo\cap B_{1/p^{4(m+1+\alpha)^2}}(\xi_i^p)}u_p(x)\rightarrow\sqrt{e},\,
\,\,\ \,i=1,\ldots,m.
\endaligned
$$
}

\vspace{1mm}

In fact,  the assumption in
Theorem  $1.2$  can be split into  the following two cases:\\
(C1) \,
$q\in\po$ is a strict
local maximum point of $a(x)$ restricted on $\po$;\\
(C2) \,
$q\in\po$ is a strict local maximum
point of $a(x)$ restricted in $\Omega$
and satisfies $\partial_{\nu}a(q):=\langle\nabla a(q),\,\nu(q)\rangle=0$.\\
Arguing exactly along the sketch of the proof of Theorem $1.2$, we can easily find that if
(C1) holds,  then problem (\ref{1.1}) has solutions with
arbitrarily many
boundary spikes which  accumulate to $q$ along $\po$; while
if (C2) holds,  then problem (\ref{1.1}) has solutions with
arbitrarily many
interior spikes which  accumulate to $q$ along the
inner normal direction of $\po$.

For the case $m=0$,  we have  the following two
results which correspond to
Theorems $1.1$ and $1.2$, respectively.

\vspace{1mm}
\vspace{1mm}
\vspace{1mm}
\vspace{1mm}

\noindent{\bf Theorem 1.3.} {\it
Let $\alpha\in(-1,+\infty)\setminus\mathbb{N}$ and $q\in\Omega$.
Then   there exists $p_0>0$
such that for any $p>p_0$,  problem {\upshape(\ref{1.1})} has
a family of   positive solutions $u_p$  such that as $p$ tends to $+\infty$,
$$
\aligned
u_p(x)=\frac{1}{\gamma\mu_0^{2/(p-1)}}\left[\,\log
\frac1{(\varepsilon^2\mu_0^2+|x-q|^{2(1+\alpha)})^{2}}
+8\pi(1+\alpha) H(x,q)\right]
+o\left(1\right),
\endaligned
$$
 uniformly on
each compact subset of $\overline{\Omega}\setminus\{q\}$,
where  the parameters $\gamma$, $\varepsilon$  and $\mu_0$  satisfy
$$
\aligned
\gamma=p^{p/(p-1)}\varepsilon^{2/(p-1)},\,\,
\quad\quad\qquad\quad\,\,
\varepsilon=e^{-p/4},
\quad\qquad\quad\quad\,\,
1/{C}<\mu_0<C,
\endaligned
$$
for some $C>0$.
In particular, for any   $d>0$,
as $p\rightarrow+\infty$,
$$
\aligned
pa(x)|x-q|^{2\alpha}u_{p}^{p+1}\rightharpoonup8\pi e(1+\alpha)a(q)\delta_{q}
\,\quad\textrm{weakly in the sense of measure in}\,\,\,\,
\overline{\Omega},
\endaligned
$$
$$
\aligned
u_{p}\rightarrow0
\,\quad\textrm{uniformly in}\,\,\,\,
\overline{\Omega}\setminus B_{d}(q),
\endaligned
$$
but
$$
\aligned
\sup\limits_{x\in B_{d}(q)}u_p(x)\rightarrow\sqrt{e}.
\endaligned
$$
}

\vspace{1mm}

\noindent{\bf Theorem 1.4.} {\it
Let $\alpha\in(-1,+\infty)\setminus\mathbb{N}$ and $q\in\partial\Omega$.
Then   there exists $p_0>0$
such that for any $p>p_0$,  problem {\upshape(\ref{1.1})} has
a family of   positive  solutions $u_p$  such that as $p$ tends to $+\infty$,
$$
\aligned
u_p(x)=\frac{1}{\gamma\mu_0^{2/(p-1)}}\left[\,\log
\frac1{(\varepsilon^2\mu_0^2+|x-q|^{2(1+\alpha)})^{2}}
+4\pi(1+\alpha) H(x,q)\right]
+o\left(1\right),
\endaligned
$$
uniformly on
each compact subset of $\overline{\Omega}\setminus\{q\}$,
where  the parameters $\gamma$, $\varepsilon$  and $\mu_0$  satisfy
$$
\aligned
\gamma=p^{p/(p-1)}\varepsilon^{2/(p-1)},\,\,
\quad\quad\qquad\quad\,\,
\varepsilon=e^{-p/4},
\quad\qquad\quad\quad\,\,
1/{C}<\mu_0<C,
\endaligned
$$
for some $C>0$.
In particular, for any  $d>0$,
as $p\rightarrow+\infty$,
$$
\aligned
pa(x)|x-q|^{2\alpha}u_{p}^{p+1}\rightharpoonup4\pi e(1+\alpha)a(q)\delta_{q}
\,\quad\textrm{weakly in the sense of measure in}\,\,\,\,
\overline{\Omega},
\endaligned
$$
$$
\aligned
u_{p}\rightarrow0
\,\quad\textrm{uniformly in}\,\,\,\,
\overline{\Omega}\setminus B_{d}(q),
\endaligned
$$
but
$$
\aligned
\sup\limits_{x\in\overline{\Omega}\cap  B_{d}(q)}u_p(x)\rightarrow\sqrt{e}.
\endaligned
$$
}

\vspace{1mm}

Let us comment  that
for the case $m=0$, by arguing simply along some
initial procedures  of the construction of solutions
of problem (\ref{1.1})
we can prove
the corresponding results in Theorems $1.3$ and $1.4$.
This shows that problem (\ref{1.1}) always admits a family of
positive solutions
concentrating only at singular source $q$ whether
$q$ is an isolated local maximum point of  $a(x)$
or not.

Finally, let us remark about the analogy and difference
existing
between our results and  those known for
the planar Neumann problem with  Hardy-H\'{e}non weight and large exponent
\begin{equation}\label{*}
\left\{\aligned
&-\Delta u+u=|x-q|^{2\alpha}u^p,\,\,\,\,\,\,
u>0\,\,\,\,\,\,\,
\textrm{in}\,\,\,\,\,
\Omega,\\[1mm]
&\frac{\partial u}{\partial\nu}=0\,\,
\qquad\qquad\qquad\qquad\qquad
\,\,\,\ \ \,\,\,\,
\textrm{on}\,\,\,
\partial\Omega,
\endaligned\right.
\end{equation}
where
$\Omega$ is a smooth  bounded   domain in $\mathbb{R}^2$,
$\alpha\in(-1,+\infty)\setminus\mathbb{N}$,
$q\in\overline{\Omega}$,
$p$ is a large exponent and $\nu$ denotes the outer unit  normal vector to  $\partial\Omega$.
Just like that in equation (\ref{1.1}),
the presence of Hardy-H\'{e}non weight
can also produce significant influence on
the existence of a solution of problem (\ref{*})
with spiky profile at each singular source $q\in\oo$,
which has been proven in
\cite{ZY} that problem (\ref{*}) always admits
a family of positive solutions $u_{p}$
 with arbitrarily
many mixed interior and boundary spikes  involving  any
$q\in\oo$
when  $p$ tends to $+\infty$.
However, due to the occurrence of Hardy-H\'{e}non weight,
it is necessary to point out that although the anisotropic planar
equation (\ref{1.1}) is seemingly similar to problem (\ref{*}),
equation (\ref{1.1})  can  not be  viewed as a special version
of (\ref{*}) in higher dimension
even if  the domain has some rotational symmetries. This seems to
imply  that unlike those for solutions of problem (\ref{*}) in
\cite{ZY} whose multiple spikes are  completely determined by the geometry of the domain,
the location of multiple spikes in solutions of equation (\ref{1.1}) may be
only characterized by anisotropic coefficient $a(x)$ and singular source $q\in\oo$.
In fact, from Theorems $1.1$ and $1.2$ it follows  that  if
$q$ is an isolated local maximum point of $a(x)$ over $\oo$ and satisfies
 $\langle\nabla a(q),\,\nu(q)\rangle=0$ if $q\in\po$, then equation (\ref{1.1})
always admits a family of positive solutions $u_p$ with arbitrarily many interior
and boundary spikes accumulating to $q$ when  $p$ tends to $+\infty$.
Thus as a by-product of our results  we  find that the presence of anisotropic coefficient
$a(x)$ can lead the concentration point $q$ of positive solutions to (\ref{1.1})
to be non-simple (or accumulated) in the sense that there exist more than one (or arbitrarily many) spiky
sequences of points converging at the same point.

The proof of  our  results  relies on a
very well-known Lyapunov-Schmidt finite dimensional reduction. In Section $2$ we
exactly describe an approximate solution of equation (\ref{1.1}).
Then we rewrite  equation (\ref{1.1}) in terms of a linearized operator for which
a solvability theory, subject to suitable orthogonality conditions,
is performed  through solving a linearized  problem and an auxiliary nonlinear problem in Section $3$.
In Section $4$ we reduce the problem of finding concentrating solutions
of (\ref{1.1}) to that of finding a critical point of a finite dimensional function
and give its asymptotic expansion.
In the last section, we give the proofs of  Theorems $1.1$ and $1.2$.

\vspace{1mm}
\vspace{1mm}

\section{Approximating solutions}
The basic cells for the construction of  an approximate solution of  problem (\ref{1.1})
are given by  two standard bubbles
\begin{equation}\label{2.1}
\aligned
U_{\delta}(x)=\log\frac{8(1+\alpha)^2\delta^{2(1+\alpha)}}{(\delta^{2(1+\alpha)}+|x|^{2(1+\alpha)})^2}
\qquad\qquad
\,\,\textrm{and}\,\,
\qquad\qquad
V_{\delta,\xi}(x)=\log\frac{8\delta^2}{(\delta^2+|x-\xi|^2)^2},
\endaligned
\end{equation}
with $\alpha\in(-1,+\infty)\setminus\mathbb{N}$,
$\delta>0$ and $\xi\in\mathbb{R}^2$,
which, respectively,  are  all  solutions of the
following two equations:
\begin{equation}\label{2.2}
\left\{\aligned
&-\Delta
u=|x|^{2\alpha}e^u\,\ \,\,\,\,
\textrm{in}\,\,\,\,
\mathbb{R}^2,\\
&\int_{\mathbb{R}^2}|x|^{2\alpha}e^u<+\infty,
\endaligned\right.
\,\,\,\quad
\qquad
\textrm{and}
\qquad
\,\,\,\quad
\left\{\aligned
&-\Delta
u=e^u\,\ \,\,\,\,
\textrm{in}\,\,\,\,
\mathbb{R}^2,\\
&\int_{\mathbb{R}^2}e^u<+\infty,
\endaligned\right.
\end{equation}
(see \cite{CL,CL1,PT}).
Let us define the configuration space in which the concentration points we
try to seek  belong to
\begin{eqnarray}\label{2.3}
\mathcal{O}_p(q):=\left\{\,\xi=(\xi_1,\ldots,\xi_m)\in\big(B_d(q)\cap\Omega\big)^l
\times\big(B_d(q)\cap\po\big)^{m-l}
\left|\,\,
\min_{i=1,\ldots,m}|\xi_i-q|>\frac{1}{p^\kappa},
\right.\right.
\nonumber&&\\[1mm]
\left.\,
\min_{i,j=1,\ldots,m,\,i\neq j}|\xi_i-\xi_j|>\frac{1}{p^\kappa},
\,\,\,\,\qquad\,\,\,
\min_{1\leq i\leq l}\dist(\xi_i,\po)>\frac{1}{p^\kappa}
\right\},
\qquad\qquad\qquad\,\,\,\,
&&
\end{eqnarray}
where
$d>0$ is a sufficiently small but fixed number independent of $p$,
$l=m$ if $q\in\Omega$  while $l=0,1,\ldots,m$ if $q\in\partial\Omega$,
and
$\kappa$ is given by
\begin{eqnarray}\label{2.4}
\kappa=2(m+1+\alpha)^2.
\end{eqnarray}
Let us fix $m\in\mathbb{N}$, $q\in\oo$  and
$\xi=(\xi_1,\ldots,\xi_m)\in\mathcal{O}_p(q)$.
For numbers  $\mu_0>0$ and $\mu_i>0$,
$i=1,\ldots,m$, yet to be chosen, we set
\begin{equation}\label{2.5}
\aligned
\gamma=p^{\frac{p}{p-1}}e^{-\frac{p}{2(p-1)}},
\,\,\qquad\,\,\varepsilon=e^{-\frac{1}{4}p},
\,\,\qquad\,\,\,\rho_0=\varepsilon^{\frac{1}{1+\alpha}},
\,\,\qquad\,\,\,\upsilon_0=\mu_0^{\frac{1}{1+\alpha}},
\,\,\qquad\,\delta_0=\rho_0\upsilon_0,
\,\,\qquad\,\delta_i=\varepsilon\mu_i,
\endaligned
\end{equation}
and
\begin{equation}\label{2.6}
\aligned
U_0(x)=\frac{1}{\gamma\mu_0^{2/(p-1)}}\left[U_{\delta_0}(x-q)
+\frac1p\omega_1\left(\frac{x-q}{\delta_0}\right)
+\frac1{p^2}\omega_2\left(\frac{x-q}{\delta_0}\right)
\right],
\endaligned
\end{equation}
and
\begin{equation}\label{2.7}
\aligned
U_i(x)=\frac{1}{\gamma\mu_i^{2/(p-1)}
\big|\xi_i-q\big|^{2\alpha/(p-1)}}\left[V_{\delta_i,\xi_i}(x)
+\frac1p\widetilde{\omega}_1\left(\frac{x-\xi_i}{\delta_i}\right)
+\frac1{p^2}\widetilde{\omega}_2\left(\frac{x-\xi_i}{\delta_i}\right)
\right].
\endaligned
\end{equation}
Here, $\omega_j$ and $\widetilde{\omega}_j$, $j=1, 2$, respectively, are radial solutions of
\begin{equation}\label{2.8}
\aligned
\Delta\omega_j+\frac{8(1+\alpha)^2|z|^{2\alpha}}{(1+|z|^{2(1+\alpha)})^2}\omega_j
=\frac{8(1+\alpha)^2|z|^{2\alpha}}{(1+|z|^{2(1+\alpha)})^2}f_j(|z|)
\quad\,\textrm{in}\,\,\,\,\mathbb{R}^2,
\endaligned
\end{equation}
and
\begin{equation}\label{2.9}
\aligned
\Delta\widetilde{\omega}_j+\frac{8}{(1+|z|^2)^2}\widetilde{\omega}_j=\frac{8}{(1+|z|^2)^2}\widetilde{f}_j(|z|)
\quad\,\textrm{in}\,\,\,\,\mathbb{R}^2,
\endaligned
\end{equation}
where
\begin{equation}\label{2.10}
\aligned
f_1=\frac12U_{1}^2,\,\quad\quad\,
f_2=\omega_1U_{1}-\frac13U_{1}^3-\frac12\omega_1^2-\frac18U_{1}^4+\frac{1}{2}\omega_1U_{1}^2,
\endaligned
\end{equation}
and
\begin{equation}\label{2.11}
\aligned
\widetilde{f}_1=\frac12V_{1,0}^2,\,\quad\quad\,
\widetilde{f}_2=\widetilde{\omega}_1V_{1,0}-\frac13V_{1,0}^3-\frac12\widetilde{\omega}_1^2-\frac18V_{1,0}^4+\frac{1}{2}\widetilde{\omega}_1V_{1,0}^2.
\endaligned
\end{equation}
Furthermore, by \cite{CI,EMP,EPW} it follows  that for each $j=1,2$ and $r=|z|$, as $r\rightarrow+\infty$,
\begin{equation}\label{2.12}
\aligned
\omega_j(r)=\frac{C_j}{2(1+\alpha)}\log\big(1+r^{2(1+\alpha)}\big)+O\left(\frac1{1+r^{1+\alpha}}\right),
\,\,\quad\,\,
\partial_{r}\omega_j(r)=\frac{C_jr^{1+2\alpha}}{1+r^{2(1+\alpha)}}+O\left(\frac{1}{1+r^{2+\alpha}}\right),
\endaligned
\end{equation}
and
\begin{equation}\label{2.13}
\aligned
\widetilde{\omega}_j(r)=\frac{\widetilde{C}_j}{2}\log(1+r^2)+O\left(\frac1{1+r}\right),
\,\,\quad\,\,
\partial_{r}\widetilde{\omega}_j(r)=\frac{\widetilde{C}_jr}{1+r^2}+O\left(\frac{1}{1+r^2}\right),
\endaligned
\end{equation}
where
\begin{equation}\label{2.14}
\aligned
C_{j}=8(1+\alpha)^2\int_0^{+\infty}t^{1+2\alpha}\frac{\,t^{2(1+\alpha)}-1\,}
{\,(t^{2(1+\alpha)}+1)^3\,}f_{j}(t)dt,
\qquad\qquad
\widetilde{C}_j=8\int_{0}^{\infty}t\frac{t^2-1}{(t^2+1)^3}\widetilde{f}_j(t)dt.
\endaligned
\end{equation}
In particular,
\begin{equation}\label{2.15}
\aligned
C_1=12(1+\alpha)-4(1+\alpha)\log8(1+\alpha)^2,
\qquad\qquad
\widetilde{C}_1=12-4\log8,
\endaligned
\end{equation}
and
\begin{eqnarray}\label{2.16}
\omega_1(z)=\frac12U_1^2(z)+6\log\big(|z|^{2(1+\alpha)}+1\big)
+\frac{2\log8(1+\alpha)^2-10}{|z|^{2(1+\alpha)}+1}+\frac{|z|^{2(1+\alpha)}-1}{|z|^{2(1+\alpha)}+1}
\times\left\{-\frac12\log^28(1+\alpha)^2\right.&&\nonumber\\
\left.
+2\log^2\big(|z|^{2(1+\alpha)}+1\big)+4\int_{|z|^{2(1+\alpha)}}^{+\infty}\frac{ds}{s+1}
\log\frac{s+1}{s}
-8(1+\alpha)\log|z|\log\big(|z|^{2(1+\alpha)}+1\big)\right\},
\,\,\,\,&&
\end{eqnarray}
and
\begin{eqnarray}\label{2.17}
\widetilde{\omega}_1(z)=\frac12V_{1,0}^2(z)+6\log(|z|^2+1)
+\frac{2\log8-10}{|z|^2+1}+\frac{|z|^2-1}{|z|^2+1}
\times\left\{-\frac12\log^28
\quad\,\,
\right.&&\nonumber\\
\left.
+2\log^2(|z|^2+1)+4\int_{|z|^2}^{+\infty}\frac{ds}{s+1}
\log\frac{s+1}{s}
-8\log|z|\log(|z|^2+1)\right\}.
&&
\end{eqnarray}

We define the approximate solution of  problem (\ref{1.1}) as
\begin{equation}\label{2.18}
\aligned
U_\xi(x):=\sum_{i=0}^mPU_i(x)=\sum_{i=0}^m\big[U_i(x)+H_i(x)\big],
\endaligned
\end{equation}
where $H_i$  is a correction term defined as the solution of
\begin{equation}\label{2.19}
\aligned
\left\{\aligned
&-\Delta_a
H_i+H_i=\nabla\log a(x)\nabla U_i-U_i\,\,\,\,\,\,
\textrm {in}\,\,\,\,\,\,\Omega,\\[1mm]
&\frac{\partial H_i}{\partial \nu}=-\frac{\partial U_i}{\partial\nu}\,
\,\qquad\qquad\qquad\qquad\qquad\,\,\,
\textrm{on}\,\,\,\,\po.
\endaligned\right.
\endaligned
\end{equation}
In order to state the asymptotic expansion of the functions $H_i$ in terms of $\xi_i$, $\delta_i$ and $p>1$ large enough, we  first use the convention
\begin{equation}\label{2.20}
\aligned
c_0=\left\{
\aligned
&8\pi(1+\alpha),
\quad\,\,\textrm{if}\quad
q\in\Omega,
\\[0.5mm]
&4\pi(1+\alpha),
\quad\,\,\textrm{if}\quad
q\in\partial\Omega,
\endaligned
\right.
\endaligned
\qquad\quad\qquad
\aligned
c_i=\left\{
\aligned
&8\pi,
\quad\,\,
\textrm{if}\quad
i=1,\ldots,l,
\\[0.5mm]
&4\pi,
\quad\,\,
\textrm{if}\quad
i=l+1,\ldots,m.
\endaligned
\right.
\endaligned
\end{equation}
Then we have the following Lemma whose proof is listed in the Appendix B.

\vspace{1mm}
\vspace{1mm}
\vspace{1mm}
\vspace{1mm}

\noindent{\bf Lemma 2.1.}\,\,{\it
For any points $\xi=(\xi_1,\ldots,\xi_m)\in\mathcal{O}_p(q)$
and any $0<\beta<1$,
then we have that
\begin{eqnarray}\label{2.21}
H_0(x)=\frac{1}{\gamma\mu_0^{2/(p-1)}}\left[
\left(
1-\frac{C_1}{4(1+\alpha)p}-\frac{C_2}{4(1+\alpha)p^2}
\right)c_0 H(x,q)
-\log\left(8(1+\alpha)^2\delta_0^{2(1+\alpha)}\right)
\right.&&\nonumber\\[1mm]
\left.
+\left(\frac{C_1}{p}+\frac{C_2}{p^2}\right)\log\delta_0
+O\left(\delta_0^{\beta/2}\right)
\right],
\qquad\qquad\qquad\qquad\qquad\qquad\qquad\qquad
\qquad\,\,\,
&&
\end{eqnarray}
and for each $i=1,\ldots,m$,
\begin{eqnarray}\label{2.22}
H_i(x)=\frac{1}{\gamma\mu_i^{2/(p-1)}\big|\xi_i-q\big|^{2\alpha/(p-1)}}\left[
\left(
1-\frac{\widetilde{C}_1}{4p}-\frac{\widetilde{C}_2}{4p^2}
\right)c_i H(x,\xi_i)
-\log(8\delta_i^2)
+\left(\frac{\widetilde{C}_1}{p}+\frac{\widetilde{C}_2}{p^2}\right)\log\delta_i
\right.
&&\nonumber\\[1mm]
\left.
+O\left(\delta_i^{\beta/2}\right)
\right],
\qquad\qquad\qquad\qquad\qquad\qquad
\qquad\qquad\qquad\qquad\qquad\qquad
\quad\,\qquad\qquad\qquad
&&
\end{eqnarray}
uniformly in $\oo$,
where
$H$ is the regular part of Green's function defined in
{\upshape(\ref{1.3})}  and for each $j=1,2$,
$C_j$ and $\widetilde{C}_j$
are the constants defined in \upshape{(\ref{2.14})}.
}

\vspace{1mm}
\vspace{1mm}
\vspace{1mm}
\vspace{1mm}

From  Lemma 2.1 we can easily check that away from singular source $q$ and each point $\xi_i$,
namely  $|x-q|\geq 1/p^{2\kappa}$ and $|x-\xi_i|\geq 1/p^{2\kappa}$ for each $i=1,\ldots,m$,
\begin{eqnarray}\label{2.26}
U_{\xi}(x)=
\frac{1}{\gamma\mu_0^{2/(p-1)}}
\left[\left(1
-\frac{C_1}{4(1+\alpha)p}-\frac{C_2}{4(1+\alpha)p^2}\right)c_0 G(x,q)
+O\left(p^{2\kappa(1+\alpha)-1}\delta_0^{1+\alpha}+\delta_0^{\beta/2}\right)\right]
\,\,&&\nonumber\\[1mm]
+\sum\limits_{i=1}^{m}
\frac{1}{\gamma\mu_i^{2/(p-1)}\big|\xi_i-q\big|^{2\alpha/(p-1)}}
\left[\left(1
-\frac{\widetilde{C}_1}{4p}-\frac{\widetilde{C}_2}{4p^2}\right)c_i G(x,\xi_i)
+O\left(\delta_i^{\beta/2}\right)\right].
\qquad\qquad\,\,
&&
\end{eqnarray}
If
$|x-q|< 1/p^{2\kappa}$,  by using  (\ref{2.1}),
(\ref{2.6}),  (\ref{2.21}) and the
fact that $H(\cdot,q)\in C^{\beta}(\oo)$ for any
$\beta\in(0,1)$
we obtain
$$
\aligned
PU_0(x)=&\frac{1}{\gamma\mu_0^{2/(p-1)}}
\left[U_{1}\left(\frac{x-q}{\delta_0}\right)
+\frac1p\omega_1\left(\frac{x-q}{\delta_0}\right)
+\frac1{p^2}\omega_2\left(\frac{x-q}{\delta_0}\right)+\left(
1-\frac{C_1}{4(1+\alpha)p}-\frac{C_2}{4(1+\alpha)p^2}
\right)c_0 H(q,q)\right.\\[1mm]
&
-\left.\log\left(8(1+\alpha)^2\delta_0^{4(1+\alpha)}\right)+\left(\frac{C_1}{p}+\frac{C_2}{p^2}\right)\log\delta_0
+O\left(
|x-q|^\beta
+\delta_0^{\beta/2}\right)
\right],
\endaligned
$$
and for any $k\neq 0$, by  (\ref{2.1}), (\ref{2.7}), (\ref{2.13}) and (\ref{2.22}),
$$
\aligned
PU_k(x)=&\frac{1}{\gamma\mu_k^{2/(p-1)}\big|\xi_k-q\big|^{2\alpha/(p-1)}}
\left[
\log\frac{8\delta_k^2}{(\delta_k^2+|x-\xi_k|^2)^2}
+\frac1p\widetilde{\omega}_1\left(\frac{x-\xi_k}{\delta_k}\right)
+\frac1{p^2}\widetilde{\omega}_2\left(\frac{x-\xi_k}{\delta_k}\right)
\right.\\[1mm]
&\left.
+\left(
1-\frac{\widetilde{C}_1}{4p}-\frac{\widetilde{C}_2}{4p^2}
\right)c_k H(x,\xi_k)
-\log(8\delta_k^2)+\left(\frac{\widetilde{C}_1}{p}+\frac{\widetilde{C}_2}{p^2}\right)\log\delta_k
+O\left(\delta_k^{\beta/2}\right)
\right]\\[1mm]
=&\frac{1}{\gamma\mu_k^{2/(p-1)}\big|\xi_k-q\big|^{2\alpha/(p-1)}}
\left[
\left(
1-\frac{\widetilde{C}_1}{4p}-\frac{\widetilde{C}_2}{4p^2}
\right)c_k G(q,\xi_k)
+O\left(
|x-q|^\beta
+\delta_k^{\beta/2}\right)
\right].
\endaligned
$$
Hence for $|x-q|< 1/p^{2\kappa}$, by (\ref{2.5}),
\begin{eqnarray}\label{2.27}
U_{\xi}(x)=\frac{1}{\gamma\mu_0^{2/(p-1)}}
\left[\,p+U_{1}\left(\frac{x-q}{\delta_0}\right)
+\frac1p\omega_1\left(\frac{x-q}{\delta_0}\right)
+\frac1{p^2}\omega_2\left(\frac{x-q}{\delta_0}\right)
+O\left(
|x-q|^\beta
+\sum_{k=0}^m\delta_k^{\beta/2}
\right)
\right]
\end{eqnarray}
is an appropriate approximation for a solution of problem (\ref{1.1})
near singular source $q$
provided that the concentration
parameter $\mu_0$  satisfies the nonlinear relation
\begin{eqnarray}\label{2.28}
\log\big(8(1+\alpha)^2\mu_0^4\big)=\left(
1-\frac{C_1}{4(1+\alpha)p}-\frac{C_2}{4(1+\alpha)p^2}
\right)c_0 H(q,q)
+\left(\frac{C_1}{p}+\frac{C_2}{p^2}\right)\log\delta_0
&&\nonumber\\[1mm]
+\left(
1-\frac{\widetilde{C}_1}{4p}-\frac{\widetilde{C}_2}{4p^2}
\right)\sum_{k=1}^m
\frac{\mu_0^{2/(p-1)}}{\,\mu_k^{2/(p-1)}\big|\xi_k-q\big|^{2\alpha/(p-1)}\,}c_k
G(q,\xi_k).
\quad
\,&&
\end{eqnarray}
Similarly, while  if
$|x-\xi_i|< 1/p^{2\kappa}$  with some $i\in\{1,\ldots,m\}$,
\begin{eqnarray}\label{2.29}
U_{\xi}(x)=\frac{1}{\gamma\mu_i^{2/(p-1)}\big|\xi_i-q\big|^{2\alpha/(p-1)}}
\left[\,p+V_{1,0}\left(\frac{x-\xi_i}{\delta_i}\right)
+\frac1p\widetilde{\omega}_1\left(\frac{x-\xi_i}{\delta_i}\right)
+\frac1{p^2}\widetilde{\omega}_2\left(\frac{x-\xi_i}{\delta_i}\right)
\right.
&&
\nonumber\\
\left. +\,O\left(
|x-\xi_i|^\beta
+
p^{2\kappa(1+\alpha)-1}\delta_0^{1+\alpha}
+\sum_{k=0}^m\delta_k^{\beta/2}
\right)
\right]
\qquad\qquad\qquad
\qquad\qquad\qquad\qquad\,\,\,\,
&&
\end{eqnarray}
is an appropriate approximation for a solution of problem (\ref{1.1})
near the point $\xi_i$
provided that for each $i=1,\ldots,m$, the
concentration
parameter $\mu_i$   satisfies the nonlinear system
\begin{eqnarray}\label{2.30}
\log\big(8\mu_i^4\big)=\left(
1-\frac{\widetilde{C}_1}{4p}-\frac{\widetilde{C}_2}{4p^2}
\right)c_i H(\xi_i,\xi_i)
+\left(\frac{\widetilde{C}_1}{p}+\frac{\widetilde{C}_2}{p^2}\right)\log\delta_i
\qquad\qquad\qquad\qquad\,\,\,\,\,
&&\nonumber\\[1mm]
+\left(
1-\frac{\widetilde{C}_1}{4p}-\frac{\widetilde{C}_2}{4p^2}
\right)\sum_{k=1,\,k\neq i}^m
\frac{\mu_i^{2/(p-1)}\big|\xi_i-q\big|^{2\alpha/(p-1)}}{\,\mu_k^{2/(p-1)}\big|\xi_k-q\big|^{2\alpha/(p-1)}\,}
c_k G(\xi_i,\xi_k)
\quad\,\,\,\,\,\,\,\,
&&\nonumber\\[1mm]
+\left(
1-\frac{C_1}{4(1+\alpha)p}-\frac{C_2}{4(1+\alpha)p^2}
\right)\frac{\,\mu_i^{2/(p-1)}\big|\xi_i-q\big|^{2\alpha/(p-1)}\,}{\,\mu_0^{2/(p-1)}\,}
c_0G(\xi_i,q).
&&
\end{eqnarray}
Indeed, the parameters  $\mu=(\mu_0,\mu_1,\ldots,\mu_m)$
are well defined in systems {\upshape(\ref{2.28})} and {\upshape(\ref{2.30})}
 under the certain region, which is stated in
the following lemma and whose proof is postponed  in the Appendix B.

\vspace{1mm}
\vspace{1mm}
\vspace{1mm}
\vspace{1mm}

\noindent{\bf Lemma 2.2.}\,\,{\it For
any points $\xi=(\xi_1,\ldots,\xi_m)\in\mathcal{O}_p(q)$
and any  $p>1$ large enough,
systems {\upshape(\ref{2.28})} and {\upshape(\ref{2.30})}
have a unique solution $\mu=(\mu_0,\mu_1,\ldots,\mu_m)$
satisfying
\begin{equation}\label{2.31}
\aligned
1/C\leq\mu_i\leq Cp^{C}
\quad\quad\,\,
\textrm{and}
\quad\quad\,\,
\big|D_{\xi}\log\mu_i\big|\leq Cp^{\kappa},
\quad\quad\,\,
\forall\,\,i=0,1,\ldots,m,
\endaligned
\end{equation}
for some  $C>0$. Moreover,
\begin{equation}\label{2.36}
\aligned
\mu_0=\mu_0(p,\xi)\equiv e\large^{-\frac{3}{4}+\frac14c_0H(q,q)
+\frac14\sum_{k=1}^m
c_k G(q,\xi_k)}\left[\,1+O\left(\frac{\log^2p}p\right)\right],
\endaligned
\end{equation}
and for any $i=1,\ldots,m$,
\begin{equation}\label{2.37}
\aligned
\mu_i=\mu_i(p,\xi)\equiv e\large^{-\frac{3}{4}+\frac14c_iH(\xi_i,\xi_i)
+\frac14
c_0G(\xi_i,q)
+\frac14\sum_{k=1,\,k\neq i}^m
c_k G(\xi_i,\xi_k)}\left[\,1+O\left(\frac{\log^2p}p\right)\right].
\endaligned
\end{equation}
}

\vspace{1mm}
\vspace{1mm}

\noindent{\bf Remark 2.3.}\,\,For any $p>1$ large enough, we see  that
if $|x-q|=\delta_0|z|< 1/p^{2\kappa}$,
by (\ref{2.1}), (\ref{2.5}), (\ref{2.12}) and (\ref{2.31}),
$$
\aligned
p+U_{1}(|z|)
+\frac1p\omega_1(|z|)
+\frac1{p^2}\omega_2(|z|)
&= p-2\log\left(1+|z|^{2(1+\alpha)}\right)+O\left(1\right)\\[0.1mm]
&\geq8\kappa(1+\alpha)\log p+
4\log\mu_0
+O\left(1\right)>7\kappa(1+\alpha)\log p,
\endaligned
$$
which, together with (\ref{2.27}),  easily implies that
$0<U_\xi\leq2\sqrt{e}$ in $B_{1/p^{2\kappa}}(q)$,
and
$\sup_{B_{1/p^{2\kappa}}(q)}U_\xi\rightarrow\sqrt{e}$
as $p\rightarrow+\infty$.
While if $|x-\xi_i|=\delta_i|\tilde{z}|< 1/p^{2\kappa}$ for each $i=1,\ldots,m$,
by (\ref{2.1}), (\ref{2.5}), (\ref{2.13}) and (\ref{2.31}),
$$
\aligned
p+V_{1,0}(|\tilde{z}|)
+\frac1p\widetilde{\omega}_1(|\tilde{z}|)
+\frac1{p^2}\widetilde{\omega}_2(|\tilde{z}|)
= p-2\log\left(1+|\tilde{z}|^2\right)+O\left(1\right)
\geq8\kappa\log p+
4\log\mu_i
+O\left(1\right)>7\kappa\log p.
\endaligned
$$
This, together with (\ref{2.3}) and (\ref{2.29}),   implies that
$0<U_\xi\leq2\sqrt{e}$ in $B_{1/p^{2\kappa}}(\xi_i)$,
and
$\sup_{B_{1/p^{2\kappa}}(\xi_i)}U_\xi\rightarrow\sqrt{e}$
as $p\rightarrow+\infty$.
Notice that by the maximum principle, it follows  that
for any $i=1,\ldots,m$, $G(x,\xi_i)>0$ and $G(x,q)>0$
over $\overline{\Omega}$, and further by (\ref{2.26}),
$U_\xi$ is a  positive, uniformly bounded function over $\overline{\Omega}$.
More precisely,
\begin{equation}\label{2.39}
\aligned
0<U_\xi(y)\leq2\sqrt{e},\,\,\,
\,\quad\forall\,\,\,y\in\overline{\Omega}.
\endaligned
\end{equation}

\vspace{1mm}
\vspace{1mm}

Consider that the scaling of solution to problem (\ref{1.1}) is as follows:
\begin{equation}\label{2.40}
\aligned
\upsilon(y)=\varepsilon^{2/(p-1)}u(\varepsilon y),\,\,\,
\,\quad\forall\,\,\,y\in\overline{\Omega}_{p},
\endaligned
\end{equation}
where $\Omega_{p}:=(e^{p/4})\Omega=(1/\varepsilon)\Omega$, then the function
 $\upsilon(y)$ satisfies
\begin{equation}\label{2.41}
\begin{array}{ll}
\left\{\aligned
&-\Delta_{a(\varepsilon y)}\upsilon+\varepsilon^2\upsilon=|\varepsilon y-q|^{2\alpha}\upsilon^p,
\,\,\,\,\,\,
\upsilon>0\,\,\,\,\,\,
\textrm{in}\,\,\,\,\,\Omega_p,\\[1mm]
&\frac{\partial \upsilon}{\partial\nu}=0
\qquad\qquad\qquad\qquad\qquad\qquad\qquad
\quad\,\,\textrm{on}\,\,\,\partial\Omega_p.
\endaligned\right.
\end{array}
\end{equation}
We write $q'=q/\varepsilon$ and
$\xi_i'=\xi_i/\varepsilon$, $i=1,\ldots,m$
and define the initial approximate solution of (\ref{2.41}) as
\begin{equation}\label{2.42}
\aligned
V_{\xi'}(y)=\varepsilon^{2/(p-1)}U_\xi(\varepsilon y),
\endaligned
\end{equation}
with
$\xi'=(\xi_1',\ldots,\xi_m')$
and
$U_\xi$  defined in (\ref{2.18}).
Let us set
$$
\aligned
S_p(\upsilon)=-\Delta_{a(\varepsilon y)}\upsilon+\varepsilon^2\upsilon-|\varepsilon y-q|^{2\alpha}\upsilon^p_{+},
\,\qquad\,\textrm{where}\,\,\,\,\upsilon_{+}=\max\{\upsilon,0\},
\endaligned
$$
and we consider the   functional
\begin{equation}\label{2.43}
\aligned
I_p(\upsilon)=\frac12\int_{\Omega_p}a(\varepsilon y)\left(
|\nabla \upsilon|^2+\varepsilon^2\upsilon^2
\right)dy-\frac1{p+1}\int_{\Omega_p}a(\varepsilon y)|\varepsilon y-q|^{2\alpha}\upsilon^{p+1}_{+}dy,
\qquad
\upsilon\in H^1(\Omega_p),
\endaligned
\end{equation}
whose nontrivial critical points are solutions of problem (\ref{2.41}).
Indeed, by virtue of the Hardy and Sobolev embedding inequalities
in \cite{GT, Y}
we have  that for all $\alpha>-1$,
the functional $I_p$ is well defined in $H^1(\Omega_p)$. Moreover,
by the maximum principle, it is easy to see that problem (\ref{2.41}) is equivalent to
$$
\aligned
S_p(\upsilon)=0,\,\quad\,\upsilon_{+}\not\equiv0
\,\quad\,\textrm{in}\,\,\,\,\Omega_p,\,\qquad\qquad\,
\frac{\partial\upsilon}{\partial\nu}=0
\,\quad\,\textrm{on}\,\,\,\partial\Omega_p.
\endaligned
$$
We will seek solutions of problem (\ref{2.41}) in the form
$\upsilon=V_{\xi'}+\phi$, where $\phi$ will represent a higher-order correction
in the expansion of $\upsilon$. Observe that
$$
\aligned
S_p(V_{\xi'}+\phi)=\mathcal{L}(\phi)+R_{\xi'}+N(\phi)=0,
\endaligned
$$
where
\begin{equation}\label{2.44}
\aligned
\mathcal{L}(\phi)=-\Delta_{a(\varepsilon y)}\phi+\varepsilon^2\phi-W_{\xi'}\phi\,\quad\,\textrm{with}\,\quad\,
W_{\xi'}=p|\varepsilon y-q|^{2\alpha}V_{\xi'}^{p-1},
\endaligned
\end{equation}
and
\begin{equation}\label{2.45}
\aligned
R_{\xi'}=-\Delta_{a(\varepsilon y)}V_{\xi'}+\varepsilon^2V_{\xi'}-|\varepsilon y-q|^{2\alpha}V_{\xi'}^p,
\quad\quad\quad
N(\phi)=-|\varepsilon y-q|^{2\alpha}\big[(V_{\xi'}+\phi)_+^p-V_{\xi'}^p-pV_{\xi'}^{p-1}\phi\big].
\endaligned
\end{equation}
In terms of $\phi$,
problem (\ref{2.41}) becomes
\begin{equation}\label{2.46}
\aligned
\left\{\aligned
&\mathcal{L}(\phi)=-\big[
R_{\xi'}+N(\phi)
\big]
\quad\textrm{in}\,\,\,\,\,\,\Omega_p,\\
&\frac{\partial \phi}{\partial\nu}=0
\quad\qquad\qquad\qquad\,\,\,\,
\textrm{on}\,\,\,\,
\partial\Omega_p.
\endaligned\right.\endaligned
\end{equation}

For any  $\xi=(\xi_1,\ldots,\xi_m)\in\mathcal{O}_{p}(q)$ and $h\in L^\infty(\Omega_p)$,
we  define a $L^\infty$-norm
$\big\|h\big\|_{*}:=\sup_{y\in\overline{\Omega}_p}\big|\mathbf{H}_{\xi'}(y)h(y)\big|$ with the weight
function
\begin{equation}\label{2.47}
\aligned
\mathbf{H}_{\xi'}(y)=\left[
\varepsilon^2+
\left(\frac{\varepsilon}{\rho_0v_0}\right)^2\frac{\big|\frac{\varepsilon y-q}{\rho_0v_0}\big|^{2\alpha}}
{\,\big(1+\big|\frac{\varepsilon y-q}{\rho_0v_0}\big|\big)^{4+2\hat{\alpha}+2\alpha}\,}
+
\sum\limits_{i=1}^m\frac{1}{\mu_i^2}\frac{1}{\big(1+\big|\frac{y-\xi'_i}{\mu_i}\big|\big)^{4+2\hat{\alpha}}}
\right]^{-1},
\endaligned
\end{equation}
where $\hat{\alpha}+1$ is a sufficiently small but fixed positive number,  independent of $p$,
such that $-1<\hat{\alpha}<\min\big\{\alpha,\,-2/3\big\}$.
With respect to the $\|\cdot\|_{*}$-norm,
we have the following.

\vspace{1mm}
\vspace{1mm}
\vspace{1mm}
\vspace{1mm}

\noindent{\bf Proposition 2.4.}\,\,{\it
Let $m$ be a non-negative  integer.
There exist constants $C>0$, $D_0>0$ and $p_m>1$ such that
\begin{equation}\label{2.48}
\aligned
\|R_{\xi'}\|_{*}\leq\frac{C}{p^4},
\endaligned
\end{equation}
and
\begin{equation}\label{2.49}
\aligned
W_{\xi'}(y)\leq D_0\left[
\left(\frac{\varepsilon}{\rho_0v_0}\right)^2
\left|\frac{\varepsilon y-q }{\rho_0v_0}\right|^{2\alpha}
e^{U_{1}\big(\frac{\varepsilon y-q }{\rho_0v_0}\big)}
+
\sum_{i=1}^{m}\frac{1}{\mu_i^2}
e^{V_{1,0}\big(\frac{y-\xi'_i}{\mu_i}\big)}\right],
\endaligned
\end{equation}
for any $\xi=(\xi_1,\ldots,\xi_m)\in\mathcal{O}_{p}(q)$ and any $p>p_m$.
Moreover,  if $|\varepsilon y-\xi_i|\leq\sqrt{\delta_i}/p^{2\kappa}$
for each $i=1,\ldots,m$,
\begin{equation}\label{2.51}
\aligned
W_{\xi'}(y)=\frac{1}{\mu_i^2}\frac{8}{\big(1+\big|\frac{y-\xi'_i}{\mu_i}\big|^2\big)^2}\left[1+
\frac1p\left(\widetilde{\omega}_1-V_{1,0}-\frac12V_{1,0}^2\right)\left(\frac{y-\xi'_i}{\mu_i}\right)
+O\left(\frac{\log^4\big(\big|\frac{y-\xi'_i}{\mu_i}\big|+2\big)}{p^2}\right)\right],
\endaligned
\end{equation}
while if $|\varepsilon y-q|\leq\sqrt{\delta_0}/p^{2\kappa}$,
\begin{equation}\label{2.50}
\aligned
W_{\xi'}(y)=
\left(\frac{\varepsilon}{\rho_0v_0}\right)^2
\frac{8(1+\alpha)^2\big|\frac{\varepsilon y-q }{\rho_0v_0}\big|^{2\alpha}}
{\big(1+\big|\frac{\varepsilon y-q }{\rho_0v_0}\big|^{2(1+\alpha)}\big)^2}
\left[1+
\frac1p\left(\omega_1-U_{1}-\frac{1}2U_{1}^2\right)\left(\frac{\varepsilon y-q }{\rho_0v_0}\right)
+O\left(\frac{\log^4\big(|\frac{\varepsilon y-q}{\rho_0v_0}|+2\big)}{p^2}\right)\right].
\endaligned
\end{equation}
}

\begin{proof}
From (\ref{2.18}), (\ref{2.19})  and (\ref{2.42}) it follows that
$$
\aligned
-\Delta_{a(\varepsilon y)}V_{\xi'}+\varepsilon^2V_{\xi'}
=\varepsilon^{2}\sum_{i=0}^m\varepsilon^{2/(p-1)}\left[-\Delta_{a}\big(U_i+H_i\big)+\big(U_i+H_i\big)\right]
=-\varepsilon^{2p/(p-1)}\sum_{i=0}^m\Delta U_i.
\endaligned
$$
Then by  (\ref{2.1}) and  (\ref{2.5})-(\ref{2.9}),
\begin{eqnarray}\label{2.52}
-\Delta_{a(\varepsilon y)}V_{\xi'}+\varepsilon^2V_{\xi'}
=\left(\frac{\varepsilon}{\rho_0v_0}\right)^2\frac{|z|^{2\alpha}e^{U_{1}(z)}}{p^{p/(p-1)}\mu_0^{2/(p-1)}}
\left(1
-\frac{1}{p}f_1
-\frac{1}{p^2}f_2
+\frac{1}{p}\omega_1
+\frac{1}{p^2}\omega_2
\right)(z)
\qquad\qquad\quad\,\,
&&\nonumber\\
+\sum_{i=1}^m\frac{e^{V_{1,0}(\tilde{z})}}{p^{p/(p-1)}\mu_i^{2p/(p-1)}\big|\xi_i-q\big|^{2\alpha/(p-1)}}
\left(1
-\frac{1}{p}\widetilde{f}_1
-\frac{1}{p^2}\widetilde{f}_2
+\frac{1}{p}\widetilde{\omega}_1
+\frac{1}{p^2}\widetilde{\omega}_2
\right)(\tilde{z})
&&
\end{eqnarray}
with  $z=(\varepsilon y-q)/\delta_0$ and $\tilde{z}=(\varepsilon y-\xi_i)/\delta_i$.
If $|\varepsilon y-q|=\delta_0|z|\geq1/p^{2\kappa}$
and $|\varepsilon y-\xi_i|=\delta_i|\tilde{z}|\geq1/p^{2\kappa}$ for each $i=1,\ldots,m$,
 by (\ref{2.1}),  (\ref{2.12}),  (\ref{2.13})  and (\ref{2.31})
we can compute that
$$
\aligned
U_{1}(z)=-p+O\left(\log p\right),\quad
V_{1,0}(\tilde{z})=-p+O\left(\log p\right),
\quad
\omega_{j}(z)=\frac{C_{j}p}{4(1+\alpha)}+O\left(\log p\right),
\quad
\widetilde{\omega}_{j}(\tilde{z})=\frac{\widetilde{C}_{j}p}{4}+O\left(\log p\right)
\endaligned
$$
with $j=1,2$,
and thus, by  (\ref{2.10}), (\ref{2.11})  and (\ref{2.52}),
\begin{equation}\label{2.53}
\aligned
-\Delta_{a(\varepsilon y)}V_{\xi'}+\varepsilon^2V_{\xi'}
=\left(\frac{\varepsilon}{\rho_0v_0}\right)^2\frac{|z|^{2\alpha}e^{U_{1}(z)}}{p^{p/(p-1)}\mu_0^{2/(p-1)}}
O\left(p^2
\right)
+
\sum_{i=1}^m\frac{e^{V_{1,0}(\tilde{z})}}{p^{p/(p-1)}\mu_i^{2p/(p-1)}\big|\xi_i-q\big|^{2\alpha/(p-1)}}
O\left(p^2
\right).
\endaligned
\end{equation}
In the same region, by (\ref{2.3}), (\ref{2.5}), (\ref{2.26}), (\ref{2.31}) and (\ref{2.42}) we obtain
\begin{equation}\label{2.54}
\aligned
|\varepsilon y-q|^{2\alpha}V_{\xi'}^p(y)
=O\left(\frac{|\varepsilon y-q|^{2\alpha}}p\left(\frac{\log p}{p}\right)^p
\right),
\endaligned
\end{equation}
which combined with  (\ref{2.45}), (\ref{2.47}) and (\ref{2.53}) easily yields that
\begin{eqnarray}\label{2.55}
\big|\mathbf{H}_{\xi'}(y)
R_{\xi'}(y)\big|
\leq Cp
\left[
\left|\frac{\varepsilon y-q}{\rho_0\upsilon_0}\right|^{2\hat{\alpha}-2\alpha}
+
\sum_{i=1}^m
\left|\frac{y-\xi'_i}{\mu_i}\right|^{2\hat{\alpha}}
+\frac{|\varepsilon y-q|^{2\alpha}}{p^2}\left(\frac{\sqrt{e}\log p}{p}\right)^p\right]
&&
\nonumber
\\[1mm]
=o\left(\max\left\{e^{p(\hat{\alpha}-\alpha)/4(1+\alpha)},\,e^{p\hat{\alpha}/4},\,e^{-p/4}\right\}
\right).
\qquad\qquad\qquad\qquad\quad\,\,\,\,\,
&&
\end{eqnarray}
On the other hand, if $|\varepsilon y-\xi_i|=\delta_i|\tilde{z}|<1/p^{2\kappa}$
for some $i\in\{1,\ldots,m\}$, then,
by (\ref{2.29}), (\ref{2.42}) and the relation
\begin{equation}\label{2.56}
\aligned
\left(
\frac{p\varepsilon^{2/(p-1)}}{\gamma\mu_i^{2/(p-1)}\big|\xi_i-q\big|^{2\alpha/(p-1)}}
\right)^p=\frac{1}{p^{p/(p-1)}\mu_i^{2p/(p-1)}\big|\xi_i-q\big|^{2\alpha p/(p-1)}},
\endaligned
\end{equation}
we find
\begin{eqnarray}\label{2.57}
|\varepsilon y-q|^{2\alpha}V^p_{\xi'}(y)=\frac{|\varepsilon y-q|^{2\alpha}}{p^{p/(p-1)}\mu_i^{2p/(p-1)}\big|\xi_i-q\big|^{2\alpha p/(p-1)}}
\left\{
1+\frac{1}pV_{1,0}(\tilde{z})+\frac{1}{p^2}\widetilde{\omega}_1(\tilde{z})
\right.
\,\qquad\,
&&\nonumber\\[0.5mm]
\left.
+\frac1{p^3}\left[\widetilde{\omega}_2(\tilde{z})+O\left(
p^2\delta_i^\beta|\tilde{z}|^\beta
+p^2\sum_{k=0}^m\delta_k^{\beta/2}\right)
\right]\right\}^p
\qquad\qquad\qquad\qquad
&&
\end{eqnarray}
with $0<\beta<\min\{1,\,2(1+\alpha)\}$.
From a Taylor expansion of the exponential and logarithmic functions
\begin{eqnarray}\label{2.58}
\left(1+\frac{a}p+\frac{b}{p^2}+\frac{c}{p^3}\right)^p
=e^a\left[1+\frac1p\left(b-\frac{a^2}2\right)+\frac1{p^2}\left(c
-a b
+\frac{a^3}3
+\frac{b^2}2-\frac{a^2b}2+\frac{a^4}8\right)
\right.
&&\nonumber\\[0.2mm]
\left.+O\left(\frac{\log^6(|\tilde{z}|+2)}{p^3}\right)\right],
\,\quad\qquad\qquad\qquad\qquad\qquad\qquad\qquad\,\,\,
&&
\end{eqnarray}
which holds for $|\tilde{z}|\leq Ce^{p/8}$
provided
$-4\log(|\tilde{z}|+2)\leq a(\tilde{z})\leq C$ and $|b(\tilde{z})|+|c(\tilde{z})|\leq C\log(|\tilde{z}|+2)$,
so we can
compute
that for $|\varepsilon y-\xi_i|=\delta_i|\tilde{z}|\leq\sqrt{\delta_i}/p^{2\kappa}$,
\begin{eqnarray}\label{2.59}
|\varepsilon y-q|^{2\alpha}V^p_{\xi'}(y)=\frac{|\varepsilon y-q|^{2\alpha}e^{V_{1,0}(\tilde{z})}}{p^{p/(p-1)}\mu_i^{2p/(p-1)}\big|\xi_i-q\big|^{2\alpha p/(p-1)}}
\left[\,
1+\frac1p\left(\widetilde{\omega}_1
-\frac12V_{1,0}^2
\right)(\tilde{z})+\frac1{p^2}\left(\widetilde{\omega}_2
-\widetilde{\omega}_1V_{1,0}
+\frac13V_{1,0}^3
\right.\right.
&&\nonumber\\[1mm]
\left.\left.
+\frac12\widetilde{\omega}_1^2
-\frac12\widetilde{\omega}_1V_{1,0}^2
+\frac18V_{1,0}^4
\right)(\tilde{z})
+O\left(\frac{\log^6(|\tilde{z}|+2)}{p^3}
+\delta_i^\beta|\tilde{z}|^\beta
+\sum_{k=0}^m\delta_k^{\beta/2}\right)
\right],
\qquad\qquad\qquad
&&
\end{eqnarray}
and then, by (\ref{2.11}), (\ref{2.45})  and (\ref{2.52}),
\begin{eqnarray}\label{2.60}
R_{\xi'}(y)=
\frac{e^{V_{1,0}(\tilde{z})}}{p^{p/(p-1)}\mu_i^{2p/(p-1)}\big|\xi_i-q\big|^{2\alpha/(p-1)}}
O\left(\frac{\log^6(|\tilde{z}|+2)}{p^3}
+\delta_i^\beta|\tilde{z}|^\beta
+\sum_{k=0}^m\delta_k^{\beta/4}\right).
\end{eqnarray}
Furthermore, in this region, by (\ref{2.47}) we obtain
\begin{eqnarray}\label{2.61}
\big|\mathbf{H}_{\xi'}(y)
R_{\xi'}(y)\big|
\leq\frac{C}{p^4}
\left(\left|\frac{y-\xi'_i}{\mu_i}\right|+1\right)^{2\hat{\alpha}}\log^6\left(\left|\frac{y-\xi'_i}{\mu_i}\right|+2\right)
=O\left(\frac1{p^4}
\right).
\end{eqnarray}
As in the  region
$\sqrt{\delta_i}/p^{2\kappa}<|\varepsilon y-\xi_i|=\delta_i|\tilde{z}|<1/p^{2\kappa}$,
by (\ref{2.11}) and (\ref{2.52}) we get
\begin{equation}\label{2.62}
\aligned
-\Delta_{a(\varepsilon y)}V_{\xi'}+\varepsilon^2V_{\xi'}=
\frac{e^{V_{1,0}(\tilde{z})}}{p^{p/(p-1)}\mu_i^{2p/(p-1)}\big|\xi_i-q\big|^{2\alpha/(p-1)}}
O\left(
p^2
\right),
\endaligned
\end{equation}
and by (\ref{2.57}),
\begin{equation}\label{2.63}
\aligned
|\varepsilon y-q|^{2\alpha}V^p_{\xi'}(y)
=\frac{|\varepsilon y-q|^{2\alpha}\,e^{V_{1,0}(\tilde{z})}}{p^{p/(p-1)}\mu_i^{2p/(p-1)}\big|\xi_i-q\big|^{2\alpha p/(p-1)}}
O\left(
1
\right),
\endaligned
\end{equation}
because  $(1+\frac{s}p)^p\leq e^s$.
Then in this region, by (\ref{2.3}), (\ref{2.31}), (\ref{2.45})   and (\ref{2.47}) we conclude
\begin{eqnarray}\label{2.64}
\big|\mathbf{H}_{\xi'}(y)
R_{\xi'}(y)\big|
\leq C p
\left|\frac{y-\xi'_i}{\mu_i}\right|^{2\hat{\alpha}}
=o\left(\frac{1}{p^4}
\right).
\end{eqnarray}
Similar to the above arguments in (\ref{2.59})-(\ref{2.60}) and (\ref{2.62})-(\ref{2.63}), we
have that if $|\varepsilon y-q|=\delta_0|z|\leq\sqrt{\delta}_0/p^{2\kappa}$,
\begin{eqnarray}\label{2.65}
R_{\xi'}(y)=
\left(\frac{\varepsilon}{\rho_0v_0}\right)^2\frac{|z|^{2\alpha}e^{U_{1}(z)}}{p^{p/(p-1)}\mu_0^{2/(p-1)}}
O\left(\frac{\log^6(|z|+2)}{p^3}
+\delta_0^\beta|z|^\beta
+\sum_{k=0}^m\delta_k^{\beta/4}\right),
\end{eqnarray}
while if $\sqrt{\delta}_0/p^{2\kappa}\leq|\varepsilon y-q|=\delta_0|z|\leq1/p^{2\kappa}$,
\begin{eqnarray}\label{2.66}
R_{\xi'}(y)=
\left(\frac{\varepsilon}{\rho_0v_0}\right)^2\frac{|z|^{2\alpha}e^{U_{1}(z)}}{p^{p/(p-1)}\mu_0^{2/(p-1)}}
O\left(p^2\right).
\end{eqnarray}
Thus in the region  $|\varepsilon y-q|=\delta_0|z|\leq1/p^{2\kappa}$, by (\ref{2.47}),
\begin{eqnarray}\label{2.67}
\big|\mathbf{H}_{\xi'}(y)
R_{\xi'}(y)\big|
\leq\frac{C}{p^4}
\left(\left|\frac{\varepsilon y -q}{\rho_0v_0}\right|+1\right)^{2\hat{\alpha}-2\alpha}\log^6\left(\left|\frac{\varepsilon y-q }{\rho_0v_0}\right|+2\right)
=O\left(\frac1{p^4}
\right).
\end{eqnarray}
As a consequence, putting    (\ref{2.55}), (\ref{2.61}), (\ref{2.64}) and   (\ref{2.67}) together,
we obtain estimate (\ref{2.48}).

Finally, it remains to prove the expansions (\ref{2.49})-(\ref{2.50}) for  $W_{\xi'}(y)=p|\varepsilon y-q|^{2\alpha}V_{\xi'}^{p-1}$.
If  $|\varepsilon y-\xi_i|=\delta_i|\tilde{z}|< 1/p^{2\kappa}$ for some $i\in\{1,\ldots,m\}$,
then by  (\ref{2.29}), (\ref{2.42}) and (\ref{2.56}),
\begin{eqnarray}\label{2.68}
W_{\xi'}(y)=p\left(
\frac{\,\,\varepsilon^{2/(p-1)}\big|\varepsilon y-q\big|^{2\alpha/(p-1)}}{\,\gamma\mu_i^{2/(p-1)}\big|\xi_i-q\big|^{2\alpha/(p-1)}\,}
\right)^{p-1}
\left[\,p+V_{1,0}(\tilde{z})
+\frac{\widetilde{\omega}_1(\tilde{z})}p
+\frac{\widetilde{\omega}_2(\tilde{z})}{p^2}
+O\left(
\delta_i^\beta|\tilde{z}|^\beta
+\sum_{k=0}^m\delta_k^{\beta/2}
\right)
\right]^{p-1}
&&\nonumber\\
=\frac{1}{\mu_i^2}
\left[
1+\frac{1}pV_{1,0}(\tilde{z})+\frac{1}{p^2}\widetilde{\omega}_1(\tilde{z})
+\frac{1}{p^3}\widetilde{\omega}_2(\tilde{z})+\frac{1}{p}O\left(
\delta_i^\beta|\tilde{z}|^\beta
+\sum_{k=0}^m\delta_k^{\beta/2}
\right)
\right]^{p-1},
\qquad\qquad\qquad\qquad\qquad
&&
\end{eqnarray}
where again we use the notation $\tilde{z}=(y-\xi'_i)/\mu_i$.
In this region,
we find
\begin{equation}\label{2.69}
\aligned
W_{\xi'}(y)\leq \frac{C}{\mu_i^2}
e^{V_{1,0}(\tilde{z})}
e^{-V_{1,0}(\tilde{z})/p}
=O\left(
\frac{1}{\mu_i^2}
e^{V_{1,0}(\tilde{z})}
\right),
\endaligned
\end{equation}
because $(1+a/p)^{p-1}\leq e^{a(p-1)/p}$
and $V_{1,0}(\tilde{z})\geq-p+O\big(\log p\big)$.
In particular, by a slight modification of formula (\ref{2.58}),
$$
\aligned
\left(1+\frac{a}p+\frac{b}{p^2}+\frac{c}{p^3}\right)^{p-1}=e^a\left[1+\frac1p\left(b-a-\frac{a^2}2\right)
+O\left(\frac{\log^4(|\tilde{z}|+2)}{p^2}\right)\right],
\endaligned
$$
we obtain that,  if $|\varepsilon y-\xi_i|=\delta_i|\tilde{z}|\leq\sqrt{\delta_i}/p^{2\kappa}$,
\begin{equation}\label{2.70}
\aligned
W_{\xi'}(y)
=\frac{1}{\mu_i^2}e^{V_{1,0}(\tilde{z})}\left[1+\frac1p\left(\widetilde{\omega}_1-V_{1,0}-\frac12V_{1,0}^2
\right)(\tilde{z})
+O\left(\frac{\log^4(|\tilde{z}|+2)}{p^2}\right)\right].
\endaligned
\end{equation}
If $|\varepsilon y-q|=\delta_0|z|\leq1/p^{2\kappa}$,
by (\ref{2.5}), (\ref{2.27}) and  (\ref{2.42}) we obtain
\begin{eqnarray}\label{2.71}
W_{\xi'}(y)=p|\varepsilon y-q|^{2\alpha}\left(
\frac{\varepsilon^{2/(p-1)}}{\gamma\mu_0^{2/(p-1)}}
\right)^{p-1}
\left[\,p+U_{1}(z)
+\frac{\omega_1(z)}p
+\frac{\omega_2(z)}{p^2}
+O\left(
\delta_0^\beta|z|^\beta
+\sum_{k=0}^m\delta_k^{\beta/2}
\right)
\right]^{p-1}
&&\nonumber\\
=\left(\frac{\varepsilon}{\rho_0v_0}\right)^2|z|^{2\alpha}
\left[
1+\frac{1}pU_{1}(z)+\frac{1}{p^2}\omega_1(z)
+\frac{1}{p^3}\omega_2(z)+\frac{1}{p}O\left(
\delta_0^\beta|z|^\beta
+\sum_{k=0}^m\delta_k^{\beta/2}
\right)
\right]^{p-1},
\qquad\,\,
&&
\end{eqnarray}
which, together with  the fact that
$U_{1}(z)\geq-p+O\left(\log p\right)$,
easily yields that in this region,
\begin{equation}\label{2.72}
\aligned
W_{\xi'}(y)\leq C\left(\frac{\varepsilon}{\rho_0v_0}\right)^2|z|^{2\alpha}
e^{U_1(z)}
e^{-U_1(z)/p}
=O\left(
\left(\frac{\varepsilon}{\rho_0v_0}\right)^2|z|^{2\alpha}
e^{U_1(z)}
\right).
\endaligned
\end{equation}
Similar to the  argument in (\ref{2.70}), by  (\ref{2.71}) we
can derive that if
$|\varepsilon y-q|=\delta_0|z|\leq\sqrt{\delta}_0/p^{2\kappa}$,
\begin{equation}\label{2.73}
\aligned
W_{\xi'}(y)=\left(\frac{\varepsilon}{\rho_0v_0}\right)^2|z|^{2\alpha}
e^{U_1(z)}\left[1+
\frac1p\left(\omega_1-U_{1}-\frac12U_{1}^2\right)(z)
+O\left(\frac{\log^4\big(|z|+2\big)}{p^2}\right)\right].
\endaligned
\end{equation}
Additionally, if
$|\varepsilon y-q|=\delta_0|z|>1/p^{2\kappa}$
and
$|\varepsilon y-\xi_i|=\delta_i|\tilde{z}|\geq1/p^{2\kappa}$  for
each $i=1,\ldots,m$, then, by (\ref{2.4}), (\ref{2.5}), (\ref{2.26}),
(\ref{2.31}) and (\ref{2.42}) we have that $W_{\xi'}(y)=O\big(\varepsilon^2p^{C}(\frac{\log p}{\gamma})^{p-1}\big)$.
This completes the proof.
\end{proof}

\vspace{1mm}
\vspace{1mm}

\noindent{\bf Remark 2.5.}\,\,\,As for $W_{\xi'}$, let us remark that
if $|\varepsilon y-q|\leq1/p^{2\kappa}$,
$$
\aligned
p|\varepsilon y-q|^{2\alpha}\left[V_{\xi'}(y)+O\left(\frac1{p^3}\right)\right]^{p-2}&\leq
Cp|\varepsilon y-q|^{2\alpha}\left(\frac{p\varepsilon^{2/(p-1)}}{\gamma\mu_0^{2/(p-1)}}
\right)^{p-2}e^{\frac{p-2}{p}U_{1}\big(\frac{\varepsilon y-q }{\rho_0\upsilon_0}\big)}\\
&=O\left(\left(\frac{\varepsilon}{\rho_0v_0}\right)^2\left|\frac{\varepsilon y-q }{\rho_0\upsilon_0}\right|^{2\alpha}
e^{U_{1}\big(\frac{\varepsilon y -q}{\rho_0\upsilon_0}\big)}\right),
\endaligned
$$
and if $|\varepsilon y-\xi_i|\leq1/p^{2\kappa}$ for some $i\in\{1,\ldots,m\}$,
$$
\aligned
p|\varepsilon y-q|^{2\alpha}\left[V_{\xi'}(y)+O\left(\frac1{p^3}\right)\right]^{p-2}
&\leq
Cp|\varepsilon y-q|^{2\alpha}\left(\frac{p\varepsilon^{2/(p-1)}}{\gamma\mu_i^{2/(p-1)}\big|\xi_i-q\big|^{2\alpha/(p-1)}}
\right)^{p-2}e^{\frac{p-2}{p}V_{1,0}\big(\frac{y-\xi'_i}{\mu_i}\big)}\\
&
=O\left(\frac{1}{\mu_i^2}e^{V_{1,0}\big(\frac{y-\xi'_i}{\mu_i}\big)}\right).
\endaligned
$$
Since these estimates are true if $|\varepsilon y-q|>1/p^{2\kappa}$
and $|\varepsilon y-\xi_i|>1/p^{2\kappa}$
for each $i=1,\ldots,m$,  we find
\begin{equation}\label{2.74}
\aligned
p|\varepsilon y-q|^{2\alpha}\left[V_{\xi'}(y)+O\left(\frac1{p^3}\right)\right]^{p-2}\leq
 C\left[
\left(\frac{\varepsilon}{\rho_0v_0}\right)^2
\left|\frac{\varepsilon y -q}{\rho_0v_0}\right|^{2\alpha}
e^{U_{1}\big(\frac{\varepsilon y-q }{\rho_0v_0}\big)}
+
\sum_{i=1}^{m}\frac{1}{\mu_i^2}
e^{V_{1,0}\big(\frac{y-\xi'_i}{\mu_i}\big)}\right].
\endaligned
\end{equation}

\vspace{1mm}
\vspace{1mm}

\section{The linearized problem and The nonlinear problem}
In this section we will first solve
the following linear problem: given $h\in C(\overline{\Omega}_p)$ and points
$\xi=(\xi_1,\ldots,\xi_m)\in\mathcal{O}_{p}(q)$, we find a function $\phi\in H^2(\Omega_p)$ and scalars
$c_{ij}\in\mathbb{R}$, $i=1,\ldots,m$, $j=1,J_i$,
such that
\begin{equation}\label{3.1}
\left\{\aligned
&\mathcal{L}(\phi)=-\Delta_{a(\varepsilon y)}\phi+\varepsilon^2\phi-W_{\xi'}\phi=h
+\frac1{a(\varepsilon y)}\sum\limits_{i=1}^m\sum\limits_{j=1}^{J_i}c_{ij}\chi_i\,Z_{ij}\,\,\ \,
\,\textrm{in}\,\,\,\,\,\,\Omega_p,\\
&\frac{\partial\phi}{\partial\nu}=0\,\,\,\,\,\,\,\,
\ \ \ \ \ \ \ \ \ \ \ \ \ \,\,
\qquad\qquad\qquad\quad\qquad\qquad\qquad\qquad
\ \,\ \ \ \,\,\,\ \,
\ \,\textrm{on}\,\,\,\,\partial\Omega_{p},\\[1mm]
&\int_{\Omega_p}\chi_i\,Z_{ij}\phi=0
\,\qquad\qquad\qquad\quad
\qquad\qquad\qquad\forall\,\,i=1,\ldots,m,\,\,\,j=1, J_i,
\endaligned\right.
\end{equation}
where $W_{\xi'}=p|\varepsilon y-q|^{2\alpha}V_{\xi'}^{p-1}$
satisfies (\ref{2.49})-(\ref{2.50}),
$J_i=1$ if $i=l+1,\ldots,m$ while
$J_i=2$ if $i=1,\ldots,l$,
and $Z_{ij}$, $\chi_i$ are defined as follows.
Let
\begin{equation}\label{3.2}
\aligned
\mathcal{Z}_{q}(z)=\frac{|z|^{2(1+\alpha)}-1}{|z|^{2(1+\alpha)}+1},
\,\,\quad\quad\quad\,\,
\mathcal{Z}_{0}(z)=\frac{|z|^2-1}{|z|^2+1},
\,\,\quad\quad\quad\,\,
\mathcal{Z}_{j}(z)=\frac{z_j}{|z|^2+1},\,\,\,\,j=1,\,2.
\endaligned
\end{equation}
It is well known (see \cite{BP,CL,CY,E,EPW}) that
\begin{itemize}
  \item any bounded
solution to
\begin{equation}\label{3.3}
\aligned
\Delta
\phi+\frac{8(1+\alpha)^2|z|^{2\alpha}}{(1+|z|^{2(1+\alpha)})^2}\phi=0
\,\quad\textrm{in}\,\,\,\mathbb{R}^2,
\endaligned
\end{equation}
where $-1<\alpha\not\in\mathbb{N}$,
is proportional to $\mathcal{Z}_q$;
  \item any bounded
solution to
\begin{equation}\label{3.22}
\aligned
\Delta
\phi+\frac{8(1+\alpha)^2|z|^{2\alpha}}{(1+|z|^{2(1+\alpha)})^2}\phi=0
\,\quad\textrm{in}\,\,\,\mathbb{R}^2_{+},
\,\qquad\qquad\,\frac{\partial\phi}{\partial\nu}=0
\,\quad\textrm{on}\,\,\,\partial\mathbb{R}^2_{+},
\endaligned
\end{equation}
where $-1<\alpha\not\in\mathbb{N}$ and $\mathbb{R}^2_{+}:=\{(z_1,z_2):\,z_2>0\}$,
is proportional to $\mathcal{Z}_q$;
  \item any bounded
solution to
\begin{equation}\label{3.4}
\aligned
\Delta
\phi+\frac{8}{(1+|z|^2)^2}\phi=0
\,\quad\textrm{in}\,\,\,\mathbb{R}^2,
\endaligned
\end{equation}
is a linear combination of $\mathcal{Z}_j$, $j=0,1,2$;
  \item any bounded
solution to
\begin{equation}\label{3.5}
\aligned
\Delta
\phi+\frac{8}{(1+|z|^2)^2}\phi=0
\,\quad\textrm{in}\,\,\,\mathbb{R}^2_{+},
\,\qquad\qquad\qquad\,\frac{\partial\phi}{\partial\nu}=0
\,\quad\textrm{on}\,\,\,\partial\mathbb{R}^2_{+},
\endaligned
\end{equation}
is a linear combination of $\mathcal{Z}_j$, $j=0,1$.
\end{itemize}

Next, let us consider a large but fixed positive  number  $R_0$ and
smooth non-increasing cut-off function
$\chi:\,\mathbb{R}\rightarrow[0,1]$  such that
$\chi(r)=1$ for $r\leq R_0$, and $\chi(r)=0$ for $r\geq R_0+1$.

For the interior spike case, i.e. $q\in\Omega$ and $\xi_i\in\Omega$ with $i=1,\ldots,l$,
we define
\begin{equation}\label{3.6}
\aligned
\chi_q(y)=\chi\left(
\frac{|\varepsilon y-q|}{\rho_0 v_0}
\right),
\,\qquad\qquad\,
Z_{q}(y)=\frac{\varepsilon}{\rho_0 v_0}\mathcal{Z}_q\left(\frac{\varepsilon y-q}{\rho_0 v_0}\right),
\endaligned
\end{equation}
and for any $i=1,\ldots,l$ and $j=0,1,2$,
\begin{equation}\label{3.7}
\aligned
\chi_i(y)=\chi\left(
\frac{|y-\xi'_i|}{\mu_i}
\right),
\,\qquad\qquad\,
Z_{ij}(y)=\frac{1}{\mu_i}\mathcal{Z}_j\left(\frac{y-\xi_i'}{\mu_i}\right).
\endaligned
\end{equation}

For the boundary spike case, i.e.
$q\in\partial\Omega$ and $\xi_i\in\partial\Omega$ with
$i=l+1,\ldots,m$, we need first to straighten
the boundary.
Namely, at each boundary points $q$ and $\xi_i$,
$i=l+1,\ldots,m$,
we
define the planar rotation maps
$A_q: \mathbb{R}^2\rightarrow\mathbb{R}^2$
and
$A_{i}: \mathbb{R}^2\rightarrow\mathbb{R}^2$
such that
$A_q\nu_{\Omega}(q)=\nu_{\mathbb{R}_+^2}(0)$
and $A_{i}\nu_{\Omega}(\xi_i)=\nu_{\mathbb{R}_+^2}(0)$.
Let $\mathcal{G}(x_1)$ be the defining function
for the boundary $A_q(\po-\{q\})$ or  $A_i(\po-\{\xi_i\})$
in a small
neighborhood $B_\delta(0)$ of the origin,
that is, there exist $R_1>0$, $\delta>0$ small and
a smooth function $\mathcal{G}:(-R_1,R_1)\mapsto\mathbb{R}$
satisfying $\mathcal{G}(0)=0$, $\mathcal{G}'(0)=0$ and such that
$A_q(\Omega-\{q\})\cap B_\delta(0)$
or $A_i(\Omega-\{\xi_i\})\cap B_\delta(0)$ can be rewritten in the form
$\{(x_1,x_2):\,-R_1<x_1<R_1,\,x_2>\mathcal{G}(x_1)\}\cap B_\delta(0)$.
Furthermore, we consider the flattening changes of variables
$F_q: \overline{A_q(\Omega-\{q\})}\cap B_\delta(0)\mapsto\mathbb{R}^2_{+}$
and
$F_i: \overline{A_i(\Omega-\{\xi_i\})}\cap B_\delta(0)\mapsto\mathbb{R}^2_{+}$,
respectively
defined by
\begin{equation}\label{3.8}
\aligned
F_q=(F_{q1}, F_{q2}),
\ \qquad
\textrm{where}
\ \qquad\,
F_{q1}=x_1+\frac{x_2-\mathcal{G}(x_1)}{\,1+|\mathcal{G}'(x_1)|^2\,}\mathcal{G}'(x_1),\ \qquad\,F_{q2}=x_2-\mathcal{G}(x_1),
\endaligned
\end{equation}
and for any $i=l+1,\ldots,m$,
\begin{equation}\label{3.9}
\aligned
F_i=(F_{i1}, F_{i2}),
\ \qquad
\textrm{where}
\ \qquad\,
F_{i1}=x_1+\frac{x_2-\mathcal{G}(x_1)}{\,1+|\mathcal{G}'(x_1)|^2\,}\mathcal{G}'(x_1),\ \qquad\,F_{i2}=x_2-\mathcal{G}(x_1).
\endaligned
\end{equation}
Let us denote
\begin{equation}\label{3.49}
\aligned
F_q^p(y)=e^{p/4}F_q\big(A_q(e^{-p/4}y-q)\big)=\frac{1}{\varepsilon}F_q\big(A_q(\varepsilon y-q)\big),
\endaligned
\end{equation}
\begin{equation}\label{3.10}
\aligned
F_i^p(y)=e^{p/4}F_i\big(A_i(e^{-p/4}y-\xi_i)\big)=\frac{1}{\varepsilon}F_i\big(A_i(\varepsilon y-\xi_i)\big),
\,\,\qquad\,\,i=l+1,\ldots,m.
\endaligned
\end{equation}
Here, recalling that $q\in\po$,  we  define
\begin{equation}\label{3.11}
\aligned
\chi_q(y)=\chi\left(
\frac{\varepsilon}{\rho_0 v_0}\big|F_q^p(y)\big|
\right),
\,\qquad\qquad\,
Z_{q}(y)=\frac{\varepsilon}{\rho_0 v_0}\mathcal{Z}_q
\left(\frac{\varepsilon}{\rho_0 v_0}F_q^p(y)\right),
\endaligned
\end{equation}
and for any $i=l+1,\ldots,m$ and $j=0,1$,
\begin{equation}\label{3.12}
\aligned
\chi_i(y)=\chi\left(
\frac{1}{\mu_i}\big|F_i^p(y)\big|
\right),
\,\qquad\quad\qquad\,
Z_{ij}(y)=\frac{1}{\mu_i}\mathcal{Z}_{j}\left(\frac{1}{\mu_i}F_i^p(y)\right).
\endaligned
\end{equation}
It is important  to note that if $q\in\po$, then
the maps $F^p_q$ and $F^p_i$, $i=l+1,\ldots,m$
preserve the  Neumann boundary condition.

\vspace{1mm}
\vspace{1mm}
\vspace{1mm}
\vspace{1mm}

\noindent {\bf Proposition 3.1.}\,\,{\it
Let $q\in\oo$ and
$m$ be a non-negative integer.
Then there exist constants $C>0$  and $p_m>1$ such
that for any  $p>p_m$,    any points
$\xi=(\xi_1,\ldots,\xi_m)\in\mathcal{O}_{p}(q)$ and any $h\in C(\overline{\Omega}_p)$,
there is a unique solution $\phi\in
H^2(\Omega_p)$  and scalars $c_{ij}\in\mathbb{R}$,  $i=1,\ldots,m$, $j=1,J_i$ to
problem {\upshape(\ref{3.1})},  which satisfy
\begin{equation}\label{3.15}
\aligned
\|\phi\|_{L^{\infty}(\Omega_p)}\leq Cp\|h\|_{*}
\qquad\qquad
\textrm{and}\,\,
\qquad\qquad
\sum_{i=1}^m\sum_{j=1}^{J_i}\mu_i|c_{ij}|\leq C\|h\|_{*}.
\endaligned
\end{equation}
}\indent We carry out the proof in the following four steps.


{\bf Step 1:} Constructing a suitable barrier.

\vspace{1mm}
\vspace{1mm}
\vspace{1mm}

\noindent{\bf Lemma 3.2.}\,\,{\it There exist constants $R_1>0$ and $C>0$,
independent of $p$, such that
for any sufficiently large
$p$,  any points $\xi=(\xi_1,\ldots,\xi_m)\in\mathcal{O}_p(q)$ and any $\hat{\alpha}\in(-1,\,\,\min\{\alpha,\,-2/3 \})$,
 there exists
$$
\aligned
\psi:\,\,\overline{\Omega}_p\setminus\left(\bigcup_{i=1}^mB_{R_1\mu_i}(\xi'_i)
\cup B_{\frac{R_1\rho_0\upsilon_0}{\varepsilon}}(q')
\right)\longmapsto\mathbb{R}
\endaligned
$$
smooth  and
positive so that
$$
\aligned
\mathcal{L}(\psi)\geq
\left(\frac{\varepsilon}{\rho_0v_0}\right)^2\frac{1}
{\,\big|\frac{\varepsilon y-q}{\rho_0v_0}\big|^{4+2\hat{\alpha}}\,}
&+
\sum\limits_{i=1}^m\frac{1}{\mu_i^2}\frac{1}{\big|\frac{y-\xi'_i}{\mu_i}\big|^{4+2\hat{\alpha}}}+\varepsilon^2
\,\,\,\quad\,
\textrm{in}\,\,\,\,\,\,\,\Omega_p\setminus\left(\bigcup_{i=1}^mB_{R_1\mu_i}(\xi'_i)
\cup B_{\frac{R_1\rho_0\upsilon_0}{\varepsilon}}(q')
\right),\\
\frac{\partial\psi}{\partial\nu}&\geq0
\,\,\,\qquad\qquad\qquad\qquad\quad\,\,
\quad
\textrm{on}\,\,\,\,\partial\Omega_p\setminus\left(\bigcup_{i=1}^mB_{R_1\mu_i}(\xi'_i)
\cup B_{\frac{R_1\rho_0\upsilon_0}{\varepsilon}}(q')
\right),\\
\psi&>0
\,\,\,\qquad\qquad\qquad\qquad\quad\,\,
\quad
\textrm{in}\,\,\,\,\,\,\,
\Omega_p\setminus\left(\bigcup_{i=1}^mB_{R_1\mu_i}(\xi'_i)
\cup B_{\frac{R_1\rho_0\upsilon_0}{\varepsilon}}(q')
\right),\\
\psi&\geq1
\,\,\,\qquad\qquad\qquad\qquad\quad\,\,
\quad
\textrm{on}\,\,\,\,
\Omega_p\cap\left(\bigcup_{i=1}^m\partial B_{R_1\mu_i}(\xi'_i)
\cup \partial B_{\frac{R_1\rho_0\upsilon_0}{\varepsilon}}(q')
\right).
\endaligned
$$
Moreover, $\psi$ is    uniformly bounded, i.e.
$$
\aligned
1<\psi\leq C\,\,\,\quad\textrm{in}\,\,\,
\,\overline{\Omega}_p\setminus\left(\bigcup_{i=1}^mB_{R_1\mu_i}(\xi'_i)
\cup B_{\frac{R_1\rho_0\upsilon_0}{\varepsilon}}(q')
\right).
\endaligned
$$
}
\begin{proof}
Let $\Psi_0(y)=\psi_0(\varepsilon y)$, where $\psi_0$ is
the solution to
$$
\aligned
\left\{
\aligned
&-\Delta_a\psi_0+\psi_0=1
\quad
\textrm{in}\,\,\,\,\,\Omega,\\
&\frac{\partial\psi_0}{\partial\nu}=1
\,\,\quad\quad\quad\quad\,\,
\,\textrm{on}
\,\,\,\partial\Omega.
\endaligned\right.
\endaligned
$$
Then
$$
\aligned
-\Delta_{a(\varepsilon y)}\Psi_0+\varepsilon^2\Psi_0=\varepsilon^2
\,\,\,\,\,
\textrm{in}\,\,\,\,\,\Omega_p,
\,\qquad\textrm{and}\quad\qquad\frac{\partial\Psi_0}{\partial\nu}=\varepsilon
\,\quad\textrm{on}
\,\,\,\,\partial\Omega_p.
\endaligned
$$
Obviously, $\Psi_0$ is a positive, uniformly bounded function over $\overline{\Omega}_p$.
Take the function
$$
\aligned
\psi=\left(1-\frac{1}
{\,\big|\frac{\varepsilon y-q}{\rho_0v_0}\big|^{2(1+\hat{\alpha})}\,}\right)
+\sum_{i=1}^m
\left(1-
\frac{1}{\,\big|\frac{y-\xi'_i}{\mu_i}\big|^{2(1+\hat{\alpha})}\,}
\right)+
C_1
\Psi_0(y).
\endaligned
$$
As a result, it is directly checked that,
choosing the positive constant $C_1$ larger if necessary, $\psi$
meets the required
conditions
of the lemma for numbers $R_1$ and $p$ large enough.
\end{proof}

\vspace{1mm}
\vspace{1mm}
\vspace{1mm}

{\bf Step 2:} Handing a linear equation. Given $h\in C^{0,\alpha}(\overline{\Omega}_p)$
and $\xi=(\xi_1,\ldots,\xi_m)\in\mathcal{O}_p(q)$, we first study the linear equation
\begin{equation}\label{3.16}
\left\{\aligned
&\mathcal{L}(\phi)=-\Delta_{a(\varepsilon y)}\phi+\varepsilon^2\phi-W_{\xi'}\phi=h
\,\,\ \,
\,\textrm{in}\,\,\,\,\,\,\Omega_p,\\
&\frac{\partial\phi}{\partial\nu}=0\,\,\,\,\,
\qquad\qquad\qquad\qquad\qquad\qquad
\ \,\textrm{on}\,\,\,\,\partial\Omega_{p}.
\endaligned\right.
\end{equation}
For the solution of (\ref{3.16}) satisfying more
orthogonality
conditions than those in (\ref{3.1}),
we prove the following a priori estimate.

\vspace{1mm}
\vspace{1mm}
\vspace{1mm}

\noindent{\bf Lemma 3.3.}\,\,{\it There exist $R_0>0$ and  $p_m>1$ such that for any $p>p_m$
and any solution $\phi$ of {\upshape (\ref{3.16})} with the orthogonality conditions
\begin{equation}\label{3.17}
\aligned
\int_{\Omega_p}\chi_qZ_q\phi=0
\,\qquad\,
\textrm{and}
\,\qquad\,
\int_{\Omega_p}\chi_iZ_{ij}\phi=0,\,\,\,\,\,
\,\,\,\,\,\,i=1,\ldots,m,\,\,j=0,1,J_i,
\endaligned
\end{equation}
we have
\begin{equation}\label{3.18}
\aligned
\|\phi\|_{L^{\infty}(\Omega_p)}\leq C
\|h\|_{*},
\endaligned
\end{equation}
where $C>0$ is independent of $p$.
}

\vspace{1mm}
\vspace{1mm}

\begin{proof}
Take $R_0=2R_1$ with
$R_1$ as the constant in the previous step.
Since  $\xi=(\xi_1,\ldots,\xi_m)\in\mathcal{O}_p(q)$,
$\rho_0 v_0=o(1/p^\kappa)$ and $\varepsilon\mu_i=o(1/p^\kappa)$
for $p$ large enough, we find $B_{R_1\rho_0 v_0/\varepsilon }(q')$ and
$B_{R_1\mu_i}(\xi'_i)$, $i=1,\ldots,m$ disjointed.
Let $h$ be bounded and $\phi$ be a bounded solution to (\ref{3.16}) satisfying (\ref{3.17}).
Let us consider the ``inner norm''
\begin{equation}\label{3.19}
\aligned
\|\phi\|_{**}=\sup_{y\in\overline{\Omega_p\cap\left(
\bigcup_{i=1}^mB_{R_1\mu_i}(\xi'_i)\cup B_{R_1\rho_0 v_0/\varepsilon}(q')\right)}}
|\phi(y)|,
\endaligned
\end{equation}
and  claim that there is a constant $C>0$ independent of  $p$ such that
\begin{equation}\label{3.20}
\aligned
\|\phi\|_{L^{\infty}(\Omega_p)}\leq C\left(\|\phi\|_{**}+
\|h\|_{*}\right).
\endaligned
\end{equation}
We will establish this estimate with the use of the barrier
$\psi$ constructed by Lemma 3.2.
In fact,  we take
$$
\aligned
\widetilde{\phi}(y)=C_1\left(\|\phi\|_{**}+
\|h\|_{*}
\right)\psi(y)
\,\qquad\forall\,\,\,y\in\overline{\Omega}_p\setminus\left(\bigcup_{i=1}^mB_{R_1\mu_i}(\xi'_i)
\cup B_{\frac{R_1\rho_0\upsilon_0}{\varepsilon}}(q')
\right),
\endaligned
$$
where
$C_1>0$ is a large   constant, independent of $p$.
Then for
$y\in\Omega_p\setminus\left(\bigcup_{i=1}^mB_{R_1\mu_i}(\xi'_i)
\cup B_{{R_1\rho_0\upsilon_0}/{\varepsilon}}(q')
\right)$,
$$
\aligned
\mathcal{L}(\widetilde{\phi}\pm\phi)(y)\geq C_{1}\,\|h\|_{*}\left\{
\left(\frac{\varepsilon}{\rho_0v_0}\right)^2\frac{1}
{\,\big|\frac{\varepsilon y-q}{\rho_0v_0}\big|^{4+2\hat{\alpha}}\,}
+
\sum\limits_{i=1}^m\frac{1}{\mu_i^2}\frac{1}{\big|\frac{y-\xi'_i}{\mu_i}\big|^{4+2\hat{\alpha}}}+\varepsilon^2
\right\}\pm
h(y)\geq|h(y)|\pm h(y)\geq0,
\endaligned
$$
for
$y\in\partial\Omega_p\setminus\left(\bigcup_{i=1}^mB_{R_1\mu_i}(\xi'_i)
\cup B_{{R_1\rho_0\upsilon_0}/{\varepsilon}}(q')
\right)$,
$$
\aligned
\frac{\partial}{\partial\nu}(\widetilde{\phi}\pm\phi)(y)\geq 0,
\endaligned
$$
and for
$y\in\Omega_p\cap\left(\bigcup_{i=1}^m\partial B_{R_1\mu_i}(\xi'_i)
\cup\partial B_{{R_1\rho_0\upsilon_0}/{\varepsilon}}(q')\right)$,
$$
\aligned
(\widetilde{\phi}\pm\phi)(y)>\|\phi\|_{**}\pm\phi(y)\geq
|\phi(y)|\pm\phi(y)\geq 0.
\endaligned
$$
 By the maximum
principle it follows   that
$-\widetilde{\phi}\leq\phi\leq\widetilde{\phi}$ on
$\overline{\Omega}_p\setminus\left(\bigcup_{i=1}^mB_{R_1\mu_i}(\xi'_i)
\cup B_{{R_1\rho_0\upsilon_0}/{\varepsilon}}(q')
\right)$, which easily
 implies that estimate (\ref{3.20}) holds.

We prove the lemma by contradiction. Assume that there are  sequences of
parameters $p_n\rightarrow+\infty$,
points $\xi^n=(\xi_1^n,\ldots,\xi_m^n)\in\mathcal{O}_{p_n}(q)$,
functions $W_{(\xi^n)'}$, $h_n$ and associated solutions $\phi_n$ of
equation (\ref{3.16}) with orthogonality conditions (\ref{3.17})
such that
\begin{equation}\label{3.21}
\aligned
\|\phi_n\|_{L^{\infty}(\Omega_{p_n})}=1
\,\,\qquad\,\,
\textrm{but}
\,\,\qquad\,\,\|h_n\|_{*}\rightarrow0
\,\,\quad\,\,\textrm{as}\,\,\,\,n\rightarrow+\infty.
\endaligned
\end{equation}
For $q\in\Omega$, we consider
$\widehat{\phi}^n_q(z)=\phi_n\big((\rho_0^nv_0^nz+q)/\varepsilon_n\big)$
and
$\widehat{h}^n_q(z)=h_n\big((\rho_0^nv_0^nz+q)/\varepsilon_n\big)$,
where  $\mu^n=\big(\mu^n_0,\mu^n_1,\ldots,\mu_m^n\big)$,
$\varepsilon_n=\exp\left\{-\frac14p_n\right\}$,
$\rho^n_0=(\varepsilon_n)^{\frac{1}{1+\alpha}}$ and
$v^n_0=(\mu_0^n)^{\frac{1}{1+\alpha}}$.
Observe that
$$
\aligned
\big(-\Delta_{a(\varepsilon_n y)}
\phi_n+\varepsilon_n^2\phi-W_{(\xi^n)'}\phi_n\big)\left|_{y=\large\frac{\rho_0^nv_0^nz+q}{\varepsilon_n}}
=\left(\frac{\varepsilon_n}{\rho_0^nv_0^n}\right)^2\left[
-\Delta_{\widehat{a}_n(z)}\widehat{\phi}_q^n
+\big(\rho_0^nv_0^n\big)^{2}\widehat{\phi}_q^n
-\left(\frac{\rho_0^nv_0^n}{\varepsilon_n}\right)^2\widehat{W}^n\widehat{\phi}_q^n
\right](z),
\right.
\endaligned
$$
where
$$
\aligned
\widehat{a}_n(z)=a(\rho_0^nv_0^nz+q),
\qquad\qquad
\widehat{W}^n(z)=W_{(\xi^n)'}\big((\rho_0^nv_0^nz+q)/\varepsilon_n\big).
\endaligned
$$
From  the expansion of $W_{(\xi^n)'}$ in (\ref{2.50}),
it follows that $\widehat{\phi}_q^n(z)$ solves
$$
\aligned
-\Delta_{\widehat{a}_n(z)}\widehat{\phi}_q^n
+\big(\rho_0^nv_0^n\big)^{2}\widehat{\phi}_q^n-\frac{8(1+\alpha)^2\big|z\big|^{2\alpha}}{\,\big(1+\big|z\big|^{2(1+\alpha)}\big)^2\,}
\left[1+O\left(\frac1p\right)
\right]
\widehat{\phi}^n_q(z)=\left(\frac{\rho^n_0v^n_0}{\varepsilon_n}\right)^2\widehat{h}^n_q(z)
\endaligned
$$
for any $z\in B_{R_0+2}(0)$.
Owing to
(\ref{3.21}) and the definition of the $\|\cdot\|_{*}$-norm
with respect to (\ref{2.47}), we find that
 for any $\theta\in\big(1, -1/\hat{\alpha}\big)$,
$\big(\frac{\rho^n_0v^n_0}{\varepsilon_n}\big)^2\widehat{h}^n_q\rightarrow 0$ in
$L^{\theta}\big(B_{R_0+2}(0)\big)$.
Since
$\frac{8(1+\alpha)^2|z|^{2\alpha}}{(1+|z|^{2(1+\alpha)})^2}$
is bounded in $L^{\theta}\big(B_{R_0+2}(0)\big)$,
standard elliptic regularity   implies that
$\widehat{\phi}^n_q$
converges uniformly over
compact subsets near the origin to a bounded solution $\widehat{\phi}^{\infty}_q$ of equation
$(\ref{3.3})$ with the property
\begin{equation}\label{3.23}
\aligned
\int_{\mathbb{R}^2}\chi \mathcal{Z}_q\widehat{\phi}_q^{\infty}=0.
\endaligned
\end{equation}
Then $\widehat{\phi}^{\infty}_q$
is proportional to $\mathcal{Z}_q$.
Since $\int_{\mathbb{R}^2}\chi \mathcal{Z}_q^2>0$, by
(\ref{3.23}) we have that $\widehat{\phi}^{\infty}_q\equiv0$ in $B_{R_1}(0)$.

For $q\in\po$, we consider
$\widehat{\phi}^n_q(z)=\phi_n\big((A^n_q)^{-1}\varepsilon_n^{-1}\rho_0^nv_0^nz+q_n'\big)$,
where $q_n'=q/\varepsilon_n$ and
$A_q^n: \mathbb{R}^2\rightarrow\mathbb{R}^2$ is a rotation map
such that $A_q^n\nu_{\Omega_{p_n}}\big(q_n'\big)=\nu_{\mathbb{R}_+^2}\big(0\big)$.
Similarly to the above argument, by using the expansion
of $W_{(\xi^n)'}$ in (\ref{2.50}) and elliptic regularity we
can derive that
$\widehat{\phi}^n_q$
converges uniformly over compact sets
to a bounded solution $\widehat{\phi}^{\infty}_q$ of equation
$(\ref{3.22})$ with the property
\begin{equation}\label{3.24}
\aligned
\int_{\mathbb{R}^2_{+}}\chi \mathcal{Z}_q\widehat{\phi}_q^{\infty}=0.
\endaligned
\end{equation}
Then $\widehat{\phi}^{\infty}_q$
is proportional to $\mathcal{Z}_q$.
Since $\int_{\mathbb{R}^2_{+}}\chi \mathcal{Z}_q^2>0$, by
(\ref{3.24}) we have that $\widehat{\phi}^{\infty}_q\equiv0$ in $B_{R_1}^{+}(0)$.

For each $i\in\{1,\ldots,l\}$, we have that $\xi_i^n\in\Omega$ and we consider
$\widehat{\phi}^n_i(z)=\phi_n\big(\mu_i^nz+(\xi^n_i)'\big)$
with $(\xi^n_i)'=\xi^n_i/\varepsilon_n$.
Notice that
$$
\aligned
h_n(y)=
\big(-\Delta_{a(\varepsilon_n y)}\phi_n+\varepsilon_n^2\phi
-W_{(\xi^n)'}\phi_n\big)\big|_{y=\mu_{i}^n z+(\xi^n_i)'}
=(\mu_i^n)^{-2}\left[
-\Delta_{\widehat{a}_n(z)}\widehat{\phi}_i^n
+\varepsilon_n^2(\mu_i^n)^{2}\widehat{\phi}_i^n
-(\mu_i^n)^{2}\widehat{W}^n\widehat{\phi}_i^n
\right](z),
\endaligned
$$
where
$$
\aligned
\widehat{a}_n(z)=a(\varepsilon_n\mu_{i}^nz+\xi^n_i),
\qquad\qquad
\widehat{W}^n(z)=W_{(\xi^n)'}\big(\mu_{i}^nz+(\xi^n_i)'\big).
\endaligned
$$
Using the expansion of $W_{(\xi^n)'}$ in (\ref{2.51})
and elliptic regularity, we find that
$\widehat{\phi}^n_i$
converges uniformly over
compact subsets  near the origin to a bounded solution $\widehat{\phi}^{\infty}_i$ of equation
$(\ref{3.4})$, which satisfies
\begin{equation}\label{3.25}
\aligned
\int_{\mathbb{R}^2}\chi \mathcal{Z}_j\widehat{\phi}_i^{\infty}=0
\quad\,\,\,\textrm{for}\,\,\,\,j=0,\,1,\,2.
\endaligned
\end{equation}
Thus $\widehat{\phi}^{\infty}_i$ is
a linear combination of $\mathcal{Z}_j$, $j=0,1,2$.
But $\int_{\mathbb{R}^2}\chi \mathcal{Z}_j^2>0$
and $\int_{\mathbb{R}^2}\chi \mathcal{Z}_j\mathcal{Z}_{t}=0$ for $j\neq t$.
Hence (\ref{3.25}) implies  that $\widehat{\phi}_i^{\infty}\equiv0$
in $B_{R_1}(0)$.

Finally,  for each $i\in\{l+1,\ldots,m\}$, we have that
$\xi_i^n\in\partial\Omega$ and we consider
$\widehat{\phi}^n_i(z)=\phi_n\big((A^n_i)^{-1}\mu_i^nz+(\xi^n_i)'\big)$,
where $A_i^n: \mathbb{R}^2\rightarrow\mathbb{R}^2$ is a rotation map
such that $A_i^n\nu_{\Omega_{p_n}}\big((\xi_i^n)'\big)=\nu_{\mathbb{R}_+^2}\big(0\big)$.
Similar to the above argument, by using
the expansion of $W_{(\xi^n)'}$ in (\ref{2.51})
and elliptic regularity
we can  derive that
$\widehat{\phi}^n_i$
converges uniformly over compact sets
to a bounded solution $\widehat{\phi}^{\infty}_i$ of equation
$(\ref{3.5})$, which satisfies
\begin{equation}\label{3.26}
\aligned
\int_{\mathbb{R}_{+}^2}\chi \mathcal{Z}_j\widehat{\phi}_i^{\infty}=0
\quad\,\,\,\textrm{for}\,\,\,\,j=0,\,1.
\endaligned
\end{equation}
Thus $\widehat{\phi}^{\infty}_i$ is
a linear combination of $\mathcal{Z}_j$, $j=0,1$.
Note that $\int_{\mathbb{R}_{+}^2}\chi \mathcal{Z}_j^2>0$
and $\int_{\mathbb{R}_{+}^2}\chi \mathcal{Z}_j\mathcal{Z}_t=0$ for $j\neq t$.
Hence (\ref{3.26}) implies that $\widehat{\phi}_i^{\infty}\equiv0$
in $B_{R_1}^{+}(0)$.
Furthermore, by (\ref{3.19}) we obtain
$\lim_{n\rightarrow+\infty}\|\phi_n\|_{**}=0$.
But (\ref{3.20}) and (\ref{3.21}) tell us
$\liminf_{n\rightarrow+\infty}\|\phi_n\|_{**}>0$,
which is a contradiction.
\end{proof}

\vspace{1mm}
\vspace{1mm}

{\bf Step 3:} Establishing uniform  an a priori estimate for solutions to
(\ref{3.16}) that satisfy orthogonality conditions with respect to
$Z_{ij}$, $j=1,J_i$ only.

\vspace{1mm}
\vspace{1mm}
\vspace{1mm}

\noindent{\bf Lemma 3.4.}\,\,{\it For $p$  large enough, if
$\phi$ solves {\upshape (\ref{3.16})} and satisfies
\begin{equation}\label{3.27}
\aligned
\int_{\Omega_p}\chi_iZ_{ij}\phi=0\,\,\,\,\,
\,\,\,\,\forall\,\,i=1,\ldots,m,\,\,j=1,J_i,
\endaligned
\end{equation}
then
\begin{equation}\label{3.28}
\aligned
\|\phi\|_{L^{\infty}(\Omega_p)}\leq Cp\, \|h\|_{*},
\endaligned
\end{equation}
where $C>0$ is independent of $p$.}

\vspace{1mm}

\begin{proof}
With no loss of generality we  prove the validity of estimate (\ref{3.28})
only  under the case $q\in\po$, because for the other case $q\in\Omega$
this estimate can also be established   in an analogous but a little bit more simple consideration.

Fix $q\in\po$ and let $R>R_0+1$ be large and fixed, $d>0$ small but fixed.
We denote that for $i=1,\ldots,m$,
\begin{equation}\label{3.29}
\aligned
\widehat{Z}_{q}(y)=Z_{q}(y)-\frac{\varepsilon}{\rho_0v_0}
+a_{q}G(\varepsilon y,q),
\,\qquad\quad\,\,\,\qquad\,
\widehat{Z}_{i0}(y)=Z_{i0}(y)-\frac1{\mu_i}
+a_{i0}G(\varepsilon y,\xi_i),
\endaligned
\end{equation}
where
\begin{equation}\label{3.30}
\aligned
a_{q}=\frac{\varepsilon}{\rho_0v_0\big[H(q,q)-\frac{4(1+\alpha)}{c_0}\log(\rho_0v_0R)\big]},
\,\qquad\quad\qquad\,
a_{i0}=\frac1{\mu_i\big[H(\xi_i,\xi_i)-\frac{4}{c_i}\log(\varepsilon \mu_iR)\big]}.
\endaligned
\end{equation}
From estimate (\ref{2.31}),
definitions (\ref{2.3}), (\ref{2.5}), (\ref{3.6}), (\ref{3.7}), (\ref{3.11}), (\ref{3.12}),
and expansions $(A3)$, $(A5)$
we obtain
\begin{equation}\label{3.31}
\aligned
C_1p\leq-\log(\rho_0v_0R)
\leq C_2p,
\qquad\qquad\qquad
C_1p\leq-\log(\varepsilon \mu_iR)
\leq C_2p,
\endaligned
\end{equation}
and
\begin{equation}\label{3.32}
\aligned
\widehat{Z}_{q}(y)=O\left(
\frac{\varepsilon G(\varepsilon y,q)}{p\rho_0v_0}
\right),
\,\,\,\qquad\qquad\qquad\,\,\,\,
\widehat{Z}_{i0}(y)=O\left(
\frac{\,G(\varepsilon y,\xi_i)\,}{p\mu_i}
\right).
\endaligned
\end{equation}
Let   $\eta_1$ and $\eta_2$ be  radial smooth cut-off functions in $\mathbb{R}^2$ such that
$$
\aligned
&0\leq\eta_1\leq1;\,\,\,\ \,\,\,|\nabla\eta_1|\leq C\,\,\ \textrm{in}\,\,\,\mathbb{R}^2;
\,\,\,\ \,\,\,\eta_1\equiv1\,\,\ \textrm{in}\,\,\,B_R(0);\,\,\,\,\ \,\,\,
\,\eta_1\equiv0\,\,\ \textrm{in}\,\,\,\mathbb{R}^2\setminus B_{R+1}(0);\\[1mm]
&0\leq\eta_2\leq1;\,\,\,\ \,\,\,|\nabla\eta_2|\leq C\,\,\ \textrm{in}\,\,\,\mathbb{R}^2;
\,\,\,\ \,\,\,\eta_2\equiv1\,\,\ \textrm{in}\,\,\,B_{3d}(0);\,\,\,
\,\,\,\,\,\,\eta_2\equiv0\,\,\ \textrm{in}\,\,\,\mathbb{R}^2\setminus B_{6d}(0).
\endaligned
$$
Set
\begin{equation}\label{3.33}
\aligned
\eta_{q1}(y)=
\eta_1\left(\frac{\varepsilon}{\rho_0 v_0}\big|F_q^p(y)\big|\right),
\,\,\qquad\qquad\,\,
\eta_{q2}(y)=
\eta_2\left(\varepsilon\big|F_q^p(y)\big|\right),
\endaligned
\end{equation}
and for any $i=1,\ldots,l$,
\begin{equation}\label{3.34}
\aligned
\eta_{i1}(y)=
\eta_1\left(\frac{1}{\mu_i}\big|y-\xi_i'\big|\right),
\,\,\qquad\quad\,\qquad\,\,
\eta_{i2}(y)=
\eta_2\left(\varepsilon\big|y-\xi'_i\big|\right),
\endaligned
\end{equation}
and for any  $i=l+1,\ldots,m$,
\begin{equation}\label{3.35}
\aligned
\eta_{i1}(y)=
\eta_1\left(\frac{1}{\mu_i}\big|F_i^p(y)\big|\right),
\,\,\qquad\quad\,\qquad\,\,
\eta_{i2}(y)=
\eta_2\left(\varepsilon\big|F_i^p(y)\big|\right).
\endaligned
\end{equation}

Let us define the two test functions
\begin{equation}\label{3.36}
\aligned
\widetilde{Z}_{q}=\eta_{q1}Z_{q}+(1-\eta_{q1})\eta_{q2}\widehat{Z}_{q},
\,\,\,\qquad\,\,\qquad\,\,\,
\widetilde{Z}_{i0}=\eta_{i1}Z_{i0}+(1-\eta_{i1})\eta_{i2}\widehat{Z}_{i0}.
\endaligned
\end{equation}
Given $\phi$ satisfying (\ref{3.16}) and (\ref{3.27}), let
\begin{equation}\label{3.37}
\aligned
\widetilde{\phi}=\phi+d_q\widetilde{Z}_{q}+\sum\limits_{i=1}^{m}d_i\widetilde{Z}_{i0}+\sum_{i=1}^m\sum\limits_{j=1}^{J_i}e_{ij}\chi_iZ_{ij}.
\endaligned
\end{equation}
We will first prove the existence of $d_q$,
$d_i$ and $e_{ij}$  such
that $\widetilde{\phi}$ satisfies
the orthogonality conditions in (\ref{3.17}).
Testing  (\ref{3.37}) against $\chi_iZ_{ij}$ and using the
orthogonality conditions in (\ref{3.17}) and (\ref{3.27}) for $j=1,J_i$
together with the fact that
$\chi_i\chi_k\equiv0$ if $i\neq k$, we find
\begin{equation}\label{3.38}
\aligned
e_{ij}=
\left(-d_q\int_{\Omega_p}\chi_iZ_{ij}\widetilde{Z}_{q}
-\sum_{k=1}^md_k\int_{\Omega_p}\chi_iZ_{ij}\widetilde{Z}_{k0}\right)
\left/\int_{\Omega_p}\chi^2_iZ^2_{ij},
\,\quad\,
i=1,\ldots,m,\,\,
j=1,J_i.
\right.
\endaligned
\end{equation}
Note that  if $i=1,\ldots,l$ and $j=1,2$,
by (\ref{3.7}),
$$
\aligned
\int_{\Omega_p}\chi_iZ_{ij}\widetilde{Z}_{i0}
=
\int_{\mathbb{R}^2}\chi\mathcal{Z}_{j}\mathcal{Z}_{0}
=0,
\,\,\,\qquad\qquad\qquad\,\,\,
\int_{\Omega_p}\chi^2_iZ^2_{ij}
=
\int_{\mathbb{R}^2}\chi^2\mathcal{Z}_{j}^2
=
C_j>0,
\endaligned
$$
while if  $i=l+1,\ldots,m$ and $j=1$,
by (\ref{3.9}), (\ref{3.10})  and (\ref{3.12}),
$$
\aligned
\int_{\Omega_p}\chi_iZ_{ij}\widetilde{Z}_{i0}
=
\int_{\mathbb{R}^2_{+}}\chi\mathcal{Z}_{j}\mathcal{Z}_{0}
[1+O\big(\varepsilon\mu_i|z|\big)]
=O\left(\varepsilon\mu_i\right),
\,\qquad\,
\int_{\Omega_p}\chi^2_iZ^2_{ij}
=
\int_{\mathbb{R}^2_{+}}\chi^2\mathcal{Z}_{j}^2[1+O\big(\varepsilon\mu_i|z|\big)]
=
\frac{C_j}2+O\left(\varepsilon\mu_i\right).
\endaligned
$$
By (\ref{3.32}) and (\ref{3.36}),
$$
\aligned
\int_{\Omega_p}\chi_iZ_{ij}\widetilde{Z}_{q}=
O\left(\frac{\varepsilon\mu_i\log p}{p\rho_0v_0}
\right),
\,\,\,\qquad\qquad\qquad\,\,\,
\int_{\Omega_p}\chi_iZ_{ij}\widetilde{Z}_{k0}dy=
O\left(\frac{\mu_i\log p}{p\mu_k}
\right),\,\,\,\,\,\,\,
\forall\,\,
k\neq i.
\endaligned
$$
Then
\begin{equation}\label{3.39}
\aligned
|e_{ij}|\leq
C\left(\varepsilon\mu_i|d_i|
+
\frac{\varepsilon\mu_i\log p}{p\rho_0v_0}|d_q|
+
\sum_{k\neq i}^m\frac{\mu_i\log p}{p\mu_k}|d_k|\right).
\endaligned
\end{equation}
So we only need  to consider $d_q$  and $d_i$.
Multiplying  definition (\ref{3.37}) by $\chi_qZ_q$ and $\chi_iZ_{i0}$,
respectively, and using
the orthogonality conditions in (\ref{3.17})
for $q$ and $j=0$
and the fact that
$\chi_q\chi_k\equiv0$ for all $k$,
 we obtain a system of $(d_q,d_1,\ldots,d_m)$
\begin{equation}\label{3.40}
\aligned
d_q\int_{\Omega_p}\chi_qZ_q\widetilde{Z}_{q}
+\sum_{k=1}^md_k\int_{\Omega_p}\chi_qZ_q\widetilde{Z}_{k0}
=-\int_{\Omega_p}\chi_qZ_q\phi,
\endaligned
\end{equation}
\begin{equation}\label{3.41}
\aligned
d_q\int_{\Omega_p}\chi_iZ_{i0}\widetilde{Z}_{q}
+\sum_{k=1}^md_k\int_{\Omega_p}\chi_iZ_{i0}\widetilde{Z}_{k0}
+\sum\limits_{t=1}^{J_i}e_{it}\int_{\Omega_p}\chi_i^2Z_{i0}Z_{it}
=-\int_{\Omega_p}\chi_iZ_{i0}\phi,
\,\quad\,i=1,\ldots,m.
\endaligned
\end{equation}
By  (\ref{3.11}),  (\ref{3.12}), (\ref{3.32}) and (\ref{3.36}) we can compute
$$
\aligned
\int_{\Omega_p}\chi_qZ_q\widetilde{Z}_{q}
=
\int_{\mathbb{R}^2_{+}}\chi\mathcal{Z}_{q}^2[1+O\big(\rho_0v_0|z|\big)]
=
\frac{C_q}2+O\left(\rho_0v_0\right),
\,\,\,\qquad\qquad\,\,\,
\int_{\Omega_p}\chi_qZ_q\widetilde{Z}_{k0}=
O\left(\frac{\rho_0v_0\log p}{p\varepsilon\mu_k}
\right),
\endaligned
$$
and
$$
\aligned
\int_{\Omega_p}\chi_iZ_{i0}\widetilde{Z}_{q}=
O\left(
\frac{\varepsilon\mu_i\log p}{p\rho_0v_0}
\right),
\qquad
\qquad\qquad
\qquad
\int_{\Omega_p}\chi_iZ_{i0}\widetilde{Z}_{k0}=
O\left(\frac{\mu_i\log p}{p\mu_k}
\right),
\quad
i\neq k.
\endaligned
$$
Moreover,  if $i=1,\ldots,l$ and $t=1,2$,
by (\ref{3.7}),
$$
\aligned
\int_{\Omega_p}\chi_iZ_{i0}\widetilde{Z}_{i0}
=
\int_{\mathbb{R}^2}\chi\mathcal{Z}_{0}^2
=C_0>0,
\qquad\qquad\quad\qquad
\int_{\Omega_p}\chi_i^2Z_{i0}Z_{it}=
\int_{\mathbb{R}^2}\chi^2\mathcal{Z}_{0}\mathcal{Z}_{t}
=0,
\endaligned
$$
while if  $i=l+1,\ldots,m$ and $t=1$,
by (\ref{3.9}), (\ref{3.10})  and (\ref{3.12}),
$$
\aligned
\int_{\Omega_p}\chi_iZ_{i0}\widetilde{Z}_{i0}
=
\int_{\mathbb{R}^2_{+}}\chi\mathcal{Z}_{0}^2
[1+O\big(\varepsilon\mu_i|z|\big)]
=
\frac{C_0}2+O\left(\varepsilon\mu_i\right),
\quad
\int_{\Omega_p}\chi_i^2Z_{i0}Z_{it}=
\int_{\mathbb{R}^2_{+}}\chi^2\mathcal{Z}_{0}\mathcal{Z}_{t}
[1+O\big(\varepsilon\mu_i|z|\big)]
=O\left(\varepsilon\mu_i\right).
\endaligned
$$
Let us denote $\mathcal{A}$ the coefficient matrix of systems (\ref{3.40})-(\ref{3.41})
with respect to $(d_q,d_1,\ldots,d_m)$.
From the
above estimates it follows that
$P^{-1}\mathcal{A}P$ is diagonally dominant and then invertible, where
$P=\diag\left(\rho_0 v_0\big/\varepsilon,\,\mu_1,\ldots,\mu_m\right)$. Hence $\mathcal{A}$ is also invertible and
$(d_q,d_1,\ldots,d_m)$ is well defined.

Estimate (\ref{3.28}) is a  direct consequence of the following two claims.

\vspace{1mm}
\vspace{1mm}

\noindent{\bf Claim 1.}\,\,{\it
Let $\mathcal{L}=-\Delta_{a(\varepsilon y)}+\varepsilon^2-W_{\xi'}$, then for any $i=1,\ldots,m$
and $j=1, J_i$,
\begin{equation}\label{3.42}
\aligned
\big\|\mathcal{L}(\widetilde{Z}_{q})\big\|_{*}\leq
\frac{C\varepsilon\log p}{p\rho_0v_0},
\endaligned
\end{equation}
and
\begin{equation}\label{3.43}
\aligned
\big\|\mathcal{L}(\chi_iZ_{ij})\big\|_{*}\leq\frac{C}{\mu_i},
\,\quad\qquad\quad\qquad\,\,\,\quad\,
\big\|\mathcal{L}(\widetilde{Z}_{i0})\big\|_{*}\leq
\frac{C\log p}{p\mu_i}.
\endaligned
\end{equation}
}

\noindent{\bf Claim 2.}\,\,{\it
For any $i=1,\ldots,m$ and $j=1,J_i$,
\begin{equation}\label{3.44}
\aligned
|d_q|\leq Cp\frac{\rho_0v_0}{\varepsilon}\|h\|_{*},
\,\quad\quad\,\,\quad\,\,\quad
|d_i|\leq Cp\mu_i\|h\|_{*},
\,\quad\qquad\quad\,\,\quad
|e_{ij}|\leq C\mu_i\log p\,\|h\|_{*}.
\endaligned
\end{equation}
}

In fact, by the definition of
$\widetilde{\phi}$ in (\ref{3.37}) we obtain
\begin{equation}\label{3.45}
\aligned\left\{\aligned
&
\mathcal{L}(\widetilde{\phi})=h
+d_q\mathcal{L}(\widetilde{Z}_{q})
+\sum\limits_{i=1}^{m}d_i\mathcal{L}(\widetilde{Z}_{i0})
+\sum_{i=1}^m\sum_{j=1}^{J_i}e_{ij}\mathcal{L}(\chi_iZ_{ij})
\,\quad\,\textrm{in}\,\,\,\ \,\,\Omega_p,\\
&\frac{\partial\widetilde{\phi}}{\partial\nu}=0
\,\qquad\qquad\qquad\qquad\quad\qquad\,
\,\qquad\qquad\qquad\qquad\qquad\quad
\,\textrm{on}\,\,\,\,\,\partial\Omega_p.
\endaligned\right.
\endaligned
\end{equation}
Since  the orthogonality conditions in  (\ref{3.17}) hold, by estimate (\ref{3.18}), Claims $1$ and $2$  we conclude
\begin{equation}\label{3.46}
\aligned
\|\widetilde{\phi}\|_{L^{\infty}(\Omega_p)}\leq
C\left\{\|h\|_{*}
+|d_q|\big\|\mathcal{L}(\widetilde{Z}_{q})\big\|_{*}
+\sum\limits_{i=1}^{m}|d_i|\big\|\mathcal{L}(\widetilde{Z}_{i0})\big\|_{*}
+\sum_{i=1}^m\sum_{j=1}^{J_i}|e_{ij}|\big\|\mathcal{L}(\chi_i Z_{ij})\big\|_{*}
\right\}\leq C \log p\,\|h\|_{*}.
\endaligned
\end{equation}
Using the definition of $\widetilde{\phi}$
again and the fact that
\begin{equation}\label{3.47}
\aligned
\big\|\widetilde{Z}_{q}\big\|_{L^{\infty}(\Omega_p)}\leq\frac{C\varepsilon}{\rho_0v_0},
\,\,\,\qquad
\big\|\widetilde{Z}_{i0}\big\|_{L^{\infty}(\Omega_p)}\leq\frac{C}{\mu_i},
\,\,\,\qquad
\big\|\chi_iZ_{ij}\big\|_{L^{\infty}(\Omega_p)}\leq\frac{C}{\mu_i},
\quad\,i=1,\ldots,m,\,\,j=1,J_i,
\endaligned
\end{equation}
estimate (\ref{3.28}) then follows from estimate (\ref{3.46}) and Claim $2$.

\vspace{1mm}
\vspace{1mm}
\vspace{1mm}

\noindent{\bf Proof of Claim 1.}
Let us try with inequality (\ref{3.42}).
Fixing $q\in\po$  and observing   that
 $F_q^p(q')=(0,0)$ and $\nabla F_q^p(q')=A_q$,
 by (\ref{3.8}) and (\ref{3.49}) we find
 \begin{equation}\label{3.48}
\aligned
z_q:=F_q^p(y)=\frac{1}{\varepsilon}F_q\big(A_q(\varepsilon y-q)\big)
=A_q(y-q')\big[\,1+O\left(\varepsilon A_q(y-q')\right)\big],
\endaligned
\end{equation}
and
\begin{equation}\label{3.50}
\aligned
\nabla_y=A_q\nabla_{z_q}+O(\varepsilon|z_q|)
\nabla_{z_q}
\,\quad\ \,\,\,\quad\,
\textrm{and}
\,\quad\ \,\,\,\quad\,
-\Delta_y=-\Delta_{z_q}+O(\varepsilon|z_q|)
\nabla_{z_q}^2+O(\varepsilon)\nabla_{z_q}.
\endaligned
\end{equation}
Furthermore,
\begin{equation}\label{3.51}
\aligned
-\Delta_{a(\varepsilon y)}
=-\frac1{a(\varepsilon y)}\nabla_y\big(a(\varepsilon y)\nabla_y(\cdot)\big)
=-\Delta_y-\varepsilon\nabla_x\log a(\varepsilon y)\nabla_y
=-\Delta_{z_q}+O(\varepsilon|z_q|)
\nabla_{z_q}^2+O(\varepsilon)\nabla_{z_q}.
\endaligned
\end{equation}
Thus by (\ref{3.2}) and  (\ref{3.11}),
\begin{equation}\label{3.14}
\aligned
Z_{q}-\frac{\varepsilon}{\rho_0 v_0}
=-
\frac{2\varepsilon}{\rho_0 v_0}
\left[1+\left|\frac{\varepsilon z_q}{\rho_0 v_0}\right|^{2(1+\alpha)}\right]^{-1}
=
O\left(\frac{\varepsilon}{\rho_0 v_0}
\left[1+\frac{|\varepsilon y-q|}{\rho_0 v_0}\right]^{-2(1+\alpha)}\right),
\endaligned
\end{equation}
and
\begin{equation}\label{3.13}
\aligned
\Delta_{a(\varepsilon y)}Z_{q}+\left(\frac{\varepsilon}{\rho_0v_0}\right)^2
\frac{8(1+\alpha)^2\big|\frac{\varepsilon y-q}{\rho_0 v_0}\big|^{2\alpha}}
{\big(1+\big|\frac{\varepsilon y-q}{\rho_0 v_0}\big|^{2(1+\alpha)}\big)^2}Z_{q}
=O\left(\varepsilon\left(\frac{\varepsilon}{\rho_0v_0}\right)^2
\left|\frac{\varepsilon y-q}{\rho_0 v_0}\right|^{2\alpha}
\left[1+\frac{|\varepsilon y-q|}{\rho_0 v_0}\right]^{-3-4\alpha}\right).
\endaligned
\end{equation}
Consider the four regions
$$
\aligned
\Omega_1=(F_q^p)^{-1}\left(\left\{\,\left|\frac{\varepsilon z_q}{\rho_0 v_0}\right|\leq R\right\}\cap\mathbb{R}_+^2\right),
\,\quad\,\quad\quad\qquad\qquad\,
\Omega_2=(F_q^p)^{-1}\left(\left\{R<\left|\frac{\varepsilon z_q}{\rho_0 v_0}\right|\leq R+1\right\}\cap\mathbb{R}_+^2\right),
\,\,\quad\qquad\\[1mm]
\Omega_3=(F_q^p)^{-1}\left(\left\{R+1<\left|\frac{\varepsilon z_q}{\rho_0 v_0}\right|\leq\frac{3d}{\rho_0v_0}\right\}\cap\mathbb{R}_+^2\right),
\quad\quad\qquad\quad
\Omega_4=(F_q^p)^{-1}\left(\left\{\frac{3d}{\rho_0v_0}<\left|\frac{\varepsilon z_q}{\rho_0 v_0}\right|\leq
\frac{6d}{\rho_0v_0}\right\}\cap\mathbb{R}_+^2\right).
\quad\quad\,\,
\endaligned
$$
Notice that for any $y\in\Omega_1\cup\Omega_2$, by (\ref{2.50}) and  (\ref{3.11}),
\begin{equation}\label{3.52}
\aligned
\left[\left(\frac{\varepsilon}{\rho_0v_0}\right)^2
\frac{8(1+\alpha)^2\big|\frac{\varepsilon y-q}{\rho_0 v_0}\big|^{2\alpha}}
{\big(1+\big|\frac{\varepsilon y-q}{\rho_0 v_0}\big|^{2(1+\alpha)}\big)^2}
-W_{\xi'}\right]Z_{q}
=
O\left(\frac{1}{p}\left(\frac{\varepsilon}{\rho_0v_0}\right)^3
\left|\frac{\varepsilon y-q}{\rho_0 v_0}\right|^{2\alpha}
\left[1+\frac{|\varepsilon y-q|}{\rho_0 v_0}\right]^{-4(1+\alpha)}\right),
\endaligned
\end{equation}
and for any
$y\in\Omega_2\cup\Omega_3\cup\Omega_4$ and any $\beta\in(0,1)$,
by (\ref{1.3}), (\ref{2.20}), (\ref{3.29}) and (\ref{3.30}),
\begin{equation}\label{3.53}
\aligned
Z_{q}-\widehat{Z}_{q}=\frac{\varepsilon}{\rho_0 v_0}
-a_{q}G(\varepsilon y,q)=
a_{q}\left[
\frac{4(1+\alpha)}{c_0}\log\frac{|\varepsilon y-q|}{R\rho_0 v_0}+O\left(
|\varepsilon y-q|^\beta\right)
\right].
\endaligned
\end{equation}
In $\Omega_1$,  by  (\ref{3.13}) and (\ref{3.52}),
\begin{equation}\label{3.54}
\aligned
\mathcal{L}(\widetilde{Z}_{q})=\mathcal{L}(Z_{q})=
O\left(\frac{1}{p}\left(\frac{\varepsilon}{\rho_0v_0}\right)^3
\left|\frac{\varepsilon y-q}{\rho_0 v_0}\right|^{2\alpha}\right).
\endaligned
\end{equation}
In $\Omega_2$, by (\ref{1.2}),  (\ref{3.29}) and (\ref{3.36}),
\begin{eqnarray*}
\mathcal{L}(\widetilde{Z}_{q})
=
\mathcal{L}(Z_{q})-(1-\eta_{q1})\mathcal{L}(Z_{q}-\widehat{Z}_{q})
-2\nabla\eta_{q1}\nabla(Z_{q}-\widehat{Z}_{q})
-(Z_{q}-\widehat{Z}_{q})\Delta_{a(\varepsilon y)}\eta_{q1}
\qquad\qquad\qquad\qquad\qquad\quad\,
&&\nonumber\\[1.5mm]
=\left[-\Delta_{a(\varepsilon y)}Z_{q}-\left(\frac{\varepsilon}{\rho_0v_0}\right)^2
\frac{8(1+\alpha)^2\big|\frac{\varepsilon y-q}{\rho_0 v_0}\big|^{2\alpha}}
{\big(1+\big|\frac{\varepsilon y-q}{\rho_0 v_0}\big|^{2(1+\alpha)}\big)^2}Z_{q}
\right]+\left[\left(\frac{\varepsilon}{\rho_0v_0}\right)^2
\frac{8(1+\alpha)^2\big|\frac{\varepsilon y-q}{\rho_0 v_0}\big|^{2\alpha}}
{\big(1+\big|\frac{\varepsilon y-q}{\rho_0 v_0}\big|^{2(1+\alpha)}\big)^2}
-W_{\xi'}\right]Z_{q}
\quad\,\,\,
&&\nonumber\\
+\varepsilon^2\left(
Z_{q}-\frac{\varepsilon}{\rho_0 v_0}
\right)
+\frac{\varepsilon^3}{\rho_0 v_0}\eta_{q1}
+W_{\xi'}(1-\eta_{q1})(Z_{q}-\widehat{Z}_{q})
-2\nabla\eta_{q1}\nabla(Z_{q}-\widehat{Z}_{q})
-(Z_{q}-\widehat{Z}_{q})\Delta_{a(\varepsilon y)}\eta_{q1}.
&&
\end{eqnarray*}
Notice that in $\Omega_2$, by  (\ref{3.30}), (\ref{3.31}), (\ref{3.48}) and (\ref{3.53}),
\begin{equation}\label{3.55}
\aligned
|Z_{q}-\widehat{Z}_{q}|=O\left(
\frac{\varepsilon}{p\rho_0v_0R}
\right)
\,\,\quad\quad\quad\,\textrm{and}\,\,\quad\quad\quad\,
|\nabla\big(Z_{q}-\widehat{Z}_{q}\big)|=O\left(
\frac{\varepsilon^2}{p\rho_0^2v_0^2R}
\right).
\endaligned
\end{equation}
Moreover, $|\nabla\eta_{q1}|=O\big(\varepsilon/(\rho_0v_0)\big)$ and
$|\Delta_{a(\varepsilon y)}\eta_{q1}|=O\big(\varepsilon^2/(\rho_0^2v_0^2)\big)$.
From (\ref{2.50}), (\ref{3.14})-(\ref{3.52}) and (\ref{3.55}) we obtain
\begin{equation}\label{3.56}
\aligned
\mathcal{L}(\widetilde{Z}_{q})
=O\left(\frac{\varepsilon^3}{p\rho_0^3v_0^3R}
\right).
\endaligned
\end{equation}
In $\Omega_3$, by  (\ref{1.2}),  (\ref{3.29}),  (\ref{3.36}), (\ref{3.14}) and (\ref{3.13}),
$$
\aligned
\mathcal{L}(\widetilde{Z}_{q})=&\mathcal{L}(\widehat{Z}_{q})
=\mathcal{L}(Z_{q})-\mathcal{L}(Z_{q}-\widehat{Z}_{q})
\\[1.6mm]
\equiv&A_1+A_2+O\left(\varepsilon\left(\frac{\varepsilon}{\rho_0v_0}\right)^2
\left[1+\frac{|\varepsilon y-q|}{\rho_0 v_0}\right]^{-3-2\alpha}\right)
+O\left(\frac{\varepsilon^3}{\rho_0 v_0}
\left[1+\frac{|\varepsilon y-q|}{\rho_0 v_0}\right]^{-2(1+\alpha)}\right),
\endaligned
$$
where
$$
\aligned
A_1=\left[\left(\frac{\varepsilon}{\rho_0v_0}\right)^2
\frac{8(1+\alpha)^2\big|\frac{\varepsilon y-q}{\rho_0 v_0}\big|^{2\alpha}}
{\big(1+\big|\frac{\varepsilon y-q}{\rho_0 v_0}\big|^{2(1+\alpha)}\big)^2}
-W_{\xi'}\right]Z_{q}
\qquad\,\,
\textrm{and}
\,\,\qquad\,\,
A_2=W_{\xi'}\left[\frac{\varepsilon}{\rho_0 v_0}
-a_{q}G(\varepsilon y,q)\right].
\endaligned
$$
For the estimates of   these two terms, we
set  $z_k:=y-\xi_k'$ for all
$k=1,\ldots,l$, but $z_k:=F_{k}^p(y)$ for all
$k=l+1,\ldots,m$, and divide $\Omega_3$ into several pieces:
$$
\aligned
\Omega_{3,k}=\Omega_3\cap\left\{\,
\left|\frac{z_k}{\mu_k}\right|
=\left|\frac{y-\xi_k'}{\mu_k}\right|
\leq
\frac{1}{3\,p^{2\kappa}\sqrt{\varepsilon\mu_k}}\right\}
\,\,\qquad\,\,\,\forall\,\,
k=1,\ldots,l,
\endaligned
$$
$$
\aligned
\Omega_{3,k}=\Omega_3\cap(F_k^p)^{-1}\left(\left\{\,
\left|\frac{z_k}{\mu_k}\right|\leq
\frac{1}{3\,p^{2\kappa}\sqrt{\varepsilon\mu_k}}\right\}\cap\mathbb{R}_+^2\right)
\,\,\,\,\,\forall\,\,
k=l+1,\ldots,m,
\endaligned
$$
and
$$
\aligned
\Omega_{q}=(F_q^p)^{-1}\left(\left\{R+1<\left|\frac{\varepsilon z_q}{\rho_0 v_0}\right|
\leq\frac{1}{3\,p^{2\kappa}\sqrt{\rho_0 v_0}}\right\}\cap\mathbb{R}_+^2\right)
\qquad\,
\textrm{and}
\qquad\quad
\widetilde{\Omega}_3=\Omega_3\setminus
\left[\bigcup_{k=1}^m\Omega_{3,k}\cup \Omega_{q}\right].
\endaligned
$$
Note that for any $k=l+1,\ldots,m$, by (\ref{3.9}) and (\ref{3.10}),
\begin{equation}\label{3.57}
\aligned
z_k=F_{k}^p(y)=\frac{1}{\varepsilon}F_k\big(A_k(\varepsilon y-\xi_k)\big)
=A_k(y-\xi'_k)
\big[\,1+O\left(\varepsilon A_k(y-\xi'_k)\right)\big].
\endaligned
\end{equation}
From  (\ref{2.1}), (\ref{2.12}), (\ref{2.49}),  (\ref{2.50}),  (\ref{3.48}) and (\ref{3.57}) we have
$$
\aligned
A_1=\left\{
\aligned
&\left(\frac{\varepsilon}{\rho_0v_0}\right)^3
\left|\frac{\varepsilon y-q}{\rho_0v_0}\right|^{-4-2\alpha}
O\left(\frac{1}{p}\log^2\left|\frac{\varepsilon y-q}{\rho_0v_0}\right|
+
\frac{1}{p^2}\log^4\left|\frac{\varepsilon y-q}{\rho_0v_0}\right|\right)
\,\,\,\quad
\textrm{in}\,\,\,\,\Omega_{q},\\[1mm]
&O\left(\varepsilon^3(\rho_0v_0)^{\alpha-1}p^{4\kappa(2+\alpha)}\right)
\quad\qquad\qquad\qquad\qquad
\qquad\qquad\qquad\qquad\quad\,\,
\textrm{in}\,\,\,\,\widetilde{\Omega}_{3}.
\endaligned
\right.
\endaligned
$$
Moreover, by  (\ref{3.30}), (\ref{3.31}) and (\ref{3.53}),
$$
\aligned
A_2=
\left(\frac{\varepsilon}{\rho_0v_0}\right)^3
\left|\frac{\varepsilon y-q}{\rho_0v_0}\right|^{-4-2\alpha}
O\left(
\frac{1}{p}\log\left|\frac{\varepsilon y-q}{R\rho_0v_0}\right|
\right)
\,\quad\,\,\textrm{in}\,\,\,\,\Omega_{q}\cup \widetilde{\Omega}_{3}.
\endaligned
$$
Thus in $\Omega_{q}\cup \widetilde{\Omega}_{3}$,
\begin{equation}\label{3.58}
\aligned
\mathcal{L}(\widetilde{Z}_{q})=\mathcal{L}(\widehat{Z}_{q})
=
\left(\frac{\varepsilon}{\rho_0v_0}\right)^3
\left|\frac{\varepsilon y-q}{\rho_0v_0}\right|^{-4-2\alpha}
O
\left(
\frac{1}{p}\log\left|\frac{\varepsilon y-q}{R\rho_0v_0}\right|
+\frac{1}{p}\log^2\left|\frac{\varepsilon y-q}{\rho_0v_0}\right|
+
\frac{1}{p^2}\log^4\left|\frac{\varepsilon y-q}{\rho_0v_0}\right|
\right).
\endaligned
\end{equation}
In $\Omega_{3,k}$ with $k=1,\ldots,m$, by (\ref{2.51}), (\ref{3.32}),
(\ref{3.14}), (\ref{3.13}) and (\ref{3.57}),
\begin{eqnarray}\label{3.59}
\mathcal{L}(\widetilde{Z}_{q})=\mathcal{L}(\widehat{Z}_{q})
=
\left(\frac{\varepsilon}{\rho_0v_0}\right)^2
\frac{8(1+\alpha)^2\big|\frac{\varepsilon y-q}{\rho_0 v_0}\big|^{2\alpha}}
{\big(1+\big|\frac{\varepsilon y-q}{\rho_0 v_0}\big|^{2(1+\alpha)}\big)^2}Z_{q}
-\left[
\Delta_{a(\varepsilon y)} Z_{q}+\left(\frac{\varepsilon}{\rho_0v_0}\right)^2
\frac{8(1+\alpha)^2\big|\frac{\varepsilon y-q}{\rho_0 v_0}\big|^{2\alpha}}
{\big(1+\big|\frac{\varepsilon y-q}{\rho_0 v_0}\big|^{2(1+\alpha)}\big)^2}Z_{q}
\right]
&&\nonumber\\[0.1mm]
+\,
\varepsilon^2\left(
Z_{q}-\frac{\varepsilon}{\rho_0 v_0}
\right)-W_{\xi'}\widehat{Z}_{q}
\qquad\qquad\qquad\qquad\qquad\qquad\qquad\qquad
\qquad\qquad\qquad\quad\,\,\,\,
&&\nonumber\\[0.1mm]
=\frac{1}{\mu_k^2}\frac{8}{\big(1+\big|\frac{y-\xi'_k}{\mu_k}\big|^2\big)^2}
\,O\left(
\frac{\varepsilon\log p}{p\rho_0v_0}
\right).
\qquad\qquad\qquad\qquad\qquad\qquad\qquad\quad
\qquad\quad\qquad\quad\quad
&&
\end{eqnarray}
In $\Omega_4$,  by  (\ref{1.2}),  (\ref{3.29}) and (\ref{3.36}),
$$
\aligned
\mathcal{L}(\widetilde{Z}_{q})
=&
\left(\frac{\varepsilon}{\rho_0v_0}\right)^2
\frac{8(1+\alpha)^2\big|\frac{\varepsilon y-q}{\rho_0 v_0}\big|^{2\alpha}}
{\big(1+\big|\frac{\varepsilon y-q}{\rho_0 v_0}\big|^{2(1+\alpha)}\big)^2}
\eta_{q2}Z_{q}
-\eta_{q2}\left[
\Delta_{a(\varepsilon y)} Z_{q}+\left(\frac{\varepsilon}{\rho_0v_0}\right)^2
\frac{8(1+\alpha)^2\big|\frac{\varepsilon y-q}{\rho_0 v_0}\big|^{2\alpha}}
{\big(1+\big|\frac{\varepsilon y-q}{\rho_0 v_0}\big|^{2(1+\alpha)}\big)^2}Z_{q}
\right]
\\[0.1mm]
&+
\varepsilon^2\eta_{q2}\left(
Z_{q}-\frac{\varepsilon}{\rho_0 v_0}
\right)-W_{\xi'}\eta_{q2}\widehat{Z}_{q}
-2\nabla\eta_{q2}\nabla\widehat{Z}_{q}-\widehat{Z}_{q}\Delta_{a(\varepsilon y)}\eta_{q2}.
\endaligned
$$
From the previous choice of the number
$d$ we have  that  for any $y\in\Omega_4$ and any $k=1,\ldots,m$,
by  (\ref{2.5}) and (\ref{3.48}),
$$
\aligned
|y-\xi'_k|\geq|y-q'|-|q'-\xi_k'|\geq\frac{2d}{\varepsilon}
-\frac{d}{\varepsilon}=\frac{d}{\varepsilon}>\frac{\sqrt{\delta_i}}
{\,\varepsilon p^{2\kappa}\,}.
\endaligned
$$
Then by  (\ref{2.31}) and (\ref{2.49}) we find $W_{\xi'}=O(\varepsilon^{4-\beta})$ in $\Omega_4$.
In addition,
$|\nabla\eta_{q2}|=O\left(\varepsilon/d\right)$,
$|\Delta_{a(\varepsilon y)}\eta_{q2}|=O\left(\varepsilon^2/d^2\right)$,
\begin{equation}\label{3.60}
\aligned
|\widehat{Z}_{q}|=O\left(
\frac{\varepsilon|\log d|}{p\rho_0v_0}
\right)\,\,\quad\quad\,\textrm{and}\,\,\quad\quad\,
|\nabla\widehat{Z}_{q}|=O\left(
\frac{\varepsilon^2}{pd\rho_0v_0}
\right)
\,\qquad\,
\textrm{in}
\,\,\,\,
\Omega_4.
\endaligned
\end{equation}
Furthermore, by (\ref{3.14}) and (\ref{3.13}),
\begin{equation}\label{3.61}
\aligned
\mathcal{L}(\widetilde{Z}_{q})
=O\left(
\frac{\varepsilon^3|\log d|}{pd^2\rho_0v_0}
\right).
\endaligned
\end{equation}
Hence by (\ref{3.54}), (\ref{3.56}),  (\ref{3.58}), (\ref{3.59}), (\ref{3.61}) and
the definition of  $\|\cdot\|_*$ with respect to (\ref{2.47}),  we arrive at
$$
\aligned
\big\|\mathcal{L}(\widetilde{Z}_{q})\big\|_{*}=O\left(
\frac{\varepsilon\log p}{p\rho_0v_0}
\right).
\endaligned
$$

The  inequalities in (\ref{3.43}) are easy to establish as they are very
similar to the consideration of inequality (\ref{3.42}), so we  leave the detailed
proof to readers.

\vspace{1mm}
\vspace{1mm}
\vspace{1mm}

\noindent{\bf Proof of Claim 2.}
Testing equation (\ref{3.45}) against
$a(\varepsilon y)\widetilde{Z}_{q}$  and using estimates (\ref{3.46})-(\ref{3.47}),
we can derive that
$$
\aligned
d_q
\int_{\Omega_p}&a(\varepsilon y)\widetilde{Z}_{q}\mathcal{L}(\widetilde{Z}_{q})
+
\sum_{k=1}^md_k
\int_{\Omega_p}a(\varepsilon y)\widetilde{Z}_{k0}\mathcal{L}(\widetilde{Z}_{q})
\\
=&-\int_{\Omega_p} a(\varepsilon y)h\widetilde{Z}_{q}
+\int_{\Omega_p}a(\varepsilon y)\widetilde{\phi}\mathcal{L}(\widetilde{Z}_{q})-\sum_{k=1}^m\sum_{t=1}^{J_k}e_{kt}
\int_{\Omega_p}a(\varepsilon y)\chi_kZ_{kt}\mathcal{L}(\widetilde{Z}_{q})
\\[1mm]
\leq&\frac{C\varepsilon}{\,\rho_0v_0\,}\|h\|_{*}
+C\big\|\mathcal{L}(\widetilde{Z}_{q})\big\|_{*}
\left(\|\widetilde{\phi}\|_{L^{\infty}(\Omega_p)}
+\sum_{k=1}^m\sum_{t=1}^{J_k}\frac{1}{\mu_k}|e_{kt}|\right)
\\[1mm]
\leq&\frac{C\varepsilon}{\,\rho_0v_0\,}\|h\|_{*}
+
C\big\|\mathcal{L}(\widetilde{Z}_{q})\big\|_{*}
\left[\|h\|_{*}
+|d_q|\big\|\mathcal{L}(\widetilde{Z}_{q})\big\|_{*}
+
\sum\limits_{k=1}^{m}|d_k|\big\|\mathcal{L}(\widetilde{Z}_{k0})\big\|_{*}
+\sum_{k=1}^m\sum_{t=1}^{J_k}|e_{kt}|\left(\frac{1}{\mu_k}
+\big\|\mathcal{L}(\chi_kZ_{kt})\big\|_{*}
\right)
\right],
\endaligned
$$
where we have  applied  the following  two inequalities:
\begin{equation*}\label{3.83}
\aligned
\left(\frac{\varepsilon}{\rho_0v_0}\right)^2\int_{\Omega_p}
\frac{\big|\frac{\varepsilon y-q}{\rho_0v_0}\big|^{2\alpha}}
{\,\big(1+\big|\frac{\varepsilon y-q}{\rho_0v_0}\big|\big)^{4+2\hat{\alpha}+2\alpha}\,}
dy\leq C
\qquad\,\,
\textrm{and}
\,\,\qquad\,\,
\int_{\Omega_p}
\frac{1}{\mu_i^2}\frac{1}{\,\big(1+\big|\frac{y-\xi'_i}{\mu_i}\big|\big)^{4+2\hat{\alpha}}\,}
dy\leq C,
\,\,\,\,\,\,\,i=1,\ldots,m.
\endaligned
\end{equation*}
This, combined with inequality (\ref{3.39})  and Claim $1$, gives
\begin{equation}\label{3.62}
\aligned
|d_q|\left|
\int_{\Omega_p}a(\varepsilon y)\widetilde{Z}_{q}\mathcal{L}(\widetilde{Z}_{q})
\right|
\leq
\frac{C\varepsilon}{\rho_0v_0}\|h\|_{*}
+\frac{C\varepsilon\log^2 p}{p\rho_0v_0}\left(
\frac{\varepsilon|d_q|}{p\rho_0v_0}
+
\sum_{k=1}^m\frac{|d_k|}{p\mu_k}
\right)
+\sum_{k=1}^m\left|d_k
\int_{\Omega_p}a(\varepsilon y)\widetilde{Z}_{k0}\mathcal{L}(\widetilde{Z}_{q})
\right|.
\endaligned
\end{equation}
Similarly, testing  (\ref{3.45}) against $a(\varepsilon y)\widetilde{Z}_{i0}$,
$i=1,\ldots,m$
and using
(\ref{3.39}), (\ref{3.46}), (\ref{3.47}) and Claim $1$, we find
\begin{eqnarray}\label{3.63}
|d_i|\left|
\int_{\Omega_p}a(\varepsilon y)\widetilde{Z}_{i0}\mathcal{L}(\widetilde{Z}_{i0})
\right|
\leq
\frac{C}{\mu_i}\|h\|_{*}
+\frac{C\log^2 p}{p\mu_i}\left(
\frac{\varepsilon|d_q|}{p\rho_0v_0}
+
\sum_{k=1}^m\frac{|d_k|}{p\mu_k}
\right)
+
\left|d_q
\int_{\Omega_p}a(\varepsilon y)\widetilde{Z}_{i0}\mathcal{L}(\widetilde{Z}_{q})
\right|
&&\nonumber\\
+
\sum_{k\neq i}^m\left|d_k
\int_{\Omega_p}a(\varepsilon y)\widetilde{Z}_{k0}\mathcal{L}(\widetilde{Z}_{i0})
\right|.
\qquad\qquad\qquad\qquad\qquad\,
\qquad\quad\quad\,\,\,\,
&&
\end{eqnarray}
To achieve  the estimates  of $d_q$, $d_i$ and $e_{ij}$ in (\ref{3.44}),  we have the following claim.

\vspace{1mm}
\vspace{1mm}
\vspace{1mm}

\noindent{\bf Claim 3.}\,\,{\it
If  $d$ is sufficiently small,
but $R$  is  sufficiently large, then we have that
for any $i,k=1,\ldots,m$ with $i\neq k$,
\begin{equation}\label{3.64}
\aligned
\int_{\Omega_p}a(\varepsilon y)\widetilde{Z}_{i0}\mathcal{L}(\widetilde{Z}_{i0})
=\frac{2c_i a(\xi_i)}{p\mu_i^2}\left[1+O\left(\frac1{R^2}\right)\right],
\qquad
\quad
\qquad
\int_{\Omega_p}a(\varepsilon y)\widetilde{Z}_{k0}\mathcal{L}(\widetilde{Z}_{i0})
=O\left(\frac{\log^2p}{p^2\mu_i\mu_k}
\right),
\endaligned
\end{equation}
and
\begin{equation}\label{3.65}
\aligned
\int_{\Omega_p}a(\varepsilon y)\widetilde{Z}_{q}\mathcal{L}(\widetilde{Z}_{q})
=\frac{2c_0 a(q)}{p}\left(\frac{\varepsilon}{\rho_0v_0}\right)^2
\left[1+O\left(\frac1{R^{2(1+\alpha)}}\right)\right],
\quad\,
\quad\,
\int_{\Omega_p}a(\varepsilon y)\widetilde{Z}_{k0}\mathcal{L}(\widetilde{Z}_{q})
=O\left(\frac{\varepsilon\log^2p}{p^2\rho_0v_0\mu_k}
\right),
\endaligned
\end{equation}
where the coefficients  $c_0$ and $c_i$, $i=1,\ldots,m$ are defined in  {\upshape(\ref{2.20})}.
}

\vspace{1mm}
\vspace{1mm}
\vspace{1mm}
\vspace{1mm}

In fact, once Claim $3$ is valid,  then inserting (\ref{3.64}) and (\ref{3.65})
into (\ref{3.63}) and (\ref{3.62}), respectively,  we give
\begin{equation}\label{3.66}
\aligned
\frac{\varepsilon|d_q|}{\,p\rho_0v_0\,}
\leq
C\|h\|_{*}
+\frac{\,C\log^2 p\,}{\,p\,}\left(
\frac{\varepsilon|d_q|}{\,p\rho_0v_0\,}
+
\sum_{k=1}^m\frac{|d_k|}{\,p\mu_k\,}
\right),
\endaligned
\end{equation}
and  for any $i=1,\ldots,m$,
\begin{equation}\label{3.67}
\aligned
\frac{|d_i|}{\,p\mu_i\,}
\leq
C\|h\|_{*}
+\frac{\,C\log^2p\,}{p}
\left(
\frac{\varepsilon|d_q|}{\,p\rho_0v_0\,}
+\sum_{k=1}^m\frac{|d_k|}{\,p\mu_k\,}
\right).
\endaligned
\end{equation}
As a result, using linear algebra arguments for (\ref{3.66})-(\ref{3.67}), we can prove  Claim 2 for
$d_q$ and $d_i$,  and then
complete the proof by inequality (\ref{3.39}).

\vspace{1mm}
\vspace{1mm}
\vspace{1mm}

\noindent{\bf Proof of Claim 3.}
Let us first establish the validity of the two expansions in (\ref{3.65}).
Observe that
\begin{eqnarray}\label{3.68}
&&\mathcal{L}(\widetilde{Z}_{q})=\eta_{q1}\mathcal{L}(Z_{q}-\widehat{Z}_{q})
+\eta_{q2}\mathcal{L}(\widehat{Z}_{q})
-(Z_{q}-\widehat{Z}_{q})\Delta_{a(\varepsilon y)}\eta_{q1}
-2\nabla\eta_{q1}\nabla(Z_{q}-\widehat{Z}_{q})
\nonumber\\[1mm]
&&\qquad\qquad\,-2\nabla\eta_{q2}\nabla\widehat{Z}_{q}-\widehat{Z}_{q}\Delta_{a(\varepsilon y)}\eta_{q2}.
\end{eqnarray}
Thus by (\ref{3.29}) and (\ref{3.36}),
$$
\aligned
\int_{\Omega_p}a(\varepsilon y)\widetilde{Z}_{q}\mathcal{L}(\widetilde{Z}_{q})
=K+L,
\endaligned
$$
where
$$
\aligned
K=\int_{\Omega_p}a(\varepsilon y)\widetilde{Z}_{q}\left[
-(Z_{q}-\widehat{Z}_{q})\Delta_{a(\varepsilon y)}\eta_{q1}
-2\nabla\eta_{q1}\nabla(Z_{q}-\widehat{Z}_{q})
-2\nabla\eta_{q2}\nabla\widehat{Z}_{q}-\widehat{Z}_{q}\Delta_{a(\varepsilon y)}\eta_{q2}
\right],
\endaligned
$$
and
$$
\aligned
L=&\int_{\Omega_p}a(\varepsilon y)\widetilde{Z}_{q}\left[
\eta_{q1}\mathcal{L}(Z_{q}-\widehat{Z}_{q})
+\eta_{q2}\mathcal{L}(\widehat{Z}_{q})
\right]\\[1mm]
=&\int_{\Omega_p}a(\varepsilon y)\eta_{q2}^2\left[Z_{q}-\big(1-\eta_{q1}\big)
\left(\frac{\varepsilon}{\rho_0v_0}-a_{q}G(\varepsilon y,q)\right)\right]\left\{
\big(1-\eta_{q1}\big)W_{\xi'}\left(\frac{\varepsilon}{\rho_0v_0}-a_{q}G(\varepsilon y,q)\right)
+\varepsilon^2\left(
Z_{q}-\frac{\varepsilon}{\rho_0 v_0}
\right)
\right.\\[1mm]
&\left.
-\left[\Delta_{a(\varepsilon y)}Z_{q}+\left(\frac{\varepsilon}{\rho_0v_0}\right)^2
\frac{8(1+\alpha)^2\big|\frac{\varepsilon y-q}{\rho_0 v_0}\big|^{2\alpha}}
{\big(1+\big|\frac{\varepsilon y-q}{\rho_0 v_0}\big|^{2(1+\alpha)}\big)^2}Z_{q}
\right]+\left[\left(\frac{\varepsilon}{\rho_0v_0}\right)^2
\frac{8(1+\alpha)^2\big|\frac{\varepsilon y-q}{\rho_0 v_0}\big|^{2\alpha}}
{\big(1+\big|\frac{\varepsilon y-q}{\rho_0 v_0}\big|^{2(1+\alpha)}\big)^2}
-W_{\xi'}\right]Z_{q}
+\frac{\varepsilon^3}{\rho_0 v_0}\eta_{q1}
\right\}.
\endaligned
$$
Integrating by parts the first  term and the last term of $K$, respectively, we obtain
$$
\aligned
K=&-\int_{\Omega_2}a(\varepsilon y)Z_{q}\nabla\eta_{q1}\nabla(Z_{q}-\widehat{Z}_{q})
+\int_{\Omega_2}a(\varepsilon y)\big(Z_{q}-\widehat{Z}_{q}\big)\nabla\eta_{q1}\nabla(Z_{q}-\widehat{Z}_{q})
\\
&+\int_{\Omega_2}a(\varepsilon y)(Z_{q}-\widehat{Z}_{q})^2|\nabla\eta_{q1}|^2+\int_{\Omega_2}a(\varepsilon y)(Z_{q}-\widehat{Z}_{q})\nabla\eta_{q1}\nabla\widehat{Z}_{q}
+\int_{\Omega_4}a(\varepsilon y)
|\widehat{Z}_{q}|^2|\nabla\eta_{q2}|^2
\\[1.2mm]
\equiv&K_{21}+K_{22}+K_{23}+K_{24}+K_4.
\endaligned
$$
From (\ref{2.5}), (\ref{2.20}), (\ref{2.31}), (\ref{3.2}), (\ref{3.11}), (\ref{3.29}), (\ref{3.30}),
(\ref{3.33}), (\ref{3.48}), (\ref{3.50})  and (\ref{3.53})
we can compute
$$
\aligned
K_{21}=&-a_{q}\left(\frac{\varepsilon}{\rho_0v_0}\right)^2
\int_{(F_q^p)^{-1}\left(\left\{R<\left|\frac{\varepsilon z_q}{\rho_0 v_0}\right|\leq R+1\right\}\cap\mathbb{R}_+^2\right)}
\frac{1}{|y-q'|}a(\varepsilon y)
\mathcal{Z}_q\left(\frac{\varepsilon z_q}{\rho_0 v_0}\right)
\eta_1'\left(\frac{\varepsilon |z_q|}{\rho_0 v_0}\right)
\left[\frac{4(1+\alpha)}{c_0}+o(1)\right]dy
\\[1.2mm]
=&-\frac{c_0 a_{q}}{4(1+\alpha)}\frac{\varepsilon}{\rho_0v_0}\int_{R}^{R+1}
a(q)\eta_1'(r)\left[
\frac{4(1+\alpha)}{c_0}+O\left(\frac1{r^{2(1+\alpha)}}\right)
\right]dr\\[1mm]
=&\frac{c_0 a(q)}{p}\left(\frac{\varepsilon}{\rho_0v_0}\right)^2
\left[1+O\left(\frac1{R^{2(1+\alpha)}}\right)\right].
\endaligned
$$
Moreover by    (\ref{3.11}),  (\ref{3.33}),  (\ref{3.48}) and  (\ref{3.53}), we find
$|\nabla\eta_{q1}|=O\big(\frac{\varepsilon}{\rho_0v_0}\big)$
and $|\nabla\widehat{Z}_{q}|=O\big(\frac{\varepsilon^2}{\rho_0^2v_0^2R^{3+2\alpha}}\big)$ in $\Omega_2$.
Furthermore,   by (\ref{3.55}),
$$
\aligned
K_{22}=O\left(\frac{\varepsilon^2}{p^2\rho_0^2v_0^2R}\right),
\,\,\quad\quad\quad\,\,
K_{23}=O\left(\frac{\varepsilon^2}{p^2\rho_0^2v_0^2R}\right),
\,\,\quad\quad\quad\,\,
K_{24}=O\left(\frac{\varepsilon^2}{p\rho_0^2v_0^2R^{3+2\alpha}}\right).
\endaligned
$$
By (\ref{3.60}),
$$
\aligned
K_4=O\left(\frac{\varepsilon^2|\log d|^2}{p^2\rho_0^2v_0^2}\right).
\endaligned
$$
Hence  for
$R$ and $p$ large enough, but  $d$ small enough,
\begin{equation}\label{3.69}
\aligned
K=\frac{c_0 a(q)}{p}\left(\frac{\varepsilon}{\rho_0v_0}\right)^2
\left[1+O\left(\frac1{R^{2(1+\alpha)}}\right)\right].
\endaligned
\end{equation}
As for  the asymptotic
behavior of $L$, by virtue of
(\ref{3.2}),
(\ref{3.11}), (\ref{3.29})-(\ref{3.33}), (\ref{3.48}) and  (\ref{3.53})
we  have,  by (\ref{2.50}),
$$
\aligned
\int\limits_{(F_q^p)^{-1}\left(\left\{\left|\frac{\varepsilon z_q}{\rho_0 v_0}\right|
\leq\frac{1}{3\,p^{2\kappa}\sqrt{\rho_0 v_0}}\right\}\cap\mathbb{R}_+^2\right)}
&a(\varepsilon y)\eta_{q2}^2\left[Z_{q}-\big(1-\eta_{q1}\big)
\left(\frac{\varepsilon}{\rho_0v_0}-a_{q}G(\varepsilon y,q)\right)\right]
\\
&\times\big(1-\eta_{q1}\big)W_{\xi'}\left(\frac{\varepsilon}{\rho_0v_0}-a_{q}G(\varepsilon y,q)\right)
dy\\[0.5mm]
=O\left(\frac{\varepsilon^2}{p\rho_0^2v_0^2R^{2(1+\alpha)}}\right)&,
\endaligned
$$
and by  (\ref{3.14}),
$$
\aligned
\int\limits_{(F_q^p)^{-1}\left(\left\{\left|\frac{\varepsilon z_q}{\rho_0 v_0}\right|
\leq\frac{1}{3\,p^{2\kappa}\sqrt{\rho_0 v_0}}\right\}\cap\mathbb{R}_+^2\right)}
&a(\varepsilon y)\eta_{q2}^2\left[Z_{q}-\big(1-\eta_{q1}\big)
\left(\frac{\varepsilon}{\rho_0v_0}-a_{q}G(\varepsilon y,q)\right)\right]\\
&\times
\left[\varepsilon^2\left(
Z_{q}-\frac{\varepsilon}{\rho_0 v_0}
\right)
+\frac{\varepsilon^3}{\rho_0 v_0}\eta_{q1}\right]dy\\
=O\,\big(p\varepsilon^2\big),
\quad\qquad&
\endaligned
$$
and by  (\ref{3.13}),
$$
\aligned
\int\limits_{(F_q^p)^{-1}\left(\left\{\left|\frac{\varepsilon z_q}{\rho_0 v_0}\right|
\leq\frac{1}{3\,p^{2\kappa}\sqrt{\rho_0 v_0}}\right\}\cap\mathbb{R}_+^2\right)}
&a(\varepsilon y)\eta_{q2}^2\left[Z_{q}-\big(1-\eta_{q1}\big)
\left(\frac{\varepsilon}{\rho_0v_0}-a_{q}G(\varepsilon y,q)\right)\right]\\
&\times
\left[\Delta_{a(\varepsilon y)}Z_{q}+\left(\frac{\varepsilon}{\rho_0v_0}\right)^2
\frac{8(1+\alpha)^2\big|\frac{\varepsilon y-q}{\rho_0 v_0}\big|^{2\alpha}}
{\big(1+\big|\frac{\varepsilon y-q}{\rho_0 v_0}\big|^{2(1+\alpha)}\big)^2}Z_{q}
\right]dy\\
=O\left(\frac{\varepsilon^2}{\rho_0v_0}\right),
\quad\qquad&
\endaligned
$$
and by (\ref{2.50}),
$$
\aligned
&\int\limits_{(F_q^p)^{-1}\left(\left\{\left|\frac{\varepsilon z_q}{\rho_0 v_0}\right|
\leq\frac{1}{3\,p^{2\kappa}\sqrt{\rho_0 v_0}}\right\}\cap\mathbb{R}_+^2\right)}
a(\varepsilon y)\eta_{q2}^2\big(1-\eta_{q1}\big)
\left(\frac{\varepsilon}{\rho_0v_0}-a_{q}G(\varepsilon y,q)\right)
\left[\left(\frac{\varepsilon}{\rho_0v_0}\right)^2
\frac{8(1+\alpha)^2\big|\frac{\varepsilon y-q}{\rho_0 v_0}\big|^{2\alpha}}
{\big(1+\big|\frac{\varepsilon y-q}{\rho_0 v_0}\big|^{2(1+\alpha)}\big)^2}
-W_{\xi'}\right]Z_{q}dy\\
&=\int\limits_{(F_q^p)^{-1}\left(\left\{R\leq\left|\frac{\varepsilon z_q}{\rho_0 v_0}\right|
\leq\frac{1}{3\,p^{2\kappa}\sqrt{\rho_0 v_0}}\right\}\cap\mathbb{R}_+^2\right)}
\left(\frac{\varepsilon}{\rho_0v_0}\right)^3
\frac{8(1+\alpha)^2\big|\frac{\varepsilon y-q}{\rho_0 v_0}\big|^{2\alpha}}
{\big(1+\big|\frac{\varepsilon y-q}{\rho_0 v_0}\big|^{2(1+\alpha)}\big)^2}
\mathcal{Z}_q
\left(\frac{\varepsilon F_q^p(y)}{\rho_0 v_0}\right)\left[\,\frac1p\left(\omega_1-U_{1}-\frac{1}2U_{1}^2\right)\left(\frac{\varepsilon y-q }{\rho_0v_0}\right)
\right.\\[1.5mm]
&\,\,\,\,\quad\quad\qquad\qquad\qquad
\left.+\,
O\left(\frac1{p^2}\log^4\left|\frac{\varepsilon y-q}{\rho_0 v_0}\right|\right)\right]
\times
a_q
\left[
\frac{4(1+\alpha)}{c_0}\log\frac{|\varepsilon y-q|}{R\rho_0 v_0}+O\left(
|\varepsilon y-q|^\beta\right)
\right]O\left(1\right)
dy\\
&=\,O\left(
\frac{\varepsilon^2}{p^2\rho_0^2v_0^2R^{1+\alpha}}
\right),
\endaligned
$$
while by     (\ref{3.32}),  (\ref{3.58}) and (\ref{3.59}),
$$
\aligned
&\int\limits_{(F_q^p)^{-1}\left(\left\{\left|\frac{\varepsilon z_q}{\rho_0 v_0}\right|
>\frac{1}{3\,p^{2\kappa}\sqrt{\rho_0 v_0}}\right\}\cap\mathbb{R}_+^2\right)}
a(\varepsilon y)\widetilde{Z}_{q}\left[
\eta_{q1}\mathcal{L}(Z_{q}-\widehat{Z}_{q})
+\eta_{q2}\mathcal{L}(\widehat{Z}_{q})
\right]dy\\[0.5mm]
&\qquad\qquad\quad\qquad=\,\sum_{k=1}^m\int_{\Omega_{3,k}}
a(\varepsilon y)\eta_{q2}^2\widehat{Z}_{q}\mathcal{L}(\widehat{Z}_{q})dy
+
\int_{\widetilde{\Omega}_{3}\cup\Omega_4}
a(\varepsilon y)\eta_{q2}^2\widehat{Z}_{q}\mathcal{L}(\widehat{Z}_{q})dy
\\[0.5mm]
&\qquad\qquad\quad\qquad=\,\sum_{k=1}^mO
\left(\int_{0}^{\mu_k/(p^{2\kappa}\sqrt{\varepsilon\mu_k})}
\frac{1}{\mu_k^2}\frac{8}{\big(1+\frac{r^2}{\mu_k^2}\big)^2}
\frac{\varepsilon^2\log^2 p}{p^2\rho_0^2v_0^2}
rdr
\right)\\
&\qquad\qquad\quad\qquad\quad\,\,\,\,\,
\,+O
\left(\int_{\sqrt{\rho_0v_0}\big/(6\varepsilon p^{2\kappa})}^{8d/\varepsilon}
\left(\frac{\varepsilon}{\rho_0v_0}\right)^4
\frac{\log\varepsilon r}{p}
\left(\frac{\varepsilon r}{\rho_0v_0}\right)^{-4-2\alpha}
\frac{\log\left|\frac{\varepsilon r}{R\rho_0v_0}\right|}{p}
rdr
\right)\\[0.5mm]
&\qquad\qquad\quad\qquad=\,O\left(
\frac{\varepsilon^2\log^2p}{p^2\rho_0^2v_0^2}
\right).
\endaligned
$$
Then
$$
\aligned
L=\int\limits_{(F_q^p)^{-1}\left(\left\{\left|\frac{\varepsilon z_q}{\rho_0 v_0}\right|
\leq\frac{1}{3\,p^{2\kappa}\sqrt{\rho_0 v_0}}\right\}\cap\mathbb{R}_+^2\right)}
\,a(\varepsilon y)\eta_{q2}^2Z_{q}^2
\left[\left(\frac{\varepsilon}{\rho_0v_0}\right)^2
\frac{8(1+\alpha)^2\big|\frac{\varepsilon y-q}{\rho_0 v_0}\big|^{2\alpha}}
{\big(1+\big|\frac{\varepsilon y-q}{\rho_0 v_0}\big|^{2(1+\alpha)}\big)^2}
-W_{\xi'}\right]dy
+O\left(\frac{\varepsilon^2}{p\rho_0^2v_0^2R^{2(1+\alpha)}}\right).
\endaligned
$$
In a straightforward but tedious way, by
(\ref{2.1}), (\ref{2.16}) and (\ref{3.2}) we can compute
$$
\aligned
\int_{\mathbb{R}_{+}^2}
|z|^{2\alpha}e^{U_{1}}
\mathcal{Z}_q^2\left(
\omega_1-U_{1}-\frac12U_{1}^2
\right)dz=-4\pi(1+\alpha).
\endaligned
$$
Thus by (\ref{2.20}), (\ref{2.50}),   (\ref{3.11}), (\ref{3.48})
and (\ref{3.69}), we conclude that
for $R$ and $p$ large enough, but  $d$ small enough,
\begin{equation}\label{3.70}
\aligned
\int_{\Omega_p}a(\varepsilon y)\widetilde{Z}_{q}\mathcal{L}(\widetilde{Z}_{q})
=K+L=\frac{2c_0 a(q)}{p}\left(\frac{\varepsilon}{\rho_0v_0}\right)^2
\left[1+O\left(\frac1{R^{2(1+\alpha)}}\right)\right].
\endaligned
\end{equation}

Let us calculate
$\int_{\Omega_p}a(\varepsilon y)\widetilde{Z}_{k0}\mathcal{L}(\widetilde{Z}_{q})$
for all $k=1,\ldots,m$.
From the previous estimates of
$\mathcal{L}(\widetilde{Z}_{q})$  and $\widetilde{Z}_{k0}$,   we can easily prove that
$$
\aligned
\int_{\Omega_1}a(\varepsilon y)\widetilde{Z}_{k0}\mathcal{L}(\widetilde{Z}_{q})
=O\left(\frac{
\varepsilon R^{2(1+\alpha)}\log p}{p^2\rho_0v_0\mu_k}\right),
\,\qquad\quad\quad\,\,
\int_{\Omega_2}a(\varepsilon y)\widetilde{Z}_{k0}\mathcal{L}(\widetilde{Z}_{q})
=O\left(\frac{
\varepsilon \log p}{p^2\rho_0v_0\mu_k}\right),
\endaligned
$$
$$
\aligned
\,\,
\int_{\Omega_4}a(\varepsilon y)\widetilde{Z}_{k0}\mathcal{L}(\widetilde{Z}_{q})=
O\left(\frac{
\varepsilon |\log d|^2}{p^2\rho_0v_0\mu_k}\right),
\quad\quad\qquad\,\,\,\,\,\,\,\,
\int_{\Omega_{q}\cup\widetilde{\Omega}_{3}}a(\varepsilon y)\widetilde{Z}_{k0}\mathcal{L}(\widetilde{Z}_{q})
=O\left(\frac{
\varepsilon \log p}{p^2\rho_0v_0\mu_k}\right),
\endaligned
$$
and
$$
\aligned
\int_{\Omega_{3,l}}a(\varepsilon y)\widetilde{Z}_{k0}\mathcal{L}(\widetilde{Z}_{q})
=O\left(\frac{\varepsilon \log^2 p}{p^2\rho_0v_0\mu_k}\right)
\quad\quad\textrm{for all}\,\,\,l\neq k.
\endaligned
$$
It remains to consider the integral over $\Omega_{3,k}$. Using (\ref{3.36}) and  an integration by parts,  we obtain
$$
\aligned
\int_{\Omega_{3,k}}a(\varepsilon y)\widetilde{Z}_{k0}\mathcal{L}(\widetilde{Z}_{q})
=\int_{\Omega_{3,k}}a(\varepsilon y)\widetilde{Z}_{q}\mathcal{L}(\widetilde{Z}_{k0})
-\int_{\partial\Omega_{3,k}}a(\varepsilon y)\widehat{Z}_{k0}\frac{\partial\widehat{Z}_{q}}{\partial\nu}
+\int_{\partial\Omega_{3,k}}a(\varepsilon y)\widehat{Z}_{q}\frac{\partial\widehat{Z}_{k0}}{\partial\nu}.
\endaligned
$$
Notice that
$$
\aligned
\int_{\Omega_{3,k}}a(\varepsilon y)\widetilde{Z}_{q}\mathcal{L}(\widetilde{Z}_{k0})
=
\left(\int_{\left\{
\left|\frac{z_k}{\mu_k}\right|\leq
R
\right\}}
+
\int_{\left\{
R<
\left|\frac{z_k}{\mu_k}\right|\leq
R+1\right\}}
+
\int_{\left\{
R+1<
\left|\frac{z_k}{\mu_k}\right|\leq
\frac{1}{3\,p^{2\kappa}\sqrt{\varepsilon\mu_k}}\right\}}
\right)
a(\varepsilon y)\widetilde{Z}_{q}\mathcal{L}(\widetilde{Z}_{k0}).
\endaligned
$$
By (\ref{2.51}), (\ref{3.2}), (\ref{3.7}), (\ref{3.12}),
(\ref{3.29})-(\ref{3.32}),  (\ref{3.36})   and (\ref{3.57})  we
can compute that for any $|z_k|\leq\mu_kR$,
$$
\aligned
\mathcal{L}(\widetilde{Z}_{k0})=\mathcal{L}(Z_{k0})=
O\left(\frac{1}{p\mu_k^3}\right),
\endaligned
$$
for any $\mu_kR<|z_k|\leq\mu_k(R+1)$,
$$
\aligned
\mathcal{L}(\widetilde{Z}_{k0})
=O\left(
\frac{1}{p\mu_i^3R}
\right),
\endaligned
$$
and for any $\mu_k(R+1)<|z_k|\leq
\mu_k/(3p^{2\kappa}\sqrt{\varepsilon\mu_k})$,
$$
\aligned
\mathcal{L}(\widetilde{Z}_{k0})=
\mathcal{L}(\widehat{Z}_{k0})=\frac{1}{\mu_k^3}
\left|\frac{y-\xi'_k}{\mu_{k}}\right|^{-4}
O\left(\frac1p\log\left|\frac{y-\xi'_k}
{R\mu_k}\right|
\right).
\endaligned
$$
These, combined with the estimate of $\widehat{Z}_{q}$ in (\ref{3.32}),
imply
$$
\aligned
\int_{\Omega_{3,k}}a(\varepsilon y)\widetilde{Z}_{q}\mathcal{L}(\widetilde{Z}_{k0})
=O\left(\frac{\varepsilon\log p}{p^2\rho_0v_0\mu_k}\right).
\endaligned
$$
As on $\partial\Omega_{3,k}$, by (\ref{2.3}) and (\ref{3.32}),
$$
\aligned
\,\,\,\,\,\,\,\widehat{Z}_{k0}
=O\left(\frac{1}{\mu_k}\right),
\,\,\qquad\,\qquad\qquad\,\qquad\,\,
\widehat{Z}_{q}
=O\left(\frac{\varepsilon\log p}{p\rho_0v_0}\right),
\quad
\endaligned
$$
and
$$
\aligned
|\nabla\widehat{Z}_{k0}|
=O\left(\frac{\varepsilon^{1/2}p^{2\kappa-1}}{\mu_k^{3/2}}\right),
\,\qquad\,\,\,\qquad\,
|\nabla\widehat{Z}_{q}|
=O\left(\frac{\,\varepsilon^2 p^{\kappa-1}\,}{\rho_0v_0}\right).
\endaligned
$$
Then
$$
\aligned
\int_{\Omega_{3,k}}a(\varepsilon y)\widetilde{Z}_{k0}\mathcal{L}(\widetilde{Z}_{q})
=O\left(\frac{\varepsilon \log p}{p^2\rho_0v_0\mu_k}\right).
\endaligned
$$
By the above estimates,  we readily have
\begin{equation}\label{3.71}
\aligned
\int_{\Omega_p}a(\varepsilon y)\widetilde{Z}_{k0}\mathcal{L}(\widetilde{Z}_{q})
=O\left(\frac{\varepsilon \log^2 p}{p^2\rho_0v_0\mu_k}\right),\,\,\quad\,\,k=1,\ldots,m.
\endaligned
\end{equation}

The two expansions in (\ref{3.64}) are easy to establish as they are very
similar to the above consideration for the two expansions in (\ref{3.65}), so we  leave the detailed
proof to readers.
\end{proof}

{\bf Step 4:}
Proof of Proposition 3.1. We try with establishing  the validity of the a priori estimate
\begin{equation}\label{3.72}
\aligned
\|\phi\|_{L^{\infty}(\Omega_p)}\leq Cp\|h\|_{*}
\endaligned
\end{equation}
for any $\phi$, $c_{ij}$ solutions of
problem (\ref{3.1}) and any $h\in C^{0,\alpha}(\overline{\Omega}_p)$.
The previous step gives
$$
\aligned
\|\phi\|_{L^{\infty}(\Omega_p)}\leq Cp\left(
\|h\|_{*}+\sum_{i=1}^m\sum_{j=1}^{J_i}|c_{ij}|\cdot\|\chi_iZ_{ij}\|_*
\right)\leq Cp\left(
\|h\|_{*}+\sum_{i=1}^m\sum_{j=1}^{J_i}\mu_i|c_{ij}|
\right).
\endaligned
$$
As before, arguing by contradiction to (\ref{3.72}), we
can proceed as in Step 2 and suppose further that
\begin{equation}\label{3.73}
\aligned
\|\phi_n\|_{L^{\infty}(\Omega_{p_n})}=1,
\,\,\ \,\quad\,\,p_n\|h_n\|_{*}\rightarrow0,
\,\,\ \,\quad\,\,
p_n\sum_{i=1}^m\sum_{j=1}^{J_i}\mu_i^n|c_{ij}^n|\geq d>0
\,\,\,\quad\,\,\,
\textrm{as}\,\,\,n\rightarrow+\infty.
\endaligned
\end{equation}
For simplicity of argument we
omit the dependence on $n$ and consider the case $q\in\po$ only.
It suffices to estimate the values of the constants $c_{ij}$.
Let  us consider the cut-off function $\eta_{i2}$
defined in  (\ref{3.34})-(\ref{3.35}). Testing
(\ref{3.1}) against $a(\varepsilon y)\eta_{i2}Z_{ij}$,
$i=1,\ldots,m$ and $j=1,J_i$,
 we find
\begin{equation}\label{3.74}
\aligned
\int_{\Omega_p}
a(\varepsilon y)\phi \mathcal{L}(\eta_{i2}Z_{ij})=
\int_{\Omega_p}
a(\varepsilon y)h\eta_{i2}Z_{ij}+
\sum_{k=1}^m\sum_{t=1}^{J_k}c_{kt}
\int_{\Omega_p}\chi_kZ_{kt}\eta_{i2}Z_{ij}.
\endaligned
\end{equation}
Notice that similar to  (\ref{3.13}), by (\ref{3.2}), (\ref{3.12}) and
(\ref{3.57}) we can compute that for any $i=1,\ldots,m$ and $j=1,J_i$,
$$
\aligned
\Delta_{a(\varepsilon y)}Z_{ij}+\frac{1}{\mu_i^2}\frac{8}{\big(1+\big|\frac{y-\xi_i'}{\mu_i}\big|^2\big)^2}Z_{ij}
=O\left(\frac{\varepsilon}{\mu_i^2}
\left[1+\frac{|y-\xi'_i|}{\mu_i}\right]^{-2}\right).
\endaligned
$$
Then
$$
\aligned
\mathcal{L}(\eta_{i2}Z_{ij})=&
\eta_{i2}\mathcal{L}(Z_{ij})
-Z_{ij}\Delta_{a(\varepsilon y)}\eta_{i2}
-2\nabla\eta_{i2}\nabla Z_{ij}
\\[1mm]
=&
\left[
\frac{1}{\mu_i^2}\frac{8}{\big(1+\big|\frac{y-\xi_i'}{\mu_i}\big|^2\big)^2}
-W_{\xi'}\right]\eta_{i2}Z_{ij}
-\eta_{i2}\left[
\Delta_{a(\varepsilon y)}Z_{ij}
+
\frac{1}{\mu_i^2}\frac{8}{\big(1+\big|\frac{y-\xi_i'}{\mu_i}\big|^2\big)^2}
Z_{ij}\right]
+\varepsilon^2\eta_{i2}Z_{ij}
+O\left(\frac{\varepsilon^3}{d^3}\right)\\
\equiv
&
B_{ij}
+O\left(\frac{\varepsilon}{\mu_i^2}
\left[1+\frac{|y-\xi'_i|}{\mu_i}\right]^{-2}\right)
+O\left(\frac{\varepsilon^2}{\mu_i}
\left[1+\frac{|y-\xi'_i|}{\mu_i}\right]^{-1}\right)
+O\left(\frac{\varepsilon^3}{d^3}\right),
\endaligned
$$
where
$$
\aligned
B_{ij}=\left[
\frac{1}{\mu_i^2}\frac{8}{\big(1+\big|\frac{y-\xi_i'}{\mu_i}\big|^2\big)^2}
-W_{\xi'}\right]\eta_{i2}Z_{ij}.
\endaligned
$$
For the estimate of  $B_{ij}$, we split  $\supp(\eta_{i2})$ into the following  pieces:
$$
\aligned
\widehat{\Omega}_{k1}=\left\{\,
\left|\frac{z_k}{\mu_k}\right|
=\left|\frac{y-\xi_k'}{\mu_k}\right|
\leq
\frac{1}{3\,p^{2\kappa}\sqrt{\varepsilon\mu_k}}\right\}
\,\,\qquad\,\,\,\forall\,\,
k=1,\ldots,l,
\endaligned
$$
$$
\aligned
\widehat{\Omega}_{k1}=\Omega_p\cap(F_k^p)^{-1}\left(\left\{\,
\left|\frac{z_k}{\mu_k}\right|\leq
\frac{1}{3\,p^{2\kappa}\sqrt{\varepsilon\mu_k}}\right\}\cap\mathbb{R}_+^2\right)
\,\,\,\,\,\forall\,\,
k=l+1,\ldots,m,
\endaligned
$$
and
$$
\aligned
\widehat{\Omega}_{q}=(F_q^p)^{-1}\left(\left\{\left|\frac{\varepsilon z_q}{\rho_0 v_0}\right|
\leq\frac{1}{3\,p^{2\kappa}\sqrt{\rho_0 v_0}}\right\}\cap\mathbb{R}_+^2\right)
\qquad\,
\textrm{and}
\qquad\quad
\widehat{\Omega}_2=\supp(\eta_{i2})\setminus
\left[\bigcup_{k=1}^m\widehat{\Omega}_{k1}\cup \widehat{\Omega}_{q}\right].
\endaligned
$$
By (\ref{2.3}), (\ref{3.48}) and (\ref{3.57}) we have that for any $y\in\widehat{\Omega}_{q}$,
\begin{equation}\label{3.75}
\aligned
|y-\xi'_i|\geq |\xi_i'-q'|-|y-q'|\geq
 |\xi_i'-q'|
-
\frac{\sqrt{\rho_0 v_0}}{\,\,\varepsilon p^{2\kappa}\,\,}
>
\frac{1}{2\varepsilon p^\kappa},
\endaligned
\end{equation}
and for any $y\in\widehat{\Omega}_{1k}$ with $k\neq i$,
\begin{equation}\label{3.76}
\aligned
|y-\xi'_i|\geq |\xi_i'-\xi_k'|-|y-\xi_k'|\geq
|\xi_i'-\xi_k'|-\frac{\sqrt{\varepsilon\mu_k}}{\,\,\varepsilon  p^{2\kappa}\,\,}
>
\frac{1}{2\varepsilon p^\kappa}.
\endaligned
\end{equation}
In $\widehat{\Omega}_{1i}$, by using
(\ref{3.7}), (\ref{3.12}),
and the expansion of $W_{\xi'}$ in (\ref{2.51})
we give, for any  $i=1,\ldots,l$ and $j=1,2$,
$$
\aligned
B_{ij}=-
\frac{1}{\mu_i^3}\frac{8}{\big(1+\big|\frac{y-\xi_i'}{\mu_i}\big|^2\big)^2}
\mathcal{Z}_{j}\left(\frac{y-\xi'_i}{\mu_i}\right)
\left\{\frac1p\left(\widetilde{\omega}_1-V_{1,0}-\frac12V_{1,0}^2
\right)\left(\frac{y-\xi'_i}{\mu_i}\right)
+O\left(\frac{\log^4\big(\big|\frac{y-\xi'_i}{\mu_i}\big|+2\big)}{p^2}\right)\right\},
\endaligned
$$
and for any $i=l+1,\ldots,m$ and $j=1$,
$$
\aligned
B_{ij}=-
\frac{1}{\mu_i^3}\frac{8}{\big(1+\big|\frac{y-\xi_i'}{\mu_i}\big|^2\big)^2}
\mathcal{Z}_{j}\left(\frac{1}{\mu_i}F_i^p(y)\right)
\left\{\frac1p\left(\widetilde{\omega}_1-V_{1,0}-\frac12V_{1,0}^2
\right)\left(\frac{y-\xi'_i}{\mu_i}\right)
+O\left(\frac{\log^4\big(\big|\frac{y-\xi'_i}{\mu_i}\big|+2\big)}{p^2}\right)\right\}.
\endaligned
$$
In $\widehat{\Omega}_{q}$,  by  (\ref{3.2}), (\ref{3.7}),  (\ref{3.12})
(\ref{3.57}), (\ref{3.75}) and  the expansion of $W_{\xi'}$ in (\ref{2.50}),
$$
\aligned
B_{ij}=
\left[
O\left(\frac{\mu_i^2}{|y-\xi'_i|^4}\right)
+\left(\frac{\varepsilon}{\rho_0v_0}\right)^2
O\left(
\frac{8(1+\alpha)^2\big|\frac{\varepsilon y-q }{\rho_0v_0}\big|^{2\alpha}}
{\big(1+\big|\frac{\varepsilon y-q }{\rho_0v_0}\big|^{2(1+\alpha)}\big)^2}
\right)
\right]O\,\big(\varepsilon p^{\kappa}
\big).
\endaligned
$$
In $\widehat{\Omega}_{k1}$ with $k\neq i$, by  (\ref{3.2}), (\ref{3.7}),
(\ref{3.12}), (\ref{3.57}), (\ref{3.76})
and  the expansion of $W_{\xi'}$ in (\ref{2.51}),
$$
\aligned
B_{ij}=\left[
O\left(\frac{\mu_i^2}{|y-\xi'_i|^4}\right)+
O\left(\frac{1}{\mu_{k}^2}\frac{8}{\big(1+\big|\frac{y-\xi'_k}{\mu_{k}}\big|^2\big)^2}\right)
\right]O\,\big(\varepsilon p^{\kappa}
\big).
\endaligned
$$
In $\widehat{\Omega}_{2}$, by the estimate  of $W_{\xi'}$ in (\ref{2.49}),
$$
\aligned
B_{ij}=\left[
\left(\frac{\,\rho_0v_0\,}{\varepsilon}\right)^{2(1+\alpha)}
O\left(
\frac{1}{|y-q'|^{4+2\alpha}}\right)
+\sum_{k=1}^mO\left(
\frac{\mu_k^2}{|y-\xi'_k|^4}\right)\right]
O\left(\frac{ p^{2\kappa}\sqrt{\varepsilon\mu_i}}{\mu_i}\right).
\endaligned
$$
We denote that for any $i=1,\ldots,l$ and $j=1,2$,
$$
\aligned
\widehat{\phi}_i(z)=\phi\big(\xi_i'+\mu_iz\big),
\qquad\qquad
E_{j}(\widehat{\phi}_i)=\int_{B_{\frac{1}{3\,p^{2\kappa}\sqrt{\varepsilon\mu_k}}}\left(0\right)}
\frac{8z_j}
{(1+|z|^{2})^3}\widehat{\phi}_i\left(\widetilde{\omega}_1-V_{1,0}-\frac12V_{1,0}^2\right)
dz,
\endaligned
$$
and for any $i=l+1,\ldots,m$ and $j=1$,
$$
\aligned
\widehat{\phi}_i(z)=\phi\left((F_i^p)^{-1}(\mu_iz)\right),
\qquad\quad
E_{j}(\widehat{\phi}_i)=
\int_{B_{\frac{1}{3\,p^{2\kappa}\sqrt{\varepsilon\mu_k}}}\left(0\right)
\bigcap \mathbb{R}_+^2}
\frac{8z_j}
{(1+|z|^{2})^3}\widehat{\phi}_i\left(\widetilde{\omega}_1-V_{1,0}-\frac12V_{1,0}^2\right)
dz.
\endaligned
$$
Thus by (\ref{3.57}),
\begin{equation}\label{3.77}
\aligned
\int_{\Omega_p}
a(\varepsilon y)\phi \mathcal{L}(\eta_{i2}Z_{ij})=-\frac{1}{p\mu_i}a(\xi_i)E_{j}(\widehat{\phi}_i)
+O\left(\frac{1}{p^2\mu_i}\|\phi\|_{L^{\infty}(\Omega_p)}\right).
\endaligned
\end{equation}
On the other hand, since
$\|\eta_{i2}Z_{ij}\|_{L^{\infty}(\Omega_p)}\leq C\mu_i^{-1}$, we obtain
\begin{equation}\label{3.78}
\aligned
\int_{\Omega_p}a(\varepsilon y)h\eta_{i2}Z_{ij}=O\left(\frac{1}{\mu_i}\|h\|_{*}\right).
\endaligned
\end{equation}
Moreover,  by (\ref{3.2}), (\ref{3.7}),
(\ref{3.12}), (\ref{3.57}) and (\ref{3.75}) we conclude  that
if $1\leq k=i\leq l$,
\begin{equation}\label{3.79}
\aligned
\int_{\Omega_p}\chi_kZ_{kt}\eta_{i2}Z_{ij}=
\int_{\mathbb{R}^2}\chi\mathcal{Z}_{t} \mathcal{Z}_{j}dz=D_t\delta_{tj},
\endaligned
\end{equation}
and if $l+1\leq k=i\leq m$,
\begin{equation}\label{3.80}
\aligned
\int_{\Omega_p}\chi_kZ_{k1}\eta_{i2}Z_{i1}=
\int_{\mathbb{R}_{+}^2}
\chi\mathcal{Z}^2_{1}
\big[1+O\left(\varepsilon\mu_i|z|\right)\big] dz=\frac12D_1\big[1+O\left(\varepsilon\mu_i\right)\big],
\endaligned
\end{equation}
while if $k\neq i$,
\begin{equation}\label{3.81}
\aligned
\int_{\Omega_p}\chi_kZ_{kt}\eta_{i2}Z_{ij}=O\left(\varepsilon\mu_k p^\kappa\right).
\endaligned
\end{equation}
As a consequence, substituting  estimates (\ref{3.77})-(\ref{3.81}) into  (\ref{3.74}),  we have that
for any $i=1,\ldots,m$ and $j=1,J_i$,
$$
\aligned
D_jc_{ij}+O\left(\sum\limits_{k=1}^m\sum\limits_{t=1}^{J_k}\varepsilon\mu_k p^\kappa|c_{kt}|
\right)
= O\left(
\frac{1}{\mu_i}\|h\|_{*}
+\frac1{p\mu_i}\|\phi\|_{L^{\infty}(\Omega_p)}
\right),
\endaligned
$$
and then, by (\ref{2.31}),
\begin{equation}\label{3.82}
\aligned
\sum_{i=1}^m\sum_{j=1}^{J_i}\mu_i|c_{ij}|
= O\left(
\|h\|_{*}+\frac1{p}\|\phi\|_{L^{\infty}(\Omega_p)}
\right).
\endaligned
\end{equation}
Since $\sum\limits_{i=1}^m\sum\limits_{j=1}^{J_i}\mu_i|c_{ij}|=o\left(1\right)$,
as in contradiction arguments of
Step 2, we deduce that for any $i=1,\ldots,l$,
$$
\aligned
\widehat{\phi}_i\rightarrow
C_i\frac{|z|^2-1}{|z|^2+1}
\,\,\ \,\,\,
\,\,\textrm{uniformly in}\,\,\,C_{loc}^0(\mathbb{R}^2),
\endaligned
$$
but for any $i=l+1,\ldots,m$,
$$
\aligned
\widehat{\phi}_i\rightarrow
C_i\frac{|z|^2-1}{|z|^2+1}
\,\,\ \,\,\,
\,\,\textrm{uniformly in}\,\,\,C_{loc}^0(\mathbb{R}_{+}^2),
\endaligned
$$
with some constant $C_i$. Hence  we have a more delicate  estimate in (\ref{3.77}), because
by Lebesgue's theorem we find that for
any $i=1,\ldots,l$ and $j=1,2$,
$$
\aligned
E_{j}(\widehat{\phi}_i)
\longrightarrow
C_i\int_{\mathbb{R}^2}\frac{8z_j}{(|z|^2+1)^3}\frac{|z|^2-1}{|z|^2+1}
\left(\widetilde{\omega}_1-V_{1,0}-\frac12V_{1,0}^2\right)(|z|)dz=0,
\endaligned
$$
and for
any $i=l+1,\ldots,m$ and $j=1$,
$$
\aligned
E_{j}(\widehat{\phi}_i)
\longrightarrow
C_i\int_{\mathbb{R}^2_{+}}\frac{8z_j}{(|z|^2+1)^3}\frac{|z|^2-1}{|z|^2+1}
\left(\widetilde{\omega}_1-V_{1,0}-\frac12V_{1,0}^2\right)(|z|)dz=0.
\endaligned
$$
Therefore,
$$
\aligned
\sum_{i=1}^m\sum_{j=1}^{J_i}\mu_i|c_{ij}|=o\left(\frac1p\right)+O\left(\|h\|_{*}\right),
\endaligned
$$
which contradicts  (\ref{3.73}). So estimate (\ref{3.72})
is established and then by (\ref{3.82}),  we obtain
$$
\aligned
|c_{ij}|\leq C\frac{1}{\mu_i}\|h\|_{*}.
\endaligned
$$
\indent Now consider the Hilbert space
$$
\aligned
H_{\xi}=\left\{\phi\in H^1(\Omega_p)\left|\,\,
\int_{\Omega_p}\chi_iZ_{ij}\phi=0
\quad
\textrm{for any}\,\,\,i=1,\ldots,m,\,\,j=1,J_i;
\quad
\frac{\partial\phi}{\partial\nu}=0\quad\textrm{on}\,\,\,\po_p
\right.
\right\}
\endaligned
$$
with the norm
$\|\phi\|_{H^1(\Omega_p)}^2=\int_{\Omega_p}a(\varepsilon y)\big(|\nabla\phi|^2+\varepsilon^2\phi^2\big)$.
Equation (\ref{3.1}) is equivalent to find
$\phi\in H_\xi$ such that
$$
\aligned
\int_{\Omega_p}a(\varepsilon y)\big(\nabla\phi\nabla\psi+\varepsilon^2\phi\psi\big)
-\int_{\Omega_p}a(\varepsilon y)W_{\xi'}\phi\psi
=\int_{\Omega_p}a(\varepsilon  y)h\psi\,\,\quad\,\,\forall
\,\,\psi\in H_\xi.
\endaligned
$$
By Fredholm's alternative this is equivalent to the uniqueness of solutions to this
problem, which is guaranteed by estimate (\ref{3.72}).
Finally, for  $p\geq p_m$  fixed,
by density of $C^{0,\alpha}(\overline{\Omega}_p)$
in $(C(\overline{\Omega}_p),\,\|\cdot\|_{L^{\infty}(\Omega_p)})$,
we can approximate $h\in C(\overline{\Omega}_p)$
by smooth functions and, by (\ref{3.72}) and
elliptic regularity theory, we find that
for any $h\in C(\overline{\Omega}_p)$,
problem (\ref{3.1}) admits a unique solution which
belongs to $H^2(\Omega_p)$ and satisfies the a priori estimate
(\ref{3.15}).\qquad\qquad\qquad$\square$

\vspace{1mm}
\vspace{1mm}
\vspace{1mm}
\vspace{1mm}

\noindent{\bf Remark 3.5.}\,\,Given $h\in C(\overline{\Omega}_p)$ with
$\|h\|_*<\infty$, let $\phi$ be the solution to equation (\ref{3.1}) given by Proposition $3.1$.
Testing  (\ref{3.1}) against $a(\varepsilon y)\phi$, we obtain
$$
\aligned
\|\phi\|_{H^1(\Omega_p)}^2
=\int_{\Omega_p}a(\varepsilon y)W_{\xi'}\phi^2
+\int_{\Omega_p}a(\varepsilon  y)h\phi,
\endaligned
$$
and then,  by  (\ref{2.49}),
$$
\aligned
\|\phi\|_{H^1(\Omega_p)}\leq C\big(\|h\|_*+\|\phi\|_{L^{\infty}(\Omega_p)}\big).
\endaligned
$$

\vspace{1mm}
\vspace{1mm}
\vspace{1mm}

Now let us consider the nonlinear problem: for any points
$\xi=(\xi_1,\ldots,\xi_m)\in\mathcal{O}_p(q)$,
we find a function $\phi$ and scalars $c_{ij}$,
$i=1,\ldots,m$, $j=1,J_i$  such that
\begin{equation}\label{4.1}
\left\{\aligned
&\mathcal{L}(\phi)=-\big[
R_{\xi'}+N(\phi)
\big]
+\frac1{a(\varepsilon y)}\sum\limits_{i=1}^m\sum\limits_{j=1}^{J_i}c_{ij}\chi_i\,Z_{ij}
\,\,\ \,\textrm{in}\,\,\,\,\,\,\Omega_p,\\
&\frac{\partial\phi}{\partial\nu}=0\,\,\,\,\,\,\,
\qquad\qquad\qquad\quad\qquad\qquad\qquad\qquad
\qquad
\textrm{on}\,\,\,\,\partial\Omega_{p},\\[1mm]
&\int_{\Omega_p}\chi_i\,Z_{ij}\phi=0
\qquad\qquad\qquad\quad
\forall\,\,\,i=1,\ldots,m,\,\,\,\,j=1, J_i.
\endaligned\right.
\end{equation}

\vspace{1mm}
\vspace{1mm}

\noindent{\bf Proposition 3.6.}\,\,{\it
Let $q\in\oo$ and
$m$ be a non-negative integer.
Then there exist constants $C>0$ and  $p_m>1$ such
that for any $p>p_m$ and any points
$\xi=(\xi_1,\ldots,\xi_m)\in\mathcal{O}_p(q)$,
problem {\upshape(\ref{4.1})} admits
a unique solution
$\phi_{\xi'}$
for some coefficients $c_{ij}(\xi')$,
$i=1,\ldots,m$, $j=1,J_i$, such that
\begin{equation}\label{4.2}
\aligned
\|\phi_{\xi'}\|_{L^{\infty}(\Omega_p)}\leq\frac{C}{p^3},\,\,\quad\,\quad\,\,
\sum_{i=1}^m\sum_{j=1}^{J_i}\mu_i|c_{ij}(\xi')|\leq\frac{C}{p^4}\quad\,\quad\textrm{and}\quad\,\quad
\|\phi_{\xi'}\|_{H^1(\Omega_p)}\leq\frac{C}{p^3}.
\endaligned
\end{equation}
Furthermore, the map $\xi'\mapsto\phi_{\xi'}$ is a $C^1$-function in $C(\overline{\Omega}_p)$ and $H^1(\Omega_p)$.
}

\vspace{1mm}

\begin{proof}
Proposition $3.1$, Remarks $2.5$ and $3.5$ allow us to apply the contraction
mapping theorem
and the implicit function theorem
to find a solution for problem (\ref{4.1})
satisfying (\ref{4.2}) and the corresponding regularity
of the map $\xi'\mapsto\phi_{\xi'}$.
Since it is a standard procedure, we omit the detailed
proof here.
\end{proof}

\vspace{1mm}
\vspace{1mm}

\noindent{\bf Remark 3.7.}\,\,The function $V_{\xi'}+\phi_{\xi'}$, where $\phi_{\xi'}$ is the unique solution
of problem (\ref{4.1}) given by Proposition $3.6$, is positive in
$\overline{\Omega}_p$. In fact, we  notice
that $p^2\phi_{\xi'}\rightarrow0$ uniformly over
 $\overline{\Omega}_p$.
Furthermore, in the region
$|y-q'|\geq1/(\varepsilon p^{2\kappa})$ and
$|y-\xi'_i|\geq1/(\varepsilon p^{2\kappa})$ for each $i=1,\ldots,m$,
by (\ref{2.5}), (\ref{2.26}) and the definition of $V_{\xi'}$ in
(\ref{2.42}) we can derive that
$V_{\xi'}+\phi_{\xi'}$ is positive.
Outside this region,
we may conclude the same result from Remark $2.3$.

\vspace{1mm}

\section{Variational reduction}
Since problem (\ref{4.1}) has been solved, we just find a solution of problem (\ref{2.46})
with $m\geq1$
and hence to the original problem  (\ref{1.1}) if we find  $\xi'$ such that
the coefficient $c_{ij}(\xi')$ in (\ref{4.1}) satisfies
\begin{equation}\label{5.1}
\aligned
c_{ij}(\xi')=0\,\quad\,\,\textrm{for all}\,\,\,i=1,\ldots,m,\,\,\,j=1,J_i.
\endaligned
\end{equation}
\noindent Let
us consider the energy function $J_p$ associated to problem
(\ref{1.1}), namely
\begin{equation}\label{5.2}
\aligned
J_p(u)=\frac1{2}\int_{\Omega}a(x)(|\nabla
u|^2+u^2)dx-\frac1{p+1}\int_{\Omega}a(x)|x-q|^{2\alpha}u_{+}^{p+1}dx,
\,\,\quad\,\,\,\,\,u\in H^1(\Omega).
\endaligned
\end{equation}
For any integer $m\geq1$,
we can introduce the reduced energy
\begin{equation}\label{5.3}
\aligned
F_p(\xi)=J_p(U_{\xi}+\widetilde{\phi}_{\xi}),\,\
\ \,\,\ \ \,\ \,\,\,\xi=(\xi_1,\ldots,\xi_m)\in\mathcal{O}_p(q),
\endaligned
\end{equation}
where $U_\xi$ is our approximate solution defined in  (\ref{2.18}) and
\begin{equation}\label{5.4}
\aligned
\widetilde{\phi}_{\xi}(x)=\varepsilon^{-2/(p-1)}\phi_{\xi'}(\varepsilon^{-1}x),
\,\quad\,x\in\Omega,
\endaligned
\end{equation}
with $\phi_{\xi'}$ the unique solution to problem (\ref{4.1}) given by
Proposition $3.6$. Then we obtain that critical points of $F_p$ correspond to solutions of (\ref{5.1})
for large $p$. That is:

\vspace{1mm}
\vspace{1mm}
\vspace{1mm}
\vspace{1mm}

\noindent{\bf Proposition 4.1.}\,\,{\it For any integer $m\geq1$,
the function $F_p:\mathcal{O}_p(q)\mapsto\mathbb{R}$
is of class $C^1$.
Moreover, for all $p$ sufficiently large,
if  $D_{\xi}F_p(\xi)=0$,  then $\xi'=\xi/\varepsilon$ satisfies {\upshape (\ref{5.1})}.
}

\vspace{1mm}

\begin{proof}
From the result obtained in Proposition $3.6$ and the
definition of function $U_\xi$ we have clearly that
for any integer $m\geq1$, the function
$F_p:\mathcal{O}_p(q)\rightarrow\mathbb{R}$
is of class $C^1$
since the map
$\xi\mapsto\widetilde{\phi}_{\xi}$ is a $C^1$-map into $H^1(\Omega)$.

Recalling the definition of  $I_p$ in (\ref{2.43}) and
making a change of variable, we give
\begin{equation}\label{5.5}
\aligned
F_p(\xi)=J_p\big(U_{\xi}+\widetilde{\phi}_{\xi}\big)=\varepsilon^{-4/(p-1)}I_p\big(V_{\xi'}+\phi_{\xi'}\big).
\endaligned
\end{equation}
Assume that  $\phi_{\xi'}$ solves problem (\ref{4.1})
and $D_{\xi}F_p(\xi)=0$. Then  we have that for any $k=1,\ldots,m$ and $t=1,J_k$,
\begin{eqnarray}\label{5.6}
&&0
=I'_p\big(V_{\xi'}+\phi_{\xi'}\big)\partial_{(\xi'_k)_t}\big(V_{\xi'}+\phi_{\xi'}\big)\nonumber\\[1mm]
&&\,\,\,\,\,=\sum\limits_{i=1}^m\sum\limits_{j=1}^{J_i}c_{ij}(\xi')\int_{\Omega_p}\chi_i Z_{ij}
\partial_{(\xi'_k)_t}V_{\xi'}
-\sum\limits_{i=1}^m\sum\limits_{j=1}^{J_i}c_{ij}(\xi')\int_{\Omega_p}\phi_{\xi'}
\partial_{(\xi'_k)_t}\big(\chi_i Z_{ij}\big).
\end{eqnarray}
Note that $V_{\xi'}(y)=\varepsilon^{2/(p-1)}U_\xi(\varepsilon y)$.
From  (\ref{2.26}),  (\ref{2.27}),  (\ref{2.29})
and the definition of $U_\xi$ in (\ref{2.18}), we obtain
$$
\aligned
\partial_{(\xi'_k)_t}V_{\xi'}(y)
=&\,\sum_{i=1}^{m}\frac{\varepsilon^{2/(p-1)}}{\,\gamma\mu_i^{2/(p-1)}
\big|\xi_i-q\big|^{2\alpha/(p-1)}\,}
\left\{\partial_{(\xi'_k)_t}\left[V_{\delta_i,\xi_i}(\varepsilon y)
+\frac1p\widetilde{\omega}_1\left(\frac{y-\xi'_i}{\mu_i}\right)
+\frac1{p^2}\widetilde{\omega}_2\left(\frac{y-\xi'_i}{\mu_i}\right)
\right.\right.\\
&\left.
+
\left.
\gamma\mu_i^{2/(p-1)}
\big|\xi_i-q\big|^{2\alpha/(p-1)}H_i(\varepsilon y)\right]
-\frac{2\,\partial_{(\xi'_k)_t}\log\big(\mu_i\,|\xi_i-q|^{\alpha}\big)}{p-1}\big[\,
p+O\left(\log p\right)\big]
\right\}\\
&+
\frac{\varepsilon^{2/(p-1)}}{\,\,\gamma\mu_0^{2/(p-1)}\,\,}
\left\{\partial_{(\xi'_k)_t}\left[U_{\delta_0}(\varepsilon y-q)
+\frac1p\omega_1\left(\frac{\varepsilon y-q}{\rho_0v_0}\right)
+\frac1{p^2}\omega_2\left(\frac{\varepsilon y-q}{\rho_0v_0}\right)
+\gamma\mu_0^{2/(p-1)}H_0(\varepsilon y)\right]\right.\\[0.5mm]
&\left.-
\frac{\,2\,\partial_{(\xi'_k)_t}\log\mu_0\,}{p-1}\big[\,
p+O\left(\log p\right)\big]
\right\}.
\endaligned
$$
Using the fact that
$|\partial_{(\xi'_k)_t}\log\mu_i|=O\left(\varepsilon p^{\kappa}\right)$
for any $i=0,1,\ldots,m$, we have that by
(\ref{2.1}), (\ref{2.5}) and (\ref{3.2}),
$$
\aligned
\partial_{(\xi'_k)_t}U_{\delta_0}(\varepsilon y-q)=
O\left(\varepsilon p^{\kappa}\right),
\qquad\qquad
\partial_{(\xi'_k)_t}V_{\delta_i,\xi_i}(\varepsilon y)=
\frac{4\delta_{ki}}{\mu_i}
\mathcal{Z}_{t}\left(
\frac{y-\xi'_i}{\mu_i}
\right)
+O\left(\varepsilon p^{\kappa}\right),
\endaligned
$$
and for each $j=1,2$,
by (\ref{2.12})-(\ref{2.13}),
$$
\aligned
\partial_{(\xi'_k)_t}\omega_{j}
\left(\frac{\varepsilon y-q}{\rho_0v_0}\right)=
O\left(\varepsilon p^{\kappa}\right),
\qquad
\partial_{(\xi'_k)_t}\widetilde{\omega}_{j}\left(\frac{y-\xi'_i}{\mu_i}\right)=
-\frac{\delta_{ki}}{\mu_i}\left[\widetilde{C}_j\mathcal{Z}_{t}\left(
\frac{y-\xi'_i}{\mu_i}
\right)+
O\left(\frac{\mu_i^2}{|y-\xi'_i|^2+\mu_i^2}\right)
\right]
+O\left(\varepsilon p^{\kappa}\right),
\endaligned
$$
where $\delta_{ki}$ denotes the Kronecker's symbol.
Moreover,  similar to the proof in Lemma $2.1$, we can  prove that
$$
\aligned
\partial_{(\xi'_k)_t}\left[\gamma\mu_0^{2/(p-1)}H_0(\varepsilon y)\right]=
O\left(\varepsilon p^{\kappa}\right),
\qquad\qquad
\partial_{(\xi'_k)_t}
\left[\gamma\mu_i^{2/(p-1)}
\big|\xi_i-q\big|^{2\alpha/(p-1)}H_i(\varepsilon y)\right]=
O\left(\varepsilon p^{\kappa}\right).
\endaligned
$$
Then
\begin{eqnarray}\label{5.8}
\partial_{(\xi'_k)_t}V_{\xi'}(y)=\frac{\varepsilon^{2/(p-1)}}
{\,\gamma\mu_k^{2/(p-1)}
\big|\xi_k-q\big|^{2\alpha/(p-1)}\,}
\left\{
\frac{4}{\mu_k}\mathcal{Z}_{t}\left(
\frac{y-\xi'_k}{\mu_k}
\right)
+
O\left(
\frac1{p\mu_k}\right)
\right\}+
O\left(\varepsilon p^{\kappa-1}\right).
\end{eqnarray}
On the other hand,  by (\ref{3.2}), (\ref{3.7}),  (\ref{3.12})  and (\ref{3.57})
we can compute
\begin{eqnarray}\label{5.7}
\big|\partial_{(\xi'_k)_t}\big(\chi_i Z_{ij}\big)\big|
=O\left(
\frac1{\mu_i}\varepsilon p^{\kappa}
+\frac1{\mu_i^2}\delta_{ki}
\right).
\end{eqnarray}
Consequently,    (\ref{5.6}) can be written as, for
each $k=1,\ldots,m$  and  $t=1,J_k$,
$$
\aligned
&\sum_{i,j}\frac{\varepsilon^{2/(p-1)}c_{ij}(\xi')}
{\,\gamma\mu_k^{2/(p-1)}
\big|\xi_k-q\big|^{2\alpha/(p-1)}\,}
\left\{
\delta_{ki}\left[
\frac{c_i}{2\pi}\int_{\mathbb{R}^2}\chi \mathcal{Z}_j\mathcal{Z}_t
+O\left(\frac1p\right)
\right]
+(1-\delta_{ki})O\left(
\frac{\mu_i}{|\xi'_i-\xi'_k|}
\right)
\right\}\\[1mm]
&+\sum_{ij}|c_{ij}(\xi')|
\big\{\,
O\left(\mu_i
\varepsilon p^{\kappa-1}
\right)+\|\phi_{\xi'}\|_{L^{\infty}(\Omega_p)}\,
O\left(\mu_i\varepsilon p^{\kappa}
+\delta_{ki}
\right)
\big\}=0,
\endaligned
$$
and then, by (\ref{2.3}), (\ref{2.5}), (\ref{2.31}) and (\ref{4.2}),
$$
\aligned
\frac{c_kc_{kt}(\xi')}{p^{p/(p-1)}\mu_k^{2/(p-1)}
\big|\xi_k-q\big|^{2\alpha/(p-1)}}
\int_0^{R_0+1}\chi(r)\frac{r^{3}}{(\,1+r^{2})^2}
dr
+
\sum_{i=1}^m\sum_{j=1}^{J_i}\big|c_{ij}(\xi')\big|O\left(
\frac{\delta_{ki}}{p^2}
+
\varepsilon\mu_i p^{\kappa-1}
+
\varepsilon p^{\kappa}
\right)=0,
\endaligned
$$
which  implies $c_{kt}(\xi')=0$
for each $k=1,\ldots,m$  and  $t=1,J_k$.
\end{proof}

\vspace{1mm}

Moreover, in order to solve for  critical points of  $F_p$, we need to give
the following reduced energy expansion.

\vspace{1mm}
\vspace{1mm}
\vspace{1mm}

\noindent{\bf Proposition 4.2.}\,\,{\it
Let $q\in\oo$ and $m$ be a positive integer.
With the choices for the parameters $\mu_0$ and $\mu_i$, $i=1,\ldots,m$, respectively given by {\upshape(\ref{2.28})}
and  {\upshape(\ref{2.30})},
there exists $p_m>1$ such that for any $p>p_m$ and any
points $\xi=(\xi_1,\ldots,\xi_m)\in\mathcal{O}_p(q)$,
the following expansion uniformly holds
\begin{eqnarray}\label{6.1}
F_p(\xi)=
\frac{e}{2p}
\sum_{i=1}^mc_ia(\xi_i)\left\{
1
+\frac{\widetilde{\mathcal{K}}+2}p
-\frac{2\log p}{p}
-\frac{4\alpha\log|\xi_i-q|}{p}
-\frac{1}{p}\left[c_i H(\xi_i,\xi_i)
+c_0 G(\xi_i,q)
+\sum_{k\neq i}^m
c_k G(\xi_i,\xi_k)\right]
\right\}
&&\nonumber\\
+
\frac{e}{2p}
c_0a(q)\left\{
1+\frac{\mathcal{K}+2}p-\frac{2\log p}{p}
-\frac{1}{p}\left[c_0 H(q,q)
+\sum_{k=1}^m
c_k G(q,\xi_k)
\right]\right\}
+O\left(\frac{\log^2p}{p^3}\right),
\qquad\qquad\qquad\quad\,\,\,\,
&&
\end{eqnarray}
where  the coefficients $c_0$ and $c_k$, $k=1,\ldots,m$, are defined in {\upshape(\ref{2.20})} and
$$
\aligned
\mathcal{K}=\frac{1}{8\pi(1+\alpha)}\int_{\mathbb{R}^2}\left[
\frac{8(1+\alpha)^2|z|^{2\alpha}}{(1+|z|^{2(1+\alpha)})^2}U_{1}(z)-\Delta\omega_1(z)
\right]dz,
\quad\,\,\quad
\widetilde{\mathcal{K}}=\frac{1}{8\pi}\int_{\mathbb{R}^2}\left[
\frac{8}{(1+|\tilde{z}|^2)^2}V_{1,0}(\tilde{z})-\Delta\widetilde{\omega}_1(\tilde{z})
\right]d\tilde{z}.
\endaligned
$$
}

\begin{proof}
First of all, multiply   the first equation in (\ref{4.1})  by $a(\varepsilon y)(V_{\xi'}+\phi_{\xi'})$
and   integrate by parts to give
$$
\aligned
\int_{\Omega_p}a(\varepsilon y)\left[|\nabla
(V_{\xi'}+\phi_{\xi'})|^2+\varepsilon^2(V_{\xi'}+\phi_{\xi'})^2\right]=\int_{\Omega_p}a(\varepsilon y)|\varepsilon y-q|^{2\alpha}(V_{\xi'}+\phi_{\xi'})^{p+1}
+\sum_{i=1}^m\sum_{j=1}^{J_i}c_{ij}(\xi')\int_{\Omega_p}\chi_iZ_{ij}V_{\xi'}.
\endaligned
$$
Since $V_{\xi'}$ is a uniformly bounded function, by  (\ref{4.2}) we get
$$
\aligned
\int_{\Omega_p}a(\varepsilon y)\left[|\nabla
(V_{\xi'}+\phi_{\xi'})|^2+\varepsilon^2(V_{\xi'}+\phi_{\xi'})^2\right]
=\int_{\Omega_p}a(\varepsilon y)|\varepsilon y-q|^{2\alpha}(V_{\xi'}+\phi_{\xi'})^{p+1}
+O\left(\frac1{p^4}\right)
\endaligned
$$
uniformly for any points
$\xi=(\xi_1,\ldots,\xi_m)\in\mathcal{O}_p(q)$.
Then by (\ref{2.43}) and (\ref{5.5}) we have
\begin{eqnarray*}
F_p(\xi)=\left(\frac1{2}-\frac{1}{p+1}\right)\varepsilon^{-4/(p-1)}\int_{\Omega_p}a(\varepsilon y)\left[|\nabla
(V_{\xi'}+\phi_{\xi'})|^2+\varepsilon^2(V_{\xi'}+\phi_{\xi'})^2\right]dy+O\left(\frac1{p^5}\right)
\,\,\,\,\,\quad
&&\nonumber\\
=\left(\frac1{2}-\frac{1}{p+1}\right)\varepsilon^{-4/(p-1)}\left\{\int_{\Omega_p}a(\varepsilon y)\big(|\nabla
V_{\xi'}|^2+\varepsilon^2V_{\xi'}^2\big)dy
\right.
\qquad\qquad\qquad\quad\,\,
\qquad\qquad\,\,\,\quad\,
&&\nonumber\\
\left.+2\int_{\Omega_p}a(\varepsilon y)\big(\nabla
V_{\xi'}\nabla\phi_{\xi'}+\varepsilon^2V_{\xi'}\phi_{\xi'}\big)dy
+\int_{\Omega_p}a(\varepsilon y)\big(|\nabla
\phi_{\xi'}|^2+\varepsilon^2\phi_{\xi'}^2\big)dy\right\}+O\left(\frac1{p^5}\right)
\quad
&&\nonumber\\
=\left(\frac1{2}-\frac{1}{p+1}\right)\int_{\Omega}a(x)\big(|\nabla
U_{\xi}|^2+U_{\xi}^2\big)dx
+O\left(\frac1{p^3}\left|\int_{\Omega}a(x)\big(|\nabla
U_{\xi}|^2+U_{\xi}^2\big)dx\right|^{1/2}+\frac1{p^5}\right).
&&
\end{eqnarray*}
Let us expand the leading term  $\int_{\Omega}a(x)(|\nabla
U_\xi|^2+U_\xi^2)$:
in view of (\ref{2.5})-(\ref{2.7}), (\ref{2.18})-(\ref{2.20}),
(\ref{2.26})-(\ref{2.27}) and (\ref{2.29}) we obtain
$$
\aligned
\int_{\Omega}&a(x)\big(|\nabla
U_{\xi}|^2+U_{\xi}^2\big)dx=\int_{\Omega}a(x)(-\Delta_{a}
U_\xi+U_{\xi})U_\xi dx\\
=&\sum_{i=1}^m\frac{1}{\gamma\mu_i^{2/(p-1)}\big|\xi_i-q\big|^{2\alpha/(p-1)}}
\int_{\Omega\bigcap B_{\frac1{p^{2\kappa}}}(\xi_i)}
\frac{a(x)}{\delta_i^2}\left[e^{V_{1,0}\big(\frac{x-\xi_i}{\delta_i}\big)}
-\frac{1}{p}\Delta\widetilde{\omega}_1\left(
\frac{x-\xi_i}{\delta_i}\right)
-\frac{1}{p^2}\Delta\widetilde{\omega}_2\left(
\frac{x-\xi_i}{\delta_i}\right)
\right]U_\xi dx \\
&+\frac{1}{\gamma\mu_0^{2/(p-1)}}
\int_{\Omega\bigcap B_{\frac1{p^{2\kappa}}}(q)}
\frac{a(x)}{\delta_0^2}\left[\left|\frac{x-q}{\delta_0}\right|^{2\alpha}
e^{U_{1}\big(\frac{x-q}{\delta_0}\big)}
-\frac{1}{p}\Delta\omega_1\left(
\frac{x-q}{\delta_0}\right)
-\frac{1}{p^2}\Delta\omega_2\left(
\frac{x-q}{\delta_0}\right)
\right]
U_\xi dx
+o\left(\frac1{p^5}\right)\\
=&\sum_{i=1}^m\frac{1}{\gamma^2\mu_i^{4/(p-1)}\big|\xi_i-q\big|^{4\alpha/(p-1)}}
\int_{\big(\frac{1}{\delta_i}(\Omega-\{\xi_i\})\big)\bigcap B_{\frac1{\delta_ip^{2\kappa}}}(0)}
a(\xi_i+\delta_i\tilde{z})\left[\frac{8}{(1+|\tilde{z}|^2)^2}
-\frac{1}{p}\Delta\widetilde{\omega}_1(\tilde{z})
-\frac{1}{p^2}\Delta\widetilde{\omega}_2(\tilde{z})
\right]\\
&\times\left[
p+V_{1,0}(\tilde{z})+\frac1p\widetilde{\omega}_1(\tilde{z})+\frac1{p^2}\widetilde{\omega}_2(\tilde{z})+O\left(
\delta_i^\beta|\tilde{z}|^\beta
+
p^{2\kappa(1+\alpha)-1}\delta_0^{1+\alpha}
+\sum_{k=0}^m\delta_k^{\beta/2}
\right)\right]d\tilde{z}\\
&+\frac{1}{\gamma^2\mu_0^{4/(p-1)}}
\int_{\big(\frac{1}{\delta_0}(\Omega-\{q\})\big)\bigcap B_{\frac1{\delta_ip^{2\kappa}}}(0)}
a(q+\delta_0z)\left[\frac{8(1+\alpha)^2|z|^{2\alpha}}{(1+|z|^{2(1+\alpha)})^2}
-\frac{1}{p}\Delta\omega_1(z)
-\frac{1}{p^2}\Delta\omega_2(z)
\right]\\
&\times\left[
p+U_{1}(z)+\frac1p\omega_1(z)+\frac1{p^2}\omega_2(z)+O\left(
\delta_0^\beta|z|^\beta
+\sum_{k=0}^m\delta_k^{\beta/2}
\right)\right]dz
+o\left(\frac1{p^5}\right)\\
=&\sum_{i=1}^m\frac{c_ia(\xi_i)}{\gamma^2\mu_i^{4/(p-1)}\big|\xi_i-q\big|^{4\alpha/(p-1)}}
\left[\,p
+\widetilde{\mathcal{K}}
+O\left(\frac1p\right)
\right]
+\frac{c_0a(q)}{\gamma^2\mu_0^{4/(p-1)}}
\left[\,p
+\mathcal{K}
+O\left(\frac1p\right)
\right]
+o\left(\frac1{p^5}\right).
\endaligned
$$
Recalling  that $\gamma=p^{p/(p-1)}e^{-p/(2p-2)}$, we find
$$
\aligned
\frac{1}{\gamma^2}=\frac{e}{p^2}\left[
1-\frac{2\log p}{p}+\frac{1}{p}+O\left(
\frac{\log^2 p}{p^2}
\right)
\right],
\,\qquad\qquad\qquad\,
\frac{1}{\,\mu_0^{4/(p-1)}}=1-\frac{4\log\mu_0}{p}+O\left(\frac{\log^2\mu_0}{p^2}\right),
\endaligned
$$
and for each $i=1,\ldots,m$,
$$
\aligned
\frac{1}{\,\mu_i^{4/(p-1)}\big|\xi_i-q\big|^{4\alpha/(p-1)}}
=1-\frac{4\log\big(\mu_i|\xi_i-q|^\alpha\big)}{p}+O\left(\frac{\log^2\mu_i+\log^2|\xi_i-q|}{p^2}\right).
\endaligned
$$
Thus by (\ref{2.3}) and (\ref{2.31}),
$$
\aligned
\frac{1}{\gamma^2\mu_0^{4/(p-1)}\,}&=
\frac{e}{p^2}\left[
1-\frac{2\log p}{p}-\frac{4\log\mu_0}{p}
+\frac{1}{p}+O\left(
\frac{\log^2p}{p^2}
\right)
\right],\\
\frac{1}{\gamma^2\mu_i^{4/(p-1)}\big|\xi_i-q\big|^{4\alpha/(p-1)}\,}=
\frac{e}{p^2}&\left[
1-\frac{2\log p}{p}-\frac{4\log\mu_i}{p}-\frac{4\alpha\log|\xi_i-q|}{p}
+\frac{1}{p}+O\left(
\frac{\log^2p}{p^2}
\right)
\right],
\quad
i=1,\ldots,m.
\endaligned
$$
Hence
$$
\aligned
F_p(\xi)=&\frac{e}{2p}
\sum_{i=1}^mc_ia(\xi_i)\left\{
1-\frac{2\log p}{p}-\frac{4\log\mu_i}{p}
-\frac{4\alpha\log|\xi_i-q|}{p}
+\frac{\widetilde{\mathcal{K}}-1}p
\right\}\\
&+\frac{e}{2p}
c_0a(q)\left\{
1-\frac{2\log p}{p}-\frac{4\log\mu_0}{p}
+\frac{\mathcal{K}-1}p
\right\}
+O\left(\frac{\log^2p}{p^3}\right),
\endaligned
$$
which, together with the expansions of $\mu_0$, $\mu_i$ in
(\ref{2.36})-(\ref{2.37}),  implies that expansion (\ref{6.1}) holds.
\end{proof}

\vspace{1mm}

\section{Proofs of theorems}
\noindent {\bf Proof of Theorem 1.1.}\,\,
According to Proposition  4.1, the function
$u_{p }=U_{\xi}+\widetilde{\phi}_{\xi}$
is a solution to problem (\ref{1.1}) if we adjust
$\xi=(\xi_1,\ldots,\xi_m)\in\mathcal{O}_p(q)$ with $q\in\Omega$ so that it is
a critical point of $F_p$ defined in (\ref{5.3}).
For this aim, let us claim that for any integer $m\geq1$ and any $p>1$ large enough, the maximization problem
$$
\aligned
\max\limits_{(\xi_1,\ldots,\xi_m)\in\overline{\mathcal{O}}_p(q)}
F_p(\xi_1,\ldots,\xi_m)
\endaligned
$$
has a solution $\xi^p=(\xi_1^p,\ldots,\xi_m^p)\in\mathcal{O}_p^o(q)$, i.e., the interior of $\mathcal{O}_p(q)$.
Once this claim  is proven, we can easily get the qualitative properties of
solutions of problem (\ref{1.1}) as predicted in Theorem 1.1.

Let $\xi^p=(\xi^p_1,\ldots,\xi^p_m)$
be the maximizer of $F_p$ over $\overline{\mathcal{O}}_p(q)$. We need to prove that $\xi^p$
belongs to the interior of $\mathcal{O}_p(q)$.
First, we obtain a lower bound for $F_p$ over $\overline{\mathcal{O}}_p(q)$.
Let
$$
\aligned
\xi^0_i=q+
\frac{1}{\sqrt{p}}\widehat{\xi}_i,
\quad
i=1,\ldots,m,
\endaligned
$$
where $\widehat{\xi}=(\widehat{\xi}_1,\ldots,\widehat{\xi}_m)$ forms
a $m$-regular polygon in $\mathbb{R}^2$.
Obviously,
$\xi^0=(\xi^0_1,\ldots,\xi^0_m)\in\mathcal{O}_p(q)$
since $q\in\Omega$ and  $\kappa>1$.
Using (\ref{6.1}) and the fact that $q\in\Omega$ is a strict local maximum point of $a(x)$,  we obtain
\begin{eqnarray}\label{7.1}
\max\limits_{\xi\in\overline{\mathcal{O}}_p(q)}
F_p(\xi)\geq
F_p(\xi^0)\geq
\frac{e}{2p^2}
\big\{
8\pi (m+1+\alpha)pa(q)-16\pi(m+1)(m+1+\alpha)a(q)\log p
+O(1)
\big\}.
\end{eqnarray}
Next, we suppose $\xi^p=(\xi_1^p,\ldots,\xi_m^p)\in\partial\mathcal{O}_p (q)$.
Then there  exist three possibilities:\\
C1. There exists an $i_0$ such that
$\xi^p_{i_0}\in\partial B_d(q)$, in which case,
$a(\xi^p_{i_0})<a(q)-d_0$ for some  $d_0>0$ independent of $p$;\\
C2. There exist indices $i_0$, $j_0$,
$i_0\neq j_0$ such that
$|\xi_{i_0}^p-\xi_{j_0}^p|=p^{-\kappa}$;\\
C3. There exists an  $k_0$
such that
$|\xi_{k_0}^p-q|=p^{-\kappa}$.\\
For the first case,   we have
\begin{eqnarray}\label{7.2}
\max\limits_{\xi\in\overline{\mathcal{O}}_p(q)}
F_p(\xi)
\leq
\frac{e}{2p^2}
\big\{8\pi p\big[(m+1+\alpha)a(q)-d_0\big]+O\left(\log p\right)\big\},
\end{eqnarray}
which contradicts to (\ref{7.1}).
This shows that $a(\xi_i^p)\rightarrow a(q)$. By the condition over $a$,
we get $\xi_i^p\rightarrow q$ for all $i=1,\ldots,m$.\\
For the second case, we have
\begin{eqnarray}\label{7.3}
\max\limits_{\xi\in\overline{\mathcal{O}}_p(q)}
F_p(\xi)
\leq
\frac{e}{2p^2}
\left\{
8\pi(m+1+\alpha)(p-2\log p)a(q)
+32\pi \big(a(\xi^p_{i_0})
+
a(\xi^p_{j_0})\big)
\log|\xi^p_{i_0}-\xi^p_{j_0}|
+O(1)
\right\}
&&\nonumber\\
\leq
\frac{e}{2p^2}
\left\{
8\pi(m+1+\alpha)(p-2\log p)a(q)
-32\pi\kappa \big(a(\xi^p_{i_0})
+
a(\xi^p_{j_0})\big)\log p
+O(1)
\right\}.
\qquad\,\,\,\,\,
&&
\end{eqnarray}
For the last case, we have
\begin{eqnarray}\label{7.4}
\max\limits_{\xi\in\overline{\mathcal{O}}_p(q)}
F_p(\xi)
\leq
\frac{e}{2p^2}
\left\{
8\pi(m+1+\alpha)(p-2\log p)a(q)
+32\pi \big[a(\xi^p_{k_0})
+(1+\alpha)
a(q)\big]
\log|\xi^p_{k_0}-q|
+O(1)
\right\}
&&\nonumber\\
\leq
\frac{e}{2p^2}
\left\{
8\pi(m+1+\alpha)(p-2\log p)a(q)
-32\pi\kappa \big[a(\xi^p_{k_0})
+(1+\alpha)
a(q)\big]\log p
+O(1)
\right\}.
\quad\,\,\,\,\,\,
&&
\end{eqnarray}
Combining (\ref{7.3})-(\ref{7.4}) with (\ref{7.1}), we find
\begin{equation}\label{7.5}
\aligned
32\pi\kappa\max\big\{\,
a(\xi^p_{i_0})+a(\xi^p_{j_0}),
\,\,
a(\xi^p_{k_0})+(1+\alpha)a(q)\big\}\log p
\leq
16\pi m(m+1+\alpha)a(q)\log p
+O(1),
\endaligned
\end{equation}
which  is impossible by the
choice of $\kappa$ in (\ref{2.4}).
\,\qquad\qquad\qquad\qquad\qquad\qquad\qquad\qquad
\qquad\qquad\qquad\qquad\qquad\quad\,\,\,$\square$

\vspace{1mm}
\vspace{1mm}
\vspace{1mm}
\vspace{1mm}

\noindent {\bf Proof of Theorem 1.2.}
By Proposition $4.1$ we have to find
a critical point $\xi^p=(\xi^p_1,\ldots,\xi^p_m)\in\big(B_d(q)\cap\Omega\big)^l
\times\big(B_d(q)\cap\po\big)^{m-l}$
of  $F_p$ such  that points $\xi^p_1,\ldots,\xi^p_m$
accumulate to the boundary point $q$.
From  (\ref{1.3}), (\ref{2.20}), (\ref{6.1}), Lemma A.1 and the fact that
$a(q)G(q,\xi_i)=a(\xi_i)G(\xi_i,q)$
and $a(\xi_i)G(\xi_i,\xi_k)=a(\xi_k)G(\xi_k,\xi_i)$ for all
$i,k=1,\ldots,m$ with $i\neq k$,
we have that $F_p$
reduces to
\begin{eqnarray}\label{7.6}
F_p(\xi)=\frac{e}{2p^2}
\left\{8\pi\sum_{i=1}^la(\xi_i)\left[
p-2\log p
-4\alpha\log|\xi_i-q|
-8\pi H(\xi_i,\xi_i)
-8\pi \sum_{k=1,\,k\neq i}^l
G(\xi_i,\xi_k)
\right]
\quad
\qquad\qquad\qquad
\right.
&&\nonumber\\
-64\pi^2\sum_{i=1}^l\sum_{k=l+1}^m a(\xi_k)
G(\xi_k,\xi_i)+4\pi\sum_{i=l+1}^ma(\xi_i)\left[
p-2\log p-4\alpha\log|\xi_i-q|+4\sum_{k=l+1,\,k\neq i}^m
\log|\xi_i-\xi_k|
\right]
&&\nonumber\\
\left.
+8\pi(1+\alpha)
a(q)\left[
p-2\log p
-16\pi \sum_{i=1}^l
G(q,\xi_i)
+ 8\sum_{i=l+1}^m
\log|\xi_i-q|
\right]+O\left(1\right)
\right\}\,\,\,
\qquad\qquad\qquad\qquad\qquad
&&
\end{eqnarray}
$C^0$-uniformly in
$\mathcal{O}_{p}(q)$.
Here we claim that for any integers $m\geq1$, $0\leq l\leq m$ and any $p>1$ large enough, the maximization problem
$$
\aligned
\max\limits_{(\xi_1,\ldots,\xi_m)\in\overline{\mathcal{O}}_p(q)}
F_p(\xi_1,\ldots,\xi_m)
\endaligned
$$
has a solution $\xi^p=(\xi_1^p,\ldots,\xi_m^p)\in\mathcal{O}_p^o(q)$, i.e., the interior of $\mathcal{O}_p(q)$.
Once this claim  is proven, we can easily get the qualitative properties of
solutions of problem (\ref{1.1}) as predicted in Theorem 1.2.

Let $\xi^p=(\xi^p_1,\ldots,\xi^p_m)$ be the maximizer of $F_p$
over $\overline{\mathcal{O}}_p(q)$. We  need to prove that $\xi^p$
lies in the interior of $\mathcal{O}_p(q)$.
First, we obtain a lower bound for $F_p$ over $\overline{\mathcal{O}}_p(q)$.
Let us consider a smooth change of variables
$$
\aligned
H_{q}^p(y)=
e^{p/2}H_{q}(e^{-p/2} y),
\endaligned
$$
where $H_{q}:B_d(q)\mapsto\mathcal{M}$
is a  diffeomorphism and
$\mathcal{M}$  is an open neighborhood of the origin such that
$H_{q}(B_d(q)\cap\Omega)=\mathcal{M}\cap\mathbb{R}_+^2$
and
$H_{q}(B_d(q)\cap\partial\Omega)=\mathcal{M}\cap\partial\mathbb{R}_+^2$.
Let
$$
\aligned
\xi^0_i=q-\frac{t_i}{\sqrt{p}}\nu(q),
\quad i=1,\ldots,l,
\qquad\quad\textrm{but}\quad\qquad
\xi^0_i=e^{-p/2}(H_{q}^p)^{-1}\left(
\frac{e^{p/2}}{\sqrt{p}}\widehat{\xi}_i^0
\right),
\quad i=l+1,\ldots,m,
\endaligned
$$
where $t_i>0$  and
$\widehat{\xi}_i^0\in \mathcal{M}\cap\partial\mathbb{R}_+^2$  satisfy
$t_{i+1}-t_i=\sigma$,
$|\widehat{\xi}_i^0-\widehat{\xi}^0_{i+1}|=\sigma$
for all $\sigma>0$ sufficiently
small, fixed and independent of $p$.
Using the expansion
$(H_{q}^p)^{-1}(z)=e^{p/2}q+z+O(e^{-p/2}|z|)$,
we find
$$
\aligned
\xi^0_i=q+\frac{1}{\sqrt{p}}\widehat{\xi}_i^0
+O\left(\frac{e^{-p/2}}{\sqrt{p}}|\widehat{\xi}_i^0|\right),
\quad i=l+1,\ldots,m.
\endaligned
$$
Clearly, $\xi^0=(\xi^0_1,\ldots,\xi^0_m)\in\mathcal{O}_{p}(q)$
because of $q\in\po$ and $\kappa>1$.
Notice that $q\in\po$ is a strict local maximum point of $a(x)$ over $\oo$
and satisfies $\partial_{\nu}a(q)=\langle\nabla a(q),\,\nu(q)\rangle=0$. Then we can derive that
there is
a constant $C>0$ independent of $p$ such that
$$
\aligned
a(q)-\frac{C}{p}\leq a(\xi^0_i)<a(q),
\qquad i=1,\ldots,m.
\endaligned
$$
From definition (\ref{1.3}),
Lemmas A.2 and A.3 we conclude  that
for any $i=1,\ldots,l$ and $k=1,\ldots,m$ with $i\neq k$,
$$
\aligned
H(\xi^0_i,\xi^0_i)=\frac1{4\pi}\log p+O\left(1\right),
\qquad\quad
G(q,\xi^0_i)=\frac1{2\pi}\log p+O\left(1\right),
\qquad\quad
G(\xi^0_k,\xi^0_i)=\frac1{2\pi}\log p+O\left(1\right).
\endaligned
$$
By  (\ref{7.6})  we find
\begin{eqnarray}\label{7.7}
\max\limits_{\xi\in\overline{\mathcal{O}}_p(q)}
F_p(\xi)\geq
F_p(\xi^0)\geq
\frac{2\pi ea(q)}{p^2}
\big\{
 (m+l+2+2\alpha)p-[2(m+l)^2+(8+6\alpha)(m+l)+4+4\alpha]\log p
+O(1)
\big\}.
\end{eqnarray}
Next, we suppose $\xi^p=(\xi^p_1,\ldots,\xi^p_m)\in\partial\mathcal{O}_{p}(q)$. Then
there exist five cases:\\
C1. \,\,There exists an $i_0\in\{1,\ldots,l\}$ such that
$\xi^p_{i_0}\in\partial B_d(q)\cap\Omega$, in which case,
$a(\xi^p_{i_0})<a(q)-d_0$ for some $d_0>0$ \\
\indent\indent independent of $p$;\\
C2. \,\,There exists an $i_0\in\{l+1,\ldots,m\}$ such that
$\xi^p_{i_0}\in\partial B_d(q)\cap\partial\Omega$, in which case,
$a(\xi^p_{i_0})<a(q)-d_0$ for some  \\
\indent\indent  $d_0>0$ independent of $p$;\\
C3. \,\,There exists an $i_0\in\{1,\ldots,l\}$ such that
$\dist(\xi_{i_0}^p,\po)=p^{-\kappa}$;\\
C4. \,\,There exists an $i_0\in\{1,\ldots,m\}$ such that
$|\xi_{i_0}^p-q|=p^{-\kappa}$;\\
C5. \,\,There exist some indices $i_0$, $k_0$, $i_0\neq k_0$ such that
$|\xi_{i_0}^p-\xi_{k_0}^p|=p^{-\kappa}$.\\
From $(A3)$-$(A6)$, (\ref{1.2})  and the maximum principle we obtain that for all $i=1,\ldots,l$ and $k=1,\ldots,m$ with $i\neq k$,
\begin{equation}\label{7.8}
\aligned
G(q,\xi_i^p)>0,
\qquad\,\,
G(\xi^p_k,\xi^p_i)>0,
\qquad\,\,
H(\xi^p_k,\xi^p_i)>0,
\qquad\,\,
H(\xi^p_i,\xi^p_i)=-\frac1{2\pi}\log\big[\dist(\xi_i^p,\po)\big]+O\left(1\right).
\endaligned
\end{equation}
For the first and second cases, by (\ref{7.6}) and (\ref{7.8}) we have
\begin{eqnarray}\label{7.9}
\max\limits_{\xi\in\overline{\mathcal{O}}_p(q)}
F_p(\xi)\leq
\frac{2\pi e}{p^2}
\left\{p\big[(m+l+2+2\alpha)a(q)-d_0\big]+O\left(\log p\right)\right\},
&&
\end{eqnarray}
which contradicts to (\ref{7.7}).
This shows that $a(\xi_i^p)\rightarrow a(q)$. By the condition over $a$,
we get $\xi_i^p\rightarrow q$ for all $i=1,\ldots,m$.\\
For the third case, by   (\ref{7.6}) and (\ref{7.8}) we have that
 if $0<\alpha\not\in\mathbb{N}^*$,
 \begin{eqnarray}\label{7.10}
\max\limits_{\xi\in\overline{\mathcal{O}}_p(q)}
F_p(\xi)\leq
\frac{2\pi e}{p^2}
\left\{2\sum_{i=1}^la(\xi_i^p)[
p-2\log p
-4\alpha\log|\xi_i^p-q|
]+\sum_{i=l+1}^ma(\xi_i^p)[
p-2\log p-4\alpha\log|\xi_i^p-q|
]+O\big(1\big)
\right.
&&\nonumber\\
\left.
+2(1+\alpha)
a(q)\left[
p-2\log p
-16\pi \sum_{i=1}^l
G(q,\xi_i^p)
+ 8\sum_{i=l+1}^m
\log|\xi_i^p-q|
\right]
-16\pi a(\xi^p_{i_0})H(\xi^p_{i_0},\xi^p_{i_0})
\right\}
\quad\,\,\,
&&\nonumber\\
\leq
\frac{2\pi e}{p^2}
\left\{2\sum_{i=1}^la(q)[
p-2\log p
-4\alpha\log|\xi_i^p-q|
]+\sum_{i=l+1}^ma(q)[
p-2\log p-4\alpha\log|\xi_i^p-q|
]+O\big(1\big)
\right.
\quad
&&\nonumber\\
\left.
+2(1+\alpha)
a(q)\left[
p-2\log p
-16\pi \sum_{i=1}^l
G(q,\xi_i^p)
+ 8\sum_{i=l+1}^m
\log|\xi_i^p-q|
\right]
-8\kappa a(\xi^p_{i_0})\log p
\right\}
\qquad\qquad\,\,
&&\nonumber\\
\leq
\frac{2\pi e}{p^2}
\left\{
(m+l+2+2\alpha)(p-2\log p)a(q)
-8\kappa a(\xi^p_{i_0})\log p
+O(1)
\right\},
\qquad\qquad\qquad\quad
\qquad\qquad\qquad
&&
\end{eqnarray}
where the last inequality is due to the fact that
for any $i=1,\ldots,l$, by (A3) and (\ref{1.3}),
$$
\aligned
-4\alpha\log|\xi_i^p-q|
-16\pi(1+\alpha)
G(q,\xi_i^p)
=(8+4\alpha)\log|\xi_i^p-q|-16\pi(1+\alpha)
H(q,\xi_i^p)
\leq C,
\endaligned
$$
and for any $i=l+1,\ldots,m$,
$$
\aligned
-4\alpha\log|\xi_i^p-q|
+16(1+\alpha)
\log|\xi_i^p-q|
=(16+12\alpha)\log|\xi_i^p-q|
\leq C
\endaligned
$$
with some large constant $C>0$.
While if $-1<\alpha<0$,
\begin{eqnarray}\label{7.11}
\max\limits_{\xi\in\overline{\mathcal{O}}_p(q)}
F_p(\xi)
\leq
\frac{2\pi e}{p^2}
\left\{
(m+l+2+2\alpha)(p-2\log p)a(q)
-16\pi a(\xi^p_{i_0})H(\xi^p_{i_0},\xi^p_{i_0})
+O(1)
\right\}
&&\nonumber\\
\leq
\frac{2\pi e}{p^2}
\left\{
(m+l+2+2\alpha)(p-2\log p)a(q)
-8\kappa a(\xi^p_{i_0})\log p
+O(1)
\right\}.
\qquad\,\,\,
&&
\end{eqnarray}
For the fourth  case, by (\ref{7.6}) and (\ref{7.8}) we have that if
$i_0\in\{1,\ldots,l\}$,
\begin{eqnarray}\label{7.12}
\max\limits_{\xi\in\overline{\mathcal{O}}_p(q)}
F_p(\xi)
\leq
\frac{2\pi e}{p^2}
\left\{
(m+l+2+2\alpha)(p-2\log p)a(q)
-16\pi a(\xi^p_{i_0})H(\xi^p_{i_0},\xi^p_{i_0})
+O(1)
\right\}
&&\nonumber\\
\leq
\frac{2\pi e}{p^2}
\left\{
(m+l+2+2\alpha)(p-2\log p)a(q)
-8\kappa a(\xi^p_{i_0})\log p
+O(1)
\right\}
\qquad\,\,\,\,\,
&&
\end{eqnarray}
because
$q\in\po$ and $p^{-\kappa}\leq\dist(\xi_{i_0}^p,\po)\leq|\xi_{i_0}^p-q|=p^{-\kappa}$, while if
$i_0\in\{l+1,\ldots,m\}$,
\begin{eqnarray}\label{7.13}
\max\limits_{\xi\in\overline{\mathcal{O}}_p(q)}
F_p(\xi)
\leq
\frac{2\pi e}{p^2}
\left\{
(m+l+2+2\alpha)(p-2\log p)a(q)
+\big[16(1+\alpha) a(q)
-4\alpha a(\xi^p_{i_0})\big]\log|\xi^p_{i_0}-q|
+O(1)
\right\}
&&\nonumber\\
\leq
\frac{2\pi e}{p^2}
\left\{
(m+l+2+2\alpha)(p-2\log p)a(q)
-4\kappa a(\xi^p_{i_0})\log p
+O(1)
\right\}.
\qquad\qquad\qquad\qquad\quad\quad\,\,\,
&&
\end{eqnarray}
For the last case, by (\ref{7.6}) and (\ref{7.8}) we have that if
$i_0\in\{1,\ldots,m\}$ and $k_0\in\{1,\ldots,l\}$,
\begin{eqnarray}\label{7.14}
\max\limits_{\xi\in\overline{\mathcal{O}}_p(q)}
F_p(\xi)
\leq
\frac{2\pi e}{p^2}
\left\{
(m+l+2+2\alpha)(p-2\log p)a(q)
+8 a(\xi^p_{i_0})
\log|\xi^p_{i_0}-\xi^p_{k_0}|
+O(1)
\right\}
&&\nonumber\\
\leq
\frac{2\pi e}{p^2}
\left\{
(m+l+2+2\alpha)(p-2\log p)a(q)
-8\kappa a(\xi^p_{i_0})\log p
+O(1)
\right\},
\qquad\,\,\,\,
&&
\end{eqnarray}
while if
$i_0\in\{l+1,\ldots,m\}$ and $k_0\in\{l+1,\ldots,m\}$,
\begin{eqnarray}\label{7.15}
\max\limits_{\xi\in\overline{\mathcal{O}}_p(q)}
F_p(\xi)
\leq
\frac{2\pi e}{p^2}
\left\{
(m+l+2+2\alpha)(p-2\log p)a(q)
+4 a(\xi^p_{i_0})
\log|\xi^p_{i_0}-\xi^p_{k_0}|
+O(1)
\right\}
&&\nonumber\\
\leq
\frac{2\pi e}{p^2}
\left\{
(m+l+2+2\alpha)(p-2\log p)a(q)
-4\kappa a(\xi^p_{i_0})\log p
+O(1)
\right\}.
\qquad\,\,\,\,\,
&&
\end{eqnarray}
Comparing (\ref{7.10})-(\ref{7.15}) with  (\ref{7.7}), we obtain
\begin{equation}\label{7.14}
\aligned
2(m+l+2+2\alpha)a(q)
\log p+8\kappa a(\xi^p_{i_0})
\log p \leq
[2(m+l)^2+(8+6\alpha)(m+l)+4+4\alpha]a(q)
\log p
+O(1),
\endaligned
\end{equation}
which  is impossible by the choice of $\kappa$ in (\ref{2.4}).
\,\qquad\qquad\qquad\qquad\qquad\qquad\qquad\qquad
\qquad\qquad\qquad\qquad\qquad\quad\,\,\,$\square$

\vspace{1mm}

\section*{Acknowledgments}
This research is supported by
the National
Natural Science Foundation of China  under Grant
No.   11671354.

\vspace{1mm}

\section*{Appendix A}
In this appendix we list  some properties of
the regular part of Green's function and its corresponding
Robin's function, see \cite{AP} for proofs.

Let the vector function $T(x)=(T_1(x),T_2(x))$ be the solution of
\begin{eqnarray*}
\qquad\qquad\qquad\qquad
\qquad\qquad\qquad
\Delta_xT-T=\frac{x}{|x|^2}\,\,\quad\,\,\textrm{in}
\,\,\,\,\mathbb{R}^2,\,\,\quad\,\,T(x)\in L^{\infty}_{loc}(\mathbb{R}^2).
\,\qquad\qquad\qquad\qquad\quad\qquad
(A1)&&
\end{eqnarray*}
Standard elliptic regularity theory yields that
$T(x)\in W^{2,\sigma}_{loc}(\mathbb{R}^2)\cap C^{\infty}(\mathbb{R}^2\setminus\{0\})$
for any $1<\sigma<2$,
and further the Sobolev embeddings  imply  that
$T(x)\in W^{1,1/\beta}(B_r(0))\cap C^{\beta}(\overline{B_r(0)})$
for any $r>0$ and $0<\beta<1$.

\vspace{1mm}
\vspace{1mm}
\vspace{1mm}
\vspace{1mm}

\noindent{\bf Lemma A.1.}\,\,{\it
Let $T(x)$ be the function described in {\upshape (A1)}.
There exists a function $H_1(x,y)$ such that\\
\indent {\upshape (i)}  for every $x,y\in\oo$,
\begin{eqnarray*}
\,\qquad\qquad\qquad\qquad\quad
\qquad\,
H(x,y)=H_1(x,y)+\left\{
\aligned
&\frac1{2\pi}\nabla\log a(y)\cdot T(x-y),
\,\quad\,\,y\in\Omega,\\[1mm]
&\,\frac1{\pi}\,\nabla\log a(y)\cdot T(x-y),
\,\,\quad\,y\in\po,
\endaligned\right.
\,\qquad\qquad\qquad\qquad
(A2)&&
\end{eqnarray*}
\indent {\upshape (ii)} the mapping
$y\in\oo\mapsto H_1(\cdot,y)$ belongs to
$ C^1\big(\Omega, C^{1}(\overline{\Omega})\big)
\cap C^1\big(\partial\Omega, C^{1}(\overline{\Omega})\big)$.\\
In this way, $y\in\oo\mapsto H(\cdot,y)\in C\big(\Omega, C^{\beta}(\overline{\Omega})\big)
\cap C\big(\partial\Omega, C^{\beta}(\overline{\Omega})\big)$ and
$H(x,y)\in C^\beta\big(\oo\times\Omega\big)\cap C^\beta\big(\oo\times\po\big)
\cap C^1\big(\oo\times\Omega\setminus\{x=y\}\big)
\cap C^1\big(\oo\times\po\setminus\{x=y\}\big)$
for any $\beta\in(0,1)$,
and the corresponding Robin's function $y\in\oo\mapsto H(y,y)$ belongs to $C^1(\Omega)\cap C^1(\partial\Omega)$.
}

\vspace{1mm}
\vspace{1mm}
\vspace{1mm}
\vspace{1mm}

Let
$\Omega_d:=\left\{
y\in\Omega\big|\,\,
\dist(y,\partial\Omega)<d\,
\right\}$
with   $d>0$ sufficiently small but fixed. Then for any
$y\in\Omega_d$, there exists a unique reflection
of $y$ across $\po$ along the outer normal direction, $y^*\in\Omega^c$,
such  that  $|y-y^*|=2\dist(y,\partial\Omega)$.

\vspace{1mm}
\vspace{1mm}
\vspace{1mm}
\vspace{1mm}

\noindent{\bf Lemma A.2.}\,\,{\it
There exists a mapping
$y\in\Omega_d\mapsto \mathrm{z}(\cdot,y)\in C\big(\Omega_d, C^{\beta}(\overline{\Omega})\big)
\cap L^\infty\big(\Omega_d, C^{\beta}(\overline{\Omega})\big)$
for any $\beta\in(0,1)$ such that for any $x\in\oo$ and  $y\in\Omega_d$,
\begin{eqnarray*}
\qquad\qquad\qquad\qquad\qquad\qquad\qquad\qquad\,
H(x,y)=\frac1{2\pi}\log\frac{1}{|x-y^*|}+\mathrm{z}(x,y),
\,\,\qquad\qquad\qquad\qquad\qquad\qquad\quad\,\,\,
(A3)&&
\end{eqnarray*}
and
\begin{eqnarray*}
\qquad\qquad\qquad\qquad\,\,\,\,
\mathrm{z}(x,y)=\frac1{2\pi}\nabla\log a(y)\cdot T(x-y)-\frac1{2\pi}\nabla\log a(y^*)\cdot T(x-y^*)+\tilde{\mathrm{z}}(x,y),
\,\,\qquad\qquad\quad\,\,\,\,\,
(A4)&&
\end{eqnarray*}
where the mapping
$y\in\overline{\Omega_d}\mapsto \tilde{\mathrm{z}}(\cdot,y)$ belongs to
$ C^1\big(\overline{\Omega_d},\,C^{1}(\overline{\Omega})\big)$.
}

\vspace{1mm}
\vspace{1mm}
\vspace{1mm}
\vspace{1mm}

\noindent{\bf Lemma A.3.}\,\,{\it
The Robin's function
$y\in\Omega\mapsto H(y,y)$ satisfies
\begin{eqnarray*}
\qquad\qquad\qquad\qquad\qquad\qquad\qquad
H(y,y)=\frac1{2\pi}\log\frac{1}{|y-y^*|}+\mathrm{z}(y),
\,\,\quad\forall\,\,y\in\Omega_d,
\,\,\qquad\qquad\qquad\qquad\quad\quad\,\,\,\,\,\,\,
(A5)&&
\end{eqnarray*}
where $\mathrm{z}\in C^1\big(\overline{\Omega_d}\big)$ and
\begin{eqnarray*}
\qquad\quad\qquad\qquad\,\,
\mathrm{z}(y)=\frac1{2\pi}\nabla\log a(y)\cdot T(0)-\frac1{2\pi}\nabla\log a(y^*)\cdot T(y-y^*)+\tilde{\mathrm{z}}(y,y),
\,\,\quad\forall\,\,y\in\Omega_d.
\,\,\qquad\qquad\,\,\,\,
(A6)&&
\end{eqnarray*}
}

\section*{Appendix B}
\noindent{\bf Proof of Lemma 2.1.}\,\,
Observe first that, for any $0<\tau<1$,
by (\ref{2.1}), (\ref{2.6}) and (\ref{2.12}),
$$
\aligned
\left\{\aligned&
-\Delta_a
H_0+H_0=\frac{1}{\gamma\mu_0^{2/(p-1)}}\left\{
\left[-4(1+\alpha)+\frac{C_1}{p}+\frac{C_2}{p^2}
\right]\left[
\frac{|x-q|^{2\alpha}(x-q)\cdot\nabla\log a(x)}{\delta_0^{2(1+\alpha)}+\big|x-q\big|^{2(1+\alpha)}}
-\frac{\log\big(\delta_0^{2(1+\alpha)}+|x-q|^{2(1+\alpha)}\big)}{2(1+\alpha)}\right]
\right.\\[2mm]
&\left.\,\quad\qquad\qquad\qquad
-\log\big[8(1+\alpha)^2\delta_0^{2(1+\alpha)}\big]
+\left(\frac{C_1}{p}+\frac{C_2}{p^2}\right)\log\delta_0
\right.\\[1mm]
&\left.\,\quad\qquad\qquad\qquad
+\frac1p
O_{\large L^{\infty}\big(\Omega\setminus B_{\delta_0^{\tau/2}}(q)\big)}
\left(\frac{\delta_0^{1+\alpha}}{\delta_0^{1+\alpha}+\big|x-q\big|^{1+\alpha}}
+\frac{\delta_0^{1+\alpha}}{\delta_0^{2+\alpha}+\big|x-q\big|^{2+\alpha}}\right)
\right.\\[2mm]
&\left.\,\quad\qquad\qquad\qquad
+\frac1p
O_{\large L^{\infty}\big(\Omega\bigcap B_{\delta_0^{\tau/2}}(q)\big)}\left(
\frac{|x-q|^{2\alpha}\big|(x-q)\cdot\nabla\log a(x)\big|}{\delta_0^{2(1+\alpha)}+\big|x-q\big|^{2(1+\alpha)}}
+\log\frac{\delta_0^{2(1+\alpha)}+|x-q|^{2(1+\alpha)}}{\delta_0^{2(1+\alpha)}}
\right)
\right\}
\,\,\,
\textrm{in}\,\,\,\,\ \,\,\Omega,\\[2mm]
&\frac{\partial H_0}{\partial \nu}=-\frac{1}{\gamma\mu_0^{2/(p-1)}}\left\{
\left[-4(1+\alpha)+\frac{C_1}{p}+\frac{C_2}{p^2}
\right]
\frac{|x-q|^{2\alpha}(x-q)\cdot\nu(x)}{\delta_0^{2(1+\alpha)}+\big|x-q\big|^{2(1+\alpha)}}
+\frac1p
O_{\large L^{\infty}\big(\partial\Omega\setminus B_{\delta_0^{\tau/2}}(q)\big)}
\left(\frac{\delta_0^{1+\alpha}}{\delta_0^{2+\alpha}+\big|x-q\big|^{2+\alpha}}\right)
\right.\\[2mm]
&\left.\,
\quad\qquad
+\frac1p
O_{\large L^{\infty}\big(\partial\Omega\bigcap B_{\delta_0^{\tau/2}}(q)\big)}\left(
\frac{|x-q|^{2\alpha}\big|(x-q)\cdot\nu(x)\big|}{\delta_0^{2(1+\alpha)}+\big|x-q\big|^{2(1+\alpha)}}
\right)
\right\}
\,\qquad\qquad\qquad\qquad\qquad\qquad
\qquad\qquad\qquad\,\,\,\,
\textrm{on}\,\,\ \ \,\po.
\endaligned\right.
\endaligned
$$
Observing that $q\in\oo$, by (\ref{1.2}), (\ref{1.3})
and (\ref{2.20}) we find
that
the regular part of Green's function, $H(x,q)$,  satisfies
$$
\left\{\aligned
&-\Delta_a
H(x,q)+H(x,q)=
\frac{4(1+\alpha)}{c_0}\log|x-q|-
\frac{4(1+\alpha)}{c_0}\frac{\,(x-q)\cdot\nabla\log a(x)\,}{|x-q|^2}
\quad
\,\,\textrm{in}\,\,\,\,\,\Omega,\\[1mm]
&\frac{\partial H(x,q)}{\partial \nu}
=\frac{4(1+\alpha)}{c_0}\frac{\,(x-q)\cdot\nu(x)\,}{|x-q|^2}
\qquad\qquad\qquad\qquad\quad\,\,
\qquad\qquad\qquad\qquad\qquad
\textrm{on}\,\,\,\po.
\endaligned\right.
$$
Thus if we define
$$
\aligned
Z_0(x)=\gamma\mu_0^{2/(p-1)}H_0(x)
-\left(1-\frac{C_1}{4(1+\alpha)p}-\frac{C_2}{4(1+\alpha)p^2}
\right) c_0 H(x,q)
+\log\left(8(1+\alpha)^2\delta_0^{2(1+\alpha)}\right)-\left(\frac{C_1}{p}+\frac{C_2}{p^2}\right)\log\delta_0,
\endaligned
$$
then $Z_0(x)$ satisfies
$$
\aligned
\left\{\aligned&
-\Delta_a
Z_0+Z_0=
\left[-4(1+\alpha)+\frac{C_1}{p}+\frac{C_2}{p^2}
\right]\left[\,\frac{1}{2(1+\alpha)}\log
\frac{|x-q|^{2(1+\alpha)}}{\delta_0^{2(1+\alpha)}+|x-q|^{2(1+\alpha)}}
-\frac{\delta_0^{2(1+\alpha)}(x-q)\cdot\nabla\log a(x)}{|x-q|^2\big(\delta_0^{2(1+\alpha)}+|x-q|^{2(1+\alpha)}\big)}
\right]\\[2mm]
&\quad\qquad\qquad\qquad
+\frac1p
O_{\large L^{\infty}\big(\Omega\setminus B_{\delta_0^{\tau/2}}(q)\big)}
\left(\frac{\delta_0^{1+\alpha}}{\delta_0^{1+\alpha}+\big|x-q\big|^{1+\alpha}}
+\frac{\delta_0^{1+\alpha}}{\delta_0^{2+\alpha}+\big|x-q\big|^{2+\alpha}}\right)\\[2mm]
&\quad\qquad\qquad\qquad
+\frac1p
O_{\large L^{\infty}\big(\Omega\bigcap B_{\delta_0^{\tau/2}}(q)\big)}\left(
\frac{|x-q|^{2\alpha}\big|(x-q)\cdot\nabla\log a(x)\big|}{\delta_0^{2(1+\alpha)}+\big|x-q\big|^{2(1+\alpha)}}
+\log\frac{\delta_0^{2(1+\alpha)}+|x-q|^{2(1+\alpha)}}{\delta_0^{2(1+\alpha)}}
\right)
\quad\,
\textrm{in}\,\,\ \,\,\ \,\,\Omega,\\[2mm]
&\frac{\partial Z_0}{\partial \nu}=
\left[-4(1+\alpha)+\frac{C_1}{p}+\frac{C_2}{p^2}
\right]
\frac{(x-q)\cdot\nu(x)}{\big|x-q\big|^{2}}
\frac{\delta_0^{2(1+\alpha)}}{\delta_0^{2(1+\alpha)}+\big|x-q\big|^{2(1+\alpha)}}
+\frac1p
O_{\large L^{\infty}\big(\partial\Omega\setminus B_{\delta_0^{\tau/2}}(q)\big)}
\left(\frac{\delta_0^{1+\alpha}}{\delta_0^{2+\alpha}+\big|x-q\big|^{2+\alpha}}\right)
\\[2mm]
&\qquad\quad
+\frac1p
O_{\large L^{\infty}\big(\partial\Omega\bigcap B_{\delta_0^{\tau/2}}(q)\big)}\left(
\frac{|x-q|^{2\alpha}\big|(x-q)\cdot\nu(x)\big|}{\delta_0^{2(1+\alpha)}+\big|x-q\big|^{2(1+\alpha)}}
\right)
\,\qquad\qquad\qquad
\quad\qquad\qquad\quad\qquad
\quad\qquad\qquad\,\,\,\textrm{on}\,\ \,\,\,\po.
\endaligned\right.
\endaligned
$$
Applying the polar coordinates with center $q$, i.e. $r=|x-q|$, and using the change of  variables $s=r/\delta_0$, we
get that
for any  $\theta>1$,
$$
\aligned
\int_{\Omega}\left|
\log
\left(\frac{|x-q|^{2(1+\alpha)}}{\delta_0^{2(1+\alpha)}+|x-q|^{2(1+\alpha)}}\right)
\right|^\theta dx
&\leq C\int_{0}^{\diam(\Omega)}\left|\log
\left(\frac{r^{2(1+\alpha)}}{\delta_0^{2(1+\alpha)}+r^{2(1+\alpha)}}\right)
\right|^\theta rdr\\
&\leq C\delta_0^2\int_{0}^{\diam(\Omega)/\delta_0}\left|\log\left(
1+\frac1{s^{2(1+\alpha)}}
\right)\right|^\theta sds\\[1mm]
&\leq C\left[\delta_0^2+\delta_0^{2\theta(1+\alpha)}\right],
\endaligned
$$
and
$$
\aligned
\int_{\Omega\setminus B_{\delta_0^{\tau/2}}(q)}\left|
\frac{\delta_0^{1+\alpha}}{\delta_0^{1+\alpha}+\big|x-q\big|^{1+\alpha}}
+\frac{\delta_0^{1+\alpha}}{\delta_0^{2+\alpha}+\big|x-q\big|^{2+\alpha}}
\right|^\theta dx
&\leq C\int_{\delta_0^{\tau/2}}^{\diam(\Omega)}\left[\frac{\delta_0^{\theta(1+\alpha)}}{(\delta_0^{1+\alpha}+r^{1+\alpha})^\theta}
+\frac{\delta_0^{\theta (1+\alpha)}}{(\delta_0^{2+\alpha}+r^{2+\alpha})^\theta}
\right]rdr\\
&\leq C\int_{\delta_0^{\tau/2-1}}^{\diam(\Omega)/\delta_0}\left[\frac{\delta_0^2}{(1+s^{1+\alpha})^\theta}
+\frac{\delta_0^{2-\theta}}{(1+s^{2+\alpha})^\theta}
\right]sds\\[1.5mm]
&\leq C\left[\delta_0^{\theta(1+\alpha)}
+\delta_0^{\tau+\theta(1+\alpha)-\frac12\tau \theta(2+\alpha)}
\right],
\endaligned
$$
and for  any $1<\theta<2$,
$$
\aligned
\int_{\Omega\bigcap B_{\delta_0^{\tau/2}}(q)}&\left|
\frac{|x-q|^{2\alpha}\big|(x-q)\cdot\nabla\log a(x)\big|}{\delta_0^{2(1+\alpha)}+\big|x-q\big|^{2(1+\alpha)}}
+\log\frac{\delta_0^{2(1+\alpha)}+|x-q|^{2(1+\alpha)}}{\delta_0^{2(1+\alpha)}}
\right|^\theta dx\\
\leq&
C\int_{0}^{\delta_0^{\tau/2}}\left[\left(
\frac{r^{1+2\alpha}}{\delta_0^{2(1+\alpha)}+r^{2(1+\alpha)}}
\right)^\theta+
\left|\log
\left(\frac{\delta_0^{2(1+\alpha)}+r^{2(1+\alpha)}}{\delta_0^{2(1+\alpha)}}\right)
\right|^\theta\right]rdr\\
=&C\delta_0^2\int_{0}^{\delta_0^{\tau/2-1}}
\left[\delta_0^{-\theta}\left(\frac{s^{1+2\alpha}}{1+s^{2(1+\alpha)}}\right)^\theta
+\log^\theta\big(
1+s^{2(1+\alpha)}
\big)\right]sds\\[1mm]
\leq& C
\delta_0^{\tau(2-\theta)/2},
\endaligned
$$
and
$$
\aligned
\int_{\Omega}\left|
\frac{(x-q)\cdot\nabla\log a(x)}{|x-q|^2}\frac{\delta_0^{2(1+\alpha)}}{\delta_0^{2(1+\alpha)}+|x-q|^{2(1+\alpha)}}
\right|^\theta dx
&\leq C\int_{0}^{\diam(\Omega)}\left|\frac{\delta_0^{2(1+\alpha)}}{r\big(\delta_0^{2(1+\alpha)}+r^{2(1+\alpha)}\big)}
\right|^\theta rdr\\
&\leq C\delta_0^{2-\theta}\int_{0}^{\diam(\Omega)/\delta_0}
\frac{s^{1-\theta}}{\big(1+s^{2(1+\alpha)}\big)^\theta}
ds\\[1mm]
&\leq C\left(\delta_0^{2-\theta}
+\delta_0^{2\theta(1+\alpha)}
\right).
\endaligned
$$
Then for  any  $1<\theta<2$,
\begin{equation*}
\aligned
\qquad\qquad\qquad\quad\,\,\,\,\,
\qquad\qquad\qquad\quad\quad
\big\|-\Delta_a
Z_0+Z_0\big\|_{L^\theta(\Omega)}\leq C\delta_0^{\tau(1/\theta-1/2)}.
\,\qquad\qquad\qquad\qquad
\qquad\qquad\quad\,\,\,
(B1)
\endaligned
\end{equation*}
On the other hand, if $q\in\po$, then, by using
the fact that
$|(x-q)\cdot\nu(q)|\leq C|x-q|^2$ for any $x\in\po$ (see \cite{AP})
we have that for any $\theta>1$,
$$
\aligned
\int_{\partial\Omega}\left|
\frac{(x-q)\cdot\nu(x)}{\big|x-q\big|^{2}}
\frac{\delta_0^{2(1+\alpha)}}{\delta_0^{2(1+\alpha)}+\big|x-q\big|^{2(1+\alpha)}}
\right|^\theta dx
&\leq C\int_{0}^{|\po|}\frac{\delta_0^{2\theta(1+\alpha)}}{\big(\delta_0^{2(1+\alpha)}+r^{2(1+\alpha)}\big)^\theta}
dr\\
&= C\delta_0\int_{0}^{|\po|/\delta_0}
\frac{1}{\big(1+s^{2(1+\alpha)}\big)^\theta}
ds\leq C\left(\delta_0
+
\delta_0^{2\theta(1+\alpha)}
\right),
\endaligned
$$
and
$$
\aligned
\int_{\partial\Omega\bigcap B_{\delta_0^{\tau/2}}(q)}\left|
\frac{|x-q|^{2\alpha}\big|(x-q)\cdot\nu(x)\big|}{\delta_0^{2(1+\alpha)}+\big|x-q\big|^{2(1+\alpha)}}
\right|^\theta dx
&\leq
C\int_{\partial\Omega\bigcap B_{\delta_0^{\tau/2}}(q)}
\left(
\frac{|x-q|^{2(1+\alpha)}}{\delta_0^{2(1+\alpha)}+\big|x-q\big|^{2(1+\alpha)}}
\right)^\theta
dx\\[1mm]
&\leq
C\left|\partial\Omega\cap B_{\delta_0^{\tau/2}}(0)
\right|\leq C\delta_0^{\tau/2},
\endaligned
$$
and
$$
\aligned
\int_{\partial\Omega\setminus B_{\delta_0^{\tau/2}}(q)}\left|
\frac{\delta_0^{1+\alpha}}{\delta_0^{2+\alpha}+\big|x-q\big|^{2+\alpha}}
\right|^\theta dx
&\leq C\int_{\delta_0^{\tau/2}}^{|\po|}\left(\frac{\delta_0^{1+\alpha}}{\delta_0^{2+\alpha}+r^{2+\alpha}}
\right)^\theta dr=
C\int_{\delta_0^{\tau/2-1}}^{|\po|/\delta_0}
\frac{\delta_0^{1-\theta}}{\,\big(1+s^{2+\alpha}\big)^\theta\,}
ds\\
&\leq
C\int_{\delta_0^{\tau/2-1}}^{|\po|/\delta_0}
\frac{\delta_0^{1-\theta}}{\,s^{\theta(2+\alpha)}\,}
ds
\leq C\delta_0^{\theta(1+\alpha)+\frac12\tau-\frac12\tau \theta(2+\alpha)},
\endaligned
$$
which implies that for  any $q\in\partial\Omega$, $\theta>1$ and $0<\tau<\min\{1,2(1+\alpha)/(2+\alpha)\}$,
\begin{equation*}
\aligned
\,\qquad\qquad\qquad\qquad\qquad
\qquad\qquad\qquad\qquad\,\,\,
\left\|\frac{\partial Z_0}{\partial \nu}\right\|_{L^{\theta}(\partial\Omega)}\leq C\delta_0^{\tau/2\theta}.
\,\qquad\qquad\qquad\qquad\qquad
\quad\qquad\qquad\quad\,\,\,
(B2)
\endaligned
\end{equation*}
While if $q\in\Omega$, then we estimate  that for any $x\in\partial\Omega$,
$$
\aligned
\left|
\frac{(x-q)\cdot\nu(x)}{\big|x-q\big|^{2}}
\frac{\delta_0^{2(1+\alpha)}}{\delta_0^{2(1+\alpha)}+\big|x-q\big|^{2(1+\alpha)}}
\right|\leq \frac{\delta_0^{2(1+\alpha)}}{\,|x-q|^{3+2\alpha}\,}\leq C\delta_0^{2(1+\alpha)},
\endaligned
$$
$$
\aligned
\frac{\delta_0^{1+\alpha}}{\delta_0^{2+\alpha}+\big|x-q\big|^{2+\alpha}}
\leq \frac{\delta_0^{1+\alpha}}{\,\big|x-q\big|^{2+\alpha}\,}
\leq C\delta_0^{1+\alpha},
\endaligned
$$
and then,
\begin{equation*}
\aligned
\,\quad\qquad\qquad\qquad\qquad\qquad
\qquad\qquad\qquad\quad\,\,\,
\left\|\frac{\partial Z_0}{\partial \nu}\right\|_{L^{\infty}(\partial\Omega)}
\leq C\delta_0^{1+\alpha}.
\,\qquad\qquad\qquad\qquad\qquad
\quad\qquad\qquad\quad\,\,\,
(B3)
\endaligned
\end{equation*}
As a consequence, by $(B1)$-$(B3)$ and elliptic regularity theory we conclude that for any $1<\theta<2$, $0<\tau<\min\{1,2(1+\alpha)/(2+\alpha)\}$ and $0<\lambda<1/\theta$,
$$
\aligned
\left\|Z_0\right\|_{W^{1+\lambda,\theta}(\Omega)}\leq
C\left(\big\|-\Delta_a
Z_0+Z_0\big\|_{L^\theta(\Omega)}+\left\|\frac{\partial Z_0}{\partial \nu}\right\|_{L^{\theta}(\partial\Omega)}\right)
\leq C\delta_0^{\tau(1/\theta-1/2)}.
\endaligned
$$
Furthermore, by  Morrey's embedding theorem we have that for any $0<\sigma<1/2+1/\theta$,
$$
\aligned
\left\|Z_0\right\|_{C^\sigma(\overline{\Omega})}\leq C\delta_0^{\tau(1/\theta-1/2)},
\endaligned
$$
which implies that expansion (\ref{2.21}) holds with $\beta=2\tau(1/\theta-1/2)$.
In addition, expansion (\ref{2.22}) can be also  derived from these analogous arguments of (\ref{2.21}).
\qquad\qquad\qquad\qquad\qquad\qquad\qquad\qquad
\qquad\qquad\qquad\qquad\qquad\qquad\qquad$\square$

\vspace{1mm}
\vspace{1mm}
\vspace{1mm}
\vspace{1mm}

\noindent{\bf Proof of Lemma 2.2.}\,\,Making the change of variables  $s=1/p$,
we  have that relations
(\ref{2.28}) and (\ref{2.30}) are equivalent to
the homogeneous equations
$\big(S_0(s,\xi,\mu),\ldots,S_m(s,\xi,\mu)\big)=\big(0,\ldots,0\big)$, namely
$$
\aligned
S_0(s,\xi,\mu):=&
\log\mu_0-\frac14c_0H(q,q)
+\frac{\,\log8(1+\alpha)^2+\frac{1}{4(1+\alpha)}C_1+\frac{1}{4(1+\alpha)}C_2s\,}{\,4-\frac{1}{1+\alpha}C_1s-
\frac{1}{1+\alpha}C_2s^2\,}
\\
&+\frac{1-\frac{1}{4}\widetilde{C}_1s-\frac{1}{4}\widetilde{C}_2s^2}{-4+\frac{1}{1+\alpha}C_1s+\frac1{1+\alpha}C_2s^2}\sum_{k=1}^m
\left(
\frac{\mu_0}{\mu_k|\xi_k-q|^{\alpha}}
\right)^{2s/(1-s)}
c_k G(q,\xi_k)
=0,
\endaligned
$$
and for each $i=1,\ldots,m$,
$$
\aligned
S_i(s,\xi,\mu):=&
\log\mu_i-\frac14c_iH(\xi_i,\xi_i)
+\frac{\,\log8+\frac{1}{4}\widetilde{C}_1+\frac{1}{4}\widetilde{C}_2s\,}{\,4-\widetilde{C}_1s-\widetilde{C}_2s^2\,}
-\frac14\sum_{k=1,\,k\neq i}^m
\left(
\frac{\,\mu_i|\xi_i-q|^{\alpha}\,}{\mu_k|\xi_k-q|^{\alpha}}
\right)^{2s/(1-s)}
c_k G(\xi_i,\xi_k)
\\
&+\frac{\,1-\frac{1}{4(1+\alpha)}C_1s-\frac{1}{4(1+\alpha)}C_2s^2\,}{\,-4+\widetilde{C}_1s+\widetilde{C}_2s^2\,}
\left(
\frac{\,\mu_i|\xi_i-q|^{\alpha}\,}{\mu_0}
\right)^{2s/(1-s)}
c_0 G(\xi_i,q)
=0.
\endaligned
$$
Obviously,  from the definitions of the constants $C_1$ and $\widetilde{C}_1$ in (\ref{2.15}) we  deduce that for $s=0$,
\begin{equation*}
\aligned
\qquad\qquad\qquad\qquad
\qquad\qquad\qquad\quad
\mu_0(0,\xi)=e\large^{-\frac{3}{4}+\frac14c_0H(q,q)
+\frac14\sum_{k=1}^m
c_k G(q,\xi_k)},
\qquad\qquad\qquad
\qquad\qquad\quad\,\,\,\,
(B4)
\endaligned
\end{equation*}
and for each $i=1,\ldots,m$,
\begin{equation*}
\aligned
\qquad\qquad\qquad\qquad
\qquad\qquad
\mu_i(0,\xi)=e\large^{-\frac{3}{4}+\frac14c_iH(\xi_i,\xi_i)
+\frac14
c_0G(\xi_i,q)
+\frac14\sum_{k=1,\,k\neq i}^m
c_k G(\xi_i,\xi_k)}.
\qquad\qquad\quad
\qquad\quad\,\,\,\,
(B5)
\endaligned
\end{equation*}
On the other hand,  observe that for any
$s>0$ small enough,
$$
\aligned
\left(
\frac{\,C^2\,}{s^{C+\kappa\alpha}}
\right)^{2s/(1-s)}
=e\large^{\frac{2s}{1-s}\big[2\log C+(C+\kappa\alpha)\log\frac1s\big]}=1+O\left(
s\log\frac{1}{s}
\right).
\endaligned
$$
Then by (\ref{2.3}) and the first estimate of (\ref{2.31})
we find that for any $i,k=1,\ldots,m$ and $i\neq k$,
\begin{equation*}
\aligned
\quad\,\,\,
\left(
\frac{\mu_0}{\mu_k|\xi_k-q|^{\alpha}}
\right)^{2s/(1-s)}
=1+O\left(
s\log\frac{1}{s}
\right)
\qquad
\textrm{and}
\qquad
\left(
\frac{\,\mu_i|\xi_i-q|^{\alpha}\,}{\mu_k|\xi_k-q|^{\alpha}}
\right)^{2s/(1-s)}
=1+O\left(
s\log\frac{1}{s}
\right).
\quad
(B6)
\endaligned
\end{equation*}
Moreover, by  $(A3)$, (\ref{1.3}), (\ref{2.3}) and the fact that
$a(\xi_i)G(\xi_i,\xi_k)=a(\xi_k)G(\xi_k,\xi_i)$
and
$a(\xi_i)G(\xi_i,q)=a(q)G(q,\xi_i)$
for all
$i,k=1,\ldots,m$ with $i\neq k$,
we can conclude that for any points $\xi=(\xi_1,\ldots,\xi_m)\in\mathcal{O}_{1/s}(q)$,
\begin{equation*}
\aligned
\qquad\qquad\,\,\,\,\,\,
G(\xi_i,\xi_k)=O\left(\log\frac1s\right)
\,\qquad\,
\textrm{and}
\,\qquad\
G(q,\xi_i)=O\left(\log\frac1s\right),
\,\quad\,\forall\,\,i,k=1,\ldots,m,\,\,i\neq k.
\qquad\quad
(B7)
\endaligned
\end{equation*}
Furthermore, some direct computations easily deduce that for each $i,k=1,\ldots,m$ with $i\neq k$,
$$
\aligned
\frac{\partial S_0(s,\xi,\mu)}{\partial\mu_0}=\frac1{\mu_0}\left[
1+O\left(s\log\frac1s\right)
\right],
\qquad
\qquad
\frac{\partial S_0(s,\xi,\mu)}{\partial\mu_k}=\frac1{\mu_k}O\left(s\log\frac1s\right),
\endaligned
$$
and
$$
\aligned
\frac{\partial S_i(s,\xi,\mu)}{\partial\mu_i}=\frac1{\mu_i}\left[
1+O\left(s\log\frac1s\right)
\right],
\quad\,\,
\frac{\partial S_i(s,\xi,\mu)}{\partial\mu_0}=\frac1{\mu_0}O\left(s\log\frac1s\right),
\quad\,\,
\frac{\partial S_i(s,\xi,\mu)}{\partial\mu_k}=\frac1{\mu_k}O\left(s\log\frac1s\right).
\endaligned
$$
If we set $S(s,\xi,\mu)=\big(S_0(s,\xi,\mu),\ldots,S_m(s,\xi,\mu)\big)$, then the vector
function $S(s,\xi,\mu)$
satisfies
$$
\aligned
\det\big(\nabla_{\mu}S(s,\xi,\mu)
\big)=\frac1{\,\mu_0\mu_1\ldots\mu_m\,}\left[
1+O\left(s\log\frac1s\right)
\right]\neq0.
\endaligned
$$
Consequently, using the implicit function theorem, we can derive  that the homogeneous equations
$S(s,\xi,\mu)=(0,\ldots,0)$ is solvable in some neighborhood of $\big(0,\xi,\mu(0,\xi)\big)$,
namely for   any points
$\xi=(\xi_1,\ldots,\xi_m)\in\mathcal{O}_{p}(q)$
and any $p>1$ large enough, systems (\ref{2.28}) and (\ref{2.30}) have a unique solution $\mu=(\mu_0,\,\mu_1,\ldots,\mu_m)$
satisfying the first estimate in
(\ref{2.31}). From $(B4)$-$(B7)$ it follows that
$$
\aligned
\mu_0=\mu_0(p,\xi)\equiv e\large^{-\frac{3}{4}+\frac14c_0H(q,q)
+\frac14\sum_{k=1}^m
c_k G(q,\xi_k)}\left[\,1+O\left(\frac{\log^2p}p\right)\right],
\endaligned
$$
and for each $i=1,\ldots,m$,
$$
\aligned
\mu_i=\mu_i(p,\xi)\equiv e\large^{-\frac{3}{4}+\frac14c_iH(\xi_i,\xi_i)
+\frac14
c_0G(\xi_i,q)
+\frac14\sum_{k=1,\,k\neq i}^m
c_k G(\xi_i,\xi_k)}\left[\,1+O\left(\frac{\log^2p}p\right)\right].
\endaligned
$$
Moreover,  by $(A3)$, $(A5)$, $(B6)$, $(B7)$, (\ref{2.3}),  (\ref{2.28}) and (\ref{2.30})
we find that the second estimate in (\ref{2.31}) holds.$\square$

\end{document}